\documentclass[11pt,leqno]{book}

  \usepackage{hyperref}
\usepackage{makeidx}     
\usepackage{amsmath}
\usepackage{amsfonts} 
\usepackage{amssymb}
 \makeindex

\newcommand{\EE}{\mathcal{E}}

\newcommand{\ThE}{\mbox{Th}(\EE)}
\newcommand{\LL}{\mathcal{L}}

\newcommand{\MM}{\mathcal{M}}

\newcommand{\A}{ \mathbf A}
\renewcommand{\AA}{ \mathbf A}

\newcommand{\B}{ \mathbf B}
\newcommand{\C}{ \mathbf C}
\newcommand{\D}{ \mathbf D}
\newcommand{\R}{ \mathbf R}
\newcommand{\bfR}{ \R}

\newcommand{\DD}{ \mathbf D}

\newcommand{\AAA}{ \mathcal{A}}
\newcommand{\CCC}{\mathcal{C}}
\newcommand{\DDD}{\mathcal{D}}

\renewcommand{\aa}{\mathbf{a}}

\newcommand{\ThA}{\mbox{Th}(\A)}

\newcommand{\NN}{{\mathbb{N}}}  
\newcommand{\bfN}{{\mathbf{N}}}  
\newcommand{\Nops}{(\NN, +, \times)}
\newcommand{\ZZ}{\mathbb{Z}}

\newcommand{\QQ}{\mathbb{Q}}

\newcommand{\ThN}{\mbox{Th}(\N)}

\renewcommand{\phi}{\varphi}

\renewcommand{\lim}{\mbox{lim}}
\renewcommand{\sup}{\mbox{sup}}

\newcommand{\id}[1]{\index{#1}}

\newcommand{\iid}[1]{{\it #1}\index{#1}} 
\newcommand{\ird}[1]{#1\index{#1}}    

\newcommand{\citelab}[2]{\cite{#1}\label{#1.#2}} 

\newcommand{\mmbox}[1]{\ \mbox{#1}\ }

\newcommand{\ttext}[1]{\ \text{#1}\ }

\newcommand{\ID}[1]{[ {#1}]_{id}}  

\newcommand{\IDEAL}{I{\small DEAL}\ }

\newcommand{\tria}{\triangleleft}
\newcommand{\trup}{\small \triangle}
\newcommand{\RR}{\mathcal{R}}
\newcommand{\BB}{\mathcal{B}}
\newcommand{\IB}{\mathcal{I}(\BB)}
\newcommand{\FF}{\mathcal{F}}

\newcommand{\DU}{(D;\le, \vee, 0,1)}
\newcommand{\sq}{\sqsubset}

\newcommand{\jump}[1]{\ES^{(#1)}}
\newcommand{\zer}{\mathbf 0}
\newcommand{\one}{\mathbf 1}
\renewcommand{\deg}{\mbox{deg}}
\newcommand{\len}{\mbox{length}}

\newcommand{\RRT}{\mathcal{R}_T}
\newcommand{\RT}{\RRT}
\newcommand{\RRm}{\mathcal{R}_m}
\newcommand{\Rm}{\RRm}
\newcommand{\Rme}{\Rm^-}
\newcommand{\RRW}{\mathcal{R}_{wtt}}
\newcommand{\Rwtt}{\mathcal{R}_{wtt}}
\newcommand{\Rtt}{\mathcal{R}_{tt}}

\newcommand{\DTE}{\DD_T(\le \ES')}
\newcommand{\Low}{\rm{Low}}
\newcommand{\High}{\rm{High}}
\newcommand{\SII}{\Sigma^0_2}
\newcommand{\SIII}{\Sigma^0_3}

\newcommand{\vsp}{\vspace{6pt}}
\newcommand{\vspsmall}{\vspace{3pt}}
\newcommand{\nii}{\not \in}

\newcommand{\ex}{\exists}
\newcommand{\fa}{\forall}
\newcommand{\LR}{\Leftrightarrow}
\newcommand{\LLR}{\ \Leftrightarrow \ }
\newcommand{\RA}{\Rightarrow}
\newcommand{\RRA}{\ \Rightarrow\ }
\newcommand{\LA}{\Leftarrow}

\newcommand{\n}{\noindent}
\newcommand{\wt}{\widetilde}
\newcommand{\eqq}{\equiv}
\newcommand{\beqq}{/_\eqq}
\newcommand{\sub}{\subseteq}
\newcommand{\ol}{\overline}

\renewcommand{\land}{\&}
\renewcommand{\lor}{\vee}
\newcommand{\ES}{\emptyset}

\renewcommand{\hat}{\widehat}
\renewcommand{\tilde}{\widetilde}

\newcommand{\lland}{\ \land \ }

\newtheorem{theorem}{Theorem}[section]
\newtheorem{definability lemma}[theorem]{Definability Lemma}

\newtheorem{fact}[theorem]{Fact}
\newtheorem{remark}[theorem]{Remark}
\newtheorem{proposition}[theorem]{Proposition}
\newtheorem{claim}[theorem]{Claim}
\newtheorem{notation}[theorem]{Notation}
\newtheorem{definition}[theorem]{Definition}
\newtheorem{lemma}[theorem]{Lemma}
\newtheorem{corollary}[theorem]{Corollary}
\newtheorem{example}[theorem]{Example}

\newtheorem{specialtrace}[theorem]{Trace Lemma}
\newtheorem{specialmainlemma}[theorem]{Main Lemma}
\newtheorem{transferlemma}[theorem]{Transfer Lemma}
\newtheorem{Idl}[theorem]{Ideal Definability Lemma}

 \newtheorem{specialmain}[theorem]{Main Theorem}

\newcommand{\la}{\langle}
\newcommand{\ra}{\rangle}

\newcommand{\pf}{\n {\it Proof.}\  \ }
\newcommand{\demo}[1]{{\it #1.} \ }
\renewcommand{\square}{\diamondsuit}  
\newcommand {\eop}{\hfill $\square$ \vspsmall \\ }
\newcommand {\eopnospace}{\hfill $\square$}

\newcommand{\be}{\begin{enumerate}}
\newcommand{\ee}{\end{enumerate}}
\newcommand{\bi}{\begin{itemize}}
\newcommand{\ei}{\end{itemize}}

\newcommand{\beq}{\begin{equation}}
\newcommand{\eeq}{\end{equation}}

\newcommand{\bc}{\begin{center}}
\newcommand{\ec}{\end{center}}

\renewcommand{\R}{{\mathcal{R}}}

\newcommand{\N}{\mathbb{N}}
\renewcommand{\NN}{\N}
\newcommand{\NNops}{(\N, +, \times)}

\newcommand{\Mp}{\mathbf{M}_{\ol \p}}

\newcommand{\M}{\mathbf{M}}

\renewcommand{\a}{\mathbf{a}}
\renewcommand{\b}{\mathbf{b}}
\renewcommand{\c}{\mathbf{c}}
\renewcommand{\d}{\mathbf{d}}
\newcommand{\e}{\mathbf{e}}

\newcommand{\g}{\mathbf{g}}
\newcommand{\h}{\mathbf{h}}

\newcommand{\p}{\mathbf{p}}
\newcommand{\q}{\mathbf{q}}

\renewcommand{\t}{\mathbf{t}}
\renewcommand{\u}{\mathbf{u}}
\renewcommand{\v}{\mathbf{v}}

\newcommand{\x}{\mathbf{x}}
\newcommand{\y}{\mathbf{y}}
\newcommand{\z}{\mathbf{z}}

\DeclareMathAlphabet{\mathbf}{OML}{cmm}{b}{it}

\newcommand{\Sk}{\Sigma^0_k}

\newcommand{\SSk}{\Sigma_k}
\newcommand{\PPk}{\Pi_k}

\newcommand{\Th}{\ensuremath{\mathrm{Th}}}

\newcommand{\PP}{\ensuremath{\mathcal{P}}}
\newcommand{\DTIME}{\ensuremath{\mathrm{D{\scriptstyle TIME}}}}
\newcommand{\PTIME}{\ensuremath{\mathcal{P}}}
\newcommand{\EXP}{\ensuremath{\mathrm{E{\scriptstyle XPTIME}}}}

\newcommand{\NP}{\ensuremath{\mathcal{NP}}}

\newcommand{\coNP}{\ensuremath{\mathrm{CoNP}}}
\newcommand{\Co}{\ensuremath{\mathrm{Co}}}
\newcommand{\Rec}{\hbox{\it Rec}}
\newcommand{\pT}{\le^p_T}

\newcommand{\leo}{\le^p_{1\text{-}tt}}
\newcommand{\ppr}{\le^p_r}

\newcommand{\ns}{\displaystyle}
\newcommand{\pr}{\preceq}
  
\begin{document} 

\title{Coding Methods in Computability \\ Theory and Complexity Theory}
 \author{Habilitationsschrift  \\ \vspace{1in}  \\ Andr\'e Nies} 

\date{January  1998}

\maketitle  

\frontmatter

\chapter{Preface}

A major part of computability theory focuses on the analysis of a few
 structures of central importance. As a tool, the method of coding with first-order
 formulas  has  been
 applied  with great success. It was  used to determine the complexity of the elementary theory, to provide restrictions on
 automorphisms, and even to obtain definability results. As an example, consider
 $\RT$, the structure of computably enumerable (c.e.) Turing degrees. The
 analysis by coding methods began with the proof by Harrington and
 Shelah \citelab{Harrington.Shelah:82}{1}  that $\Th(\RT)$ is
 undecidable. Extending the  coding methods used, Harrington and Slaman  (unpublished)
gave an interpretation  of $\Th\Nops$, also called \iid{true
  arithmetic},
 in $\Th(\RT)$. Here an \iid{interpretation} is  a
many-one-reduction of theories based on a computable map defined
on sentences in some natural way. A different approach
to the same problem, due to Slaman and Woodin, introduced a very
versatile way of coding copies of $\Nops$ into $\RT$ with parameters,
which  was a main ingredient  for the investigations in Nies, Shore and Slaman
\citelab{NSS:96}{pref}. In the latter work, the definability of some 
important classes, including   $\Low_2$ and $\High_1$, is proved. Moreover, it is
shown that no automorphism of $\RT$ can change the second jump of a
degree, and that a  coding of $\NN$ in $\RT$ without parameters
exists. In a different direction,  Lempp, Nies and Slaman \citelab{LNS:nd}{1},
combining the Harrington-Shelah type of coding with 
algebraic methods, proved that the $\fa\ex\fa$-theory of $\RT$ (as a
partial order) is undecidable.

We will describe how a similar  program can be  carried out
for several other structures, including $\Rm$, the structure of c.e.\ many
one degrees, $\RRW$,  the structure of c.e.\ weak truth table 
 degrees and $\EE$, the lattice of c.e.\ sets under inclusion. 
In all cases we will obtain undecidability of, or even an
interpretation  of $\ThN$ in the theory
of the structure.  For $\Rm$, we also obtain definability results and restrictions on
automorphisms. Moreover, for both $\Rm$ and $\RRW$ a  coding  of $\NN$
without parameters can be given. On the other hand, for $\EE$ such stronger
coding properties must fail: no infinite linear order can be
coded without parameters. In connection with the study of $\EE$, we
also consider the  lattices  $\IB$ of c.e.\  ideals for certain  c.e.\
boolean algebras $\BB$ and prove that their theories are undecidable. These
lattices, besides being of intrinsic interest in effective algebra,
can be coded  into
many important structures, like degree structures from complexity
theory ``low down''. Thereby they provide a tool to prove undecidability  for
theories from very different contexts.

While so far most structures from computability theory (and complexity
theory) were studied in isolation, our approach has a unifying aspect,
since first general tools and concepts of a model theoretic flavor are
developed, which then re-emerge again and again. For instance, for
$\Rm$, $\EE$ and to some extent the lattices $\IB$, we will prove
definability lemmas\id{definability lemma} which give a
way to pass from arithmetical definability  of subsets of a structure to definability with
parameters in the structure. These definability lemmas constitute the main tool for
our analyses by coding methods of the structures in question. 

\vspsmall

The first chapter and to some extend the second chapter are  of an
introductory nature. The methods  in Section \ref{sectionITA} and Section
\ref{fragments} appeared first
in  Nies \cite{Nies:95} and Nies \cite{Nies:96*1},
respectively. The first three sections from Chapter \ref{RRRmmm} are
also 
from \cite{Nies:95}. Chapter \ref{EEEEE} is based on Harrington and
Nies \citelab{Harrington.Nies:nd}{1}, but contains substantial improvements
in Section \ref{CodeRG} which lead to new results about fragments of
$\Th(\EE^*)$ in Section \ref{Efragments}.  Chapter   \ref{CHIB}
appeared in Nies \cite{Nies:97*1}, as did Section \ref{IntE}. Section
\ref{EIBC} is based on Downey and Nies \cite{Downey.Nies:97}, while
 Chapter \ref{ChapterRwtt} contains very recent results of the author.  An extended version
of this work containing also material about $\RT$ will appear as a
book in the {\it Oxford Logic Guides}.

\vspace{22pt}

\hfill{Heidelberg-Madison-Ithaca-Chicago, 1992-1997}

\tableofcontents

\mainmatter

\chapter{The objects of investigation}

We introduce the structures we will study and  discuss  some of their basic
properties.


\section{Structures based on computably enumerable sets} 

 A central notion in
logic is the notion of a computably enumerable (c.e., or  r.e.) set of natural numbers.
In this section we review structures based on c.e.\ sets.  The
study of global and local properties of these structures is regarded as a central
topic in computability theory. 

\subsection{Degree structures}

The  relative computational complexity of c.e.\ sets  is investigated through the study
of the uppersemilattices $\RR_m$ and $\RR_T$ of enumerable many-one
degrees\id{degree!many-one} and of enumerable Turing ($T$-)degrees\id{degree!Turing}, and also of the degree structures
$\RR_{wtt}$\id{degree!$wtt$}  and  $\RR_{tt}$\id{degree!$tt$} which arise from reducibilities
between $\le_m$ and $\le_T$.
These   reducibilities  are obtained
from Turing-reducibility by more and more restricting the underlying concept
of oracle computation: for subsets $X,Y$ of $\NN$, $X\le_{wtt}Y$ if 
 $X\le_TY$ via an oracle computation procedure where the largest  oracle
question asked is recursively bounded in the input;  $X\le_{tt}Y$ if 
such an oracle computation procedure is total for every oracle.
Finally, $X\le_mY$ if there is a computable function $f$ such that $n
\in X \LR f(n) \in Y$ for all $n$. To avoid trivialities, we actually
allow TRUE and FALSE as values of $f$. Since  each reducibility $\le_r$ is a
preordering, we obtain a
\ird{degree structure} $\RR_r$ of $r$-degrees of c.e.\ sets, which is an
\ird{upper semilattice}  with  a least element, denoted by $\zer$, and a
largest element, denoted by $\one$. The degree $\zer$ consists of the
computable sets, and $\one$ is the $r$-degree of the halting problem.
The $r$-degree\id{degree} of a set $X\sub \NN$ is denoted by $\deg_r(X)$.

Ever since Post's problem was formulated
\cite{Post:44}\label{Postpaper}  which asks
whether $\zer, \one$ are the only c.e.\ $T$-degrees, the study of
$\RT$ has been a mainstay of computability theory. A wide range of
facts about $\RT$,  all formalizable within first-order logic, were
found. The structure is dense (Sacks   \cite{Sacks:64}), has pairs with
infimum $\zer$ (called {\it minimal pairs}\id{minimal pair}; Yates \cite{Yates:66*}) but also
nonzero degrees which don't bound any minimal pairs (Lachlan; see
\cite{Soare:87}). Such properties seem to reflect pathological rather
than orderly behavior of $\RT$. The structure $\Rm$,  on the other hand, is
much more homogeneous and well-behaved, and in fact is the only c.e.\
degree structure which permits a characterization (Denisov \citelab{Denisov:72}{1},
see also Section \ref{CharRm}). While  $\Rtt$ exhibits quite 
a pathological behavior as well, $\RRW$ is at the borderline. For
instance,  the
same theorems about minimal pairs hold as for $\RT$, but it shares
with $\Rm$ the property of being {\it
  distributive}\id{distributive!upper semilattice}\id{upper semilattice!distributive} as an upper
semilattice, namely

 \begin{equation} \label{ldistr}\fa \x \fa \a \fa \b [\x \le \a \vee \b \ \RA \
\ex \a_0 \le \a \ex \b_0 \le \b \ \x = \a_0 \vee \b_0] \end{equation} 

(see  
Lachlan \citelab{Lachlan:72}{1} for a proof).
 
For the study of enumerable sets, the reducibilities refining $T$-reducibility
are interesting  partially   because they are more closely related to structural
properties of an enumerable set than $T$-reducibility is. For instance, a 
maximal enumerable set must have minimal many-one degree, and a hypersimple set is 
necessarily $wtt$-incomplete, but not always $T$-incomplete (see
Odifreddi \cite[p. 338]{Odifreddi:81}).

\subsection{C.e. sets under inclusion and ideal lattices}
\label{EIB}
 A more algebraic aspect of computably enumerable sets is
captured by the lattice $\EE$ of computably enumerable sets under
inclusion. This view of c.e.\ sets is the most elementary one,
because no further concepts are required to relate them. Clearly $\EE$
is a distributive lattice with least and greatest elements. Moreover, $\EE$ satisfies the reduction
principle: 

\begin{equation} \label{RedPr} \fa A \fa B \ \ex \wt A\sub A \ex \wt B
  \sub B [ \wt A \cap \wt B=\ES \lland \wt A \cup \wt B = A\cup B]. \end{equation}

 Despite of the 
conceptually simple way $\EE$ is  introduced,  it is a structure  of great algebraic complexity. 
Several interrelated directions in the study of $\EE$ have been pursued: one is the investigation of 
automorphisms (Soare \citelab{Soare:74}{1}), a further one is the relationship between  
the behavior of an enumerable set as an element of $\EE$ and its computational complexity 
(see e.g. Martin \citelab{Martin:66*1}{1}  and Harrington and  Soare \citelab{Harrington.Soare:91}{1}). 
Here we follow another approach, initiated by the undecidability proofs  for  $\ThE$
due to Herrmann \cite{Herrmann:84} and Harrington: the approach of
studying coding and definability.

The next  type of structures we consider is actually based not on
c.e.\ sets but on c.e.\ ideals. 
A boolean algebra $\BB$ is {\it computably enumerable}\id{boolean
  algebra!c.e.}  if 
$\BB= \DDD/_{\ns H}$ for a
c.e.\  ideal $H$ of the 
computable dense boolean algebra $\DDD$.
 Let $\IB$\id{ideals!lattice of} be the lattice of c.e.\ ideals of a c.e.\ 
boolean algebra $\BB$ (thus, if $\BB= \DDD/_{\ns H}$, $\IB$ is the lattice of c.e.\ ideals of 
$\DDD$ containing $H$). We list   some properties of $\IB$ which show that, in a sense, $\IB$
is similar to $\EE$. First,
 $\IB$ is a distributive lattice with least and greatest elements. It
 is easy to prove that 
  $\IB$ also satisfies  the reduction principle.  All
  principal\id{ideals!principal}  ideals
  $[0,b]_\BB$ of $\BB$ are in $\IB$. The class of  principal ideals is
  definable in $\IB$:  an ideal is principal iff it is
 complemented in $\IB$. 

It is possible that $\IB \cong \BB $, even    for a dense c.e.\
$\BB$:
 one can construct a dense $\BB$ 
 such that  every c.e.\
ideal is principal (Martin and Pour-El \citelab{Martin.Pour:70}{2}).
However, the type of c.e.\ boolean algebras we  consider here  have 
a very complex lattice of c.e.\ ideals. 
 We call   a c.e.\  boolean
algebra $\BB$ \iid{effectively dense} if,   for each
 element $x$ of $\BB$, we can effectively find an   element $y \le x$
 such that $x \neq 0$ implies $0 < y < x$.
Thus e.g.\ the recursive dense boolean algebra is effectively dense,
but in fact many other c.e.\ presentations of the countable dense boolean
algebra are as well.  For instance, consider the \ird{Lindenbaum
  algebra}  of
sentences over Peano arithmetic. This c.e.\ boolean algebra
is effectively dense by \ird{Rosser's theorem}, a refinement of
G\"odel's second incompleteness theorem (Example \ref{Rosser} below).

\section{Structures from Complexity Theory}

In complexity theory, one  considers  sets of strings, mostly from
 $\{0,1\}^{<\omega}$, instead of sets of numbers.  Polynomial
time bounded analogs of the recursion theoretic reducibilities were introduced. For instance, 
for  $X,Y \sub \{0,1\}^{<\omega}$, polynomial time many-one reducibility is defined by

$$X \le^p_m Y \LR (\ex f \in \PP [ X=f^{-1}(Y)]),$$ 

(where $f$, as
before, may have {\small TRUE}  and {\small FALSE} as
values). Polynomial time Turing  reducibility is defined by
$X \le^p_T Y \LR \ $ there is a   polynomial time bounded
deterministic oracle
Turing machine taking inputs in  $\{0,1\}^{< \omega}$
which computes $X$ if the oracle is $Y$.
Analogs of some  other reducibilities, like truth-table reducibility, can be
defined in a similar way. We let 
$(\Rec^p_r, \le^p_r)$\id{degree!polynomial time}  
be the p.o.\ of polynomial time  $r$--degrees 
of computable sets, where  $\le^p_r$ is a polynomial  time reducibility
in between  (and including) $\le^p_m$ and $\le^p_T$.
As before, $\Rec^p_r$ is an u.s.l.\  which  has a least element $\zer$, the
degree consisting of sets in $\PP$. But $\Rec^p_r$ has  no greatest element. 

The fact that the base sets are computable allows for a method
radically different from the methods used in computability theory: the
delay diagonalization (or looking-back) method introduced in
Landweber, Lipton and Robertson \citelab{Landweber:81}{1}. They used
the technique  to reprove 
 Ladner's result \citelab{Ladner:75}{1}  that   $\Rec^p_r$ is dense
 (see also Balcazar e.a.\ \cite{BDG:88}).
The idea is as follows: in a construction of a computable set $A$, at
stage $s$ $A^{=s}=A \cap \Sigma^s$ is determined. If $s$ is large enough, one
can in time polynomial in $s$ see if a requirement was satisfied at a
much earlier stage (which may involve checking if some short strings  are in
given computable sets). Then at stage $s$ one can react accordingly,
e.g.\ by starting to work on a different requirement.

 Slaman and  Shinoda \citelab{Shinoda.Slaman:91}{1} gave an interpretation of  $\ThN$ in
 $\Th (\Rec^p_T)$,
but left open the case of polynomial time many-one degrees. Three
years later, Ambos Spies and Nies \citelab{Ambos.Nies:92}{1} proved that 
 $\Th(\Rec^p_m)$  is undecidable.
However, the two latter  results
use the so-called ``speed-up technique''\id{speed-up technique} introduced by Ambos-Spies, a method which 
leads to 
 computable sets of very 
high complexity (usually nonelementary sets). 
From a complexity theorist's point of view, such sets are not very relevant
because they are only computable in an ideal sense. Therefore here we
will consider degree structures based on sets of low complexity. Let

$$\mbox{\DTIME}(h):=  \{ X \sub \{0,1\}^{<\omega} : X \ \mbox{can be computed
  in time} \ O(h) \}.$$ 

 A function $h : \NN \mapsto \NN$
  is  \iid{time constructible} 
if $h(n) $ can be computed in time $O(h(n))$ (here we identify $\NN$
with $\{0\}^{<\omega}$). 
We will prove that, for each time constructible $h$ which dominates
 the polynomial $n \mapsto n^k$ for each $k$, $(\mbox{\DTIME}(h), \le^p_r)$ has an
undecidable theory (Downey and Nies \citelab{Downey.Nies:97}{1}). Thus, 
 for
instance  the polynomial time $T$--degrees of sets in exponential
time have an undecidable theory.

A set $A$ is \iid{tally} if $A \sub \{0\}^{< \omega}$. For  the
result of Downey and Nies mentioned above, we will in fact prove that each initial
interval $[\zer, \a]$ has an undecidable theory, where  $\a \neq \zer$ is the degree of very particular
type of a tally set, called a \ird{super sparse} set. This notion was
introduced by Ambos-Spies \citelab{Ambos:86}{SS}. One requires that $A \sub
\{0^{f(k)}: k \in \NN \}$, for a time constructible $f$ which
increases so fast that ``$A(0^{f(k)}=1$ ?''  can be determined in time
$O(f(k+1))$. These sets allow us some of the advantages of the
speed-up technique, while still existing in each class
$\DTIME(h)$, $h$ as above.

\chapter{Theories and coding}

A {\it  theory}\id{theory}
is a consistent set of first-order sentences in some language closed
under logical inference. Given a 
theory $T$ in an effective  first-order language,
 an important   first question 
 is whether the theory is decidable\id{theory!decidable}.
Such investigations were initiated by G\"odel (implicit in
 \cite{ Godel:31}) and Tarski \cite{Tarski:49} and
have played an important role ever since. 
If $\bf A$ is a structure whose theory is known to be
undecidable, an interesting further problem is to determine
 the 
computational complexity of $\hbox{Th}(\bf A)$.  If $\bf A$ can
be coded  in $(\NN,+,\times)$, an upper bound for its 
computational complexity is the degree of
$\hbox{Th}(\NN,+,\times)$ (this theory is also called {\it true 
arithmetic}\id{true 
arithmetic}), because there is an \ird{interpretation}  of $\hbox{Th}(\bf A)$ in true
arithmetic. For most of the structures introduced in the previous
chapter, we will give an interpretation in the other direction. So
$\hbox{Th}(\bf A)$
 has the same computational complexity as true arithmetic. 

A further
 question we will consider is which fragments of an undecidable
 theory $T$ are
 undecidable.

\section{Coding}

\label{Coding}

We explain coding with first-order formulas and introduce the central
concept of a coding scheme.
Consider first-order languages $L_0, L_1$ over finite symbols sets,
and suppose that $L_0$ is relational. We intend to code
$L_0$-structures $\C$ into $L_1$-structures $\A$, by  using
an appropriate collection 
 of
$L_1$-formulas.
We represent elements of $\C$ by elements  in an $\A$-definable
set $D$, modulo an  $\A$-definable equivalence relation $\eqq$
($\A$-definable means definable in $\A$ with parameters). Then the relations of
$\C$ give rise to corresponding relations on $D/_\eqq$, which we also
require to be $\A$-definable. The uniformity is embodied in the fact  that all the
definability requirements are satisfied via a fixed collection of formulas,
called a scheme. Thus, 
a \iid{scheme} {\it
  for  coding
$L_0$-structures into $L_1$-structures}  is given  by  a
collection of $L_1$-formulas

\begin{equation} 
  \label{general scheme}
  S=\phi_{dom}(x; \wt p), \phi_\eqq(x,y; \wt p), (\phi_R(x_1, \ldots,
  \x_n; \wt p))_{R \mmbox{\small relation symbol of } L_0},
\end{equation}

together with a \iid{correctness condition} $\alpha(\wt p)$ which
expresses to the least  that actually an $L_1$ structure is coded by $\wt
p$. The correctness condition states that 

\bi
\item $D=\{x: \phi_{dom}(x; \wt
p)\}$ is nonempty,
\item  $\eqq \ = \{x,y: \phi_\eqq(x,y; \wt p)\}$ is an equivalence relation
when restricted to $D$, and 
\item  the relations on $D$ defined by the
formulas $\phi_R$ are compatible with $\eqq$.

\ei

 We say that $\C$ is {\it
  coded}\id{coding!of structures} in   $\A$ via $S$ and a list of parameters $\wt
a$ in $\A$ if the structure defined by $S$ with these parameters on
$D\beqq$ equals $\C$. 

Coding of this kind was introduced  to prove in an indirect
way that the theory of a class of $L_1$-structures is undecidable.
For {\it uniform
coding}\id{coding!uniform}  (up to isomorphism) 
of a class $\mathcal C$ of $L_0$-structures in a class $\mathcal A$
of $L_1$-structures one requires that via a fixed scheme of formulas with
parameters  a copy of each structure from $\mathcal C$ in
some structure $\A$ from $\mathcal A$ can be coded if appropriate values in $\A$
are substituted for the parameters.

 Given a first-order language $L$,  $L$-\iid{valid} is the set of valid
 $L$-sentences. A theory $T\sub L$ is \iid{hereditarily
  undecidable} (h.u.)  if, for each $X$, 
$$ L-\text{valid} \sub X \sub T \RRA X \ \text{undecidable}.$$
The following well-known fact (see for instance Burris and
Sankappanavar \citelab{Burris.Sankappanavar:81}{1}) is used to transfer hereditary undecidability of theories of classes.

\begin{fact} \label{indundec} 
If $\Th(\CCC)$ is \text{h.u.} and $\CCC$ can be uniformly
  coded in $\AAA$, then $\Th(\AAA)$ is \text{h.u.} \end{fact}

 For instance, to show that the
theory of the structure of r.e.\ $m$-degrees is undecidable,  one can use for $\mathcal C$ 
the class of finite distributive lattices, viewed as partial orders: each such
lattice is isomorphic to an initial interval $[\zer ,{\a}]$ of the r.e.\
$m$-degrees (Lachlan \citelab{Lachlan:70*1}{fin}). So $\mathcal C$ is uniformly coded in the class
$\{\Rm\}$. Since $\Th(\CCC)$ is known to be h.u.,  $\Th(\Rm)$
 is undecidable. 

Clearly, in $\Nops$ (or in fact in any model of Peano arithmetic) one
can uniformly code, say, the class of finite undirected graphs\id{graph!undirected}. By
Theorem \ref{Lavrov}  below, this class has a h.u.\ theory, so by Fact
\ref{indundec} 
$\Th\Nops$ is h.u. So, as a special case of Fact \ref{indundec} , we obtain

\begin{corollary}[\cite{Burris.Sankappanavar:81}] If $\Nops$ can be coded in a structure $\A$ with parameters, then
  $\ThA$ is h.u.   \label{inundecN} \label{indundecN} \end{corollary}

We will give more details on this method when we discuss undecidable
fragments in Section \ref{fragments} below. Next we consider
 interpreting $\ThN$ in $\ThA$.  Here a  central notion is the
following. 

\begin{example}
\label{codepa} A scheme $S_{M}$ for coding models of some finitely
axiomatized fragment $PA^{-}$ of Peano arithmetic (in the language $L$ $%
(+,\times )$) is given by the formulas

\begin{equation}
\label{SMformulas}
\varphi_{num}(x,{\overline{p}}),\phi_{\eqq}(x,y,\ol p),\varphi _+(x,y,z;{\overline{p}}),\varphi
_\times(x,y,z;{\overline{p}}) 
\end{equation}
and a \ird{correctness condition} $\alpha_0 ({\overline{p}})$ which says that 
\bi
\item $\phi_{\eqq}(x,y,\ol p)$ defines an equivalence relation $\equiv$ on $\{x:\varphi _{num}(x;{%
\overline{p}})\}$
\item $\varphi
_+$ and $\varphi _\times$ define binary operations on the set $\{x: \varphi_{num}(x;{%
\overline{p}})\}$ which are compatible with $\equiv$
\item $\{x:\varphi _{num}(x;{%
\overline{p}})\}/_\equiv$ with the corresponding operations satisfies
the \\ finitely many axioms of $PA^{-}$.
\ei 
\end{example}

(Formally, we view $L(+, \times)$ as a language with two ternary
relation symbols.)
 In our
applications, the axioms ensure that $\M$ has a standard part. For
instance think of $PA^{-}$ as Robinson arithmetic $Q$. In some
applications it is necessary to represent numbers by equivalence
classes of tuples of a fixed length (as opposed to elements). Thus the
coding is similar to the coding of $\QQ$ in $\ZZ$ given by the
quotient field construction, where a rational is represented by a 
pair of integers  (but we may also use parameters). To adapt the
definitions, in 
the above $x,y,z$ have to be interpreted as tuples of variables.

Notice that we are
now  interested in the collection of coded  structures as they are, not only in structures
up to isomorphism. In fact we will code more general objects then structures
into $\A$: we drop the condition that there be a domain formula and a
formula $\phi_\eqq$. Thus an {\it object scheme}\id{scheme!object-}
 for coding in an $L_1$-structure
is given by a list of $L_1$-formulas

$$ \phi_1, \ldots , \phi_n$$

with a shared parameter list $\ol  p$, together with a correctness condition $\alpha(\ol p)$. 

\begin{example} 
A scheme $S_{g}$ for defining a function $g$ is given by a formula 
 $\varphi _{1}(x,y; \ol p)$ defining
the relation between inputs and outputs; and a correctness condition $\alpha
(x,y; \ol  p)$ which says that a function is  defined: $\fa x \ex^{\le
  1} \phi_1(x,y; \ol p)$
\end{example}

\begin{example} \label{rels}  We will often consider object schemes for classes of $n$-ary relations on $\A$. Such a scheme is given
  by a formula $\phi(x_1, \ldots, x_n; \ol  p)$ and a correctness
  condition $\alpha(\ol p)$.
\end{example}

For instance, if $\A$ is a \ird{linear order}  and $\mathcal{C} $ is
the set of closed intervals, then $\mathcal{C}$ is uniformly definable via the scheme
consisting  of $\phi_1(x; a,b)\LR a \le x \le b$ and the correctness condition
$\alpha(a,b)\LR a \le b$.

In general, an object  scheme $S_X$ introduces a new type of object. The
parameters $\overline{p}$ satisfying $\alpha (\overline{p})$ code an object,
and $S_X$ acts as a decoding key\id{coding!decoding key}. Using this coding, it becomes possible to
quantify over objects of the new type (a form of second order
quantification) in the first--order language of $\A$. Thus one can
quantify over uniformly definable classes in the sense of the
following definition.

\begin{definition}
\be \item[(i)] A class $\mathcal{C}$ of objects of a common type is  \ird{uniformly definable} in $%
\A$ if, for some scheme $S$, $\mathcal{C}$ is the class of
objects coded via $S$ as the parameters range over tuples in $\A$
which satisfy the correctness condition. 

\item[(ii)] $\mathcal{C}$ is {\rm weakly uniformly
  definable}\id{uniformly  definable!weakly} if $\mathcal{C}$ is
contained in a uniformly definable class.
\ee

\end{definition}

 We can  perform basic mathematical operations on objects of two
possibly different types and obtain a uniform way of coding objects of a yet
different type. For example, we can define a scheme $S$ for the
compositions $g \circ h$ of maps $g,h$ defined by schemes $S_g,S_h$.
Furthermore, we can express basic relationships between coded objects by
first order conditions on codes; for instance we can express the
relationship ``$g$ is a partial map from $\M_0$ to $\M_1$'', where $\M_0$, $\M_1$
are coded via $S_M$ and $g$ is coded via $S_g$, by formulas of $\mathcal{R}$%
.

\begin{notation} \label{letter convention} {\rm  We use the following convention throughout: If a scheme $S_X$ is given,
variables $X$, $X_0$, etc. denote objects coded by this scheme for a
particular parameter list $\overline{\mathbf{p}}$ satisfying the correctness
condition. If it is necessary to mention the parameters explicitly, we write 
$X(\overline{\mathbf{p}})$ (or $X_{\overline{\mathbf{p}}}$,
$X_0(\overline{\mathbf{p}})$, etc. We say\id{coding!via a scheme}
 that
 $\overline{\mathbf{p}}$ codes $X(\mathbf{\bar p})$ via $S_X$.}
\end{notation}

 We will use the term ``scheme'' to refer to either a coding \ird{scheme} or
 an object scheme\id{scheme!object-}. It will be clear from the context which notion is meant.

\section{Interpreting true arithmetic}
\label{sectionITA}

In the following we assume that $\A$ is a structure which can be coded
in $\NNops$. 
Note that  there
is an onto map $\gamma:\NN\to \A$ such that the preimages of the
relations and functions of $\A$ are arithmetical.  For instance, if $\A$
is $\Rm$, let $\gamma(i)=\deg_m(W_i)$. We call the preimage of a relation
$R$ under $\gamma$  the corresponding {\it index
  relation}\id{index!-relation} (index set\id{index!-set}
 if $R$ is unary) and sometimes
write $\Theta R$ for this preimage. In fact  we will often identify
$R$ and $\Theta R$.

To interpret $\ThN$ in $\ThA$, for $\Rm$ and $\EE$ we carry out the following two steps:

\begin{eqnarray} \label{code} & \mbox{Specify a scheme $S_M$ as in Example \ref{codepa} and a special
  list}\ \ol \p \ \\
 &  \mbox{such that }\ \M_{\ol \p}  \ \mbox{is standard}  \nonumber
\end{eqnarray}

\begin{eqnarray} \label{see}  & \mbox{Find  an additional
    correctness condition } \ \alpha_{st}(\ol p) \\
& \mbox{which holds iff} \  \ol \p \ \mbox{code a copy of} \ \Nops \nonumber
\end{eqnarray}\id{correctness condition}
 
The point is that the condition $\alpha_0(\ol p)$ from Example
\ref{codepa} only gives an approximation to standardness. While
(\ref{code}) can be seen as a  ``local'' coding, relying on very
special parameters, recognizing
standardness of an {\it arbitrary} $\M_{\ol p}$ in a first-order way depends on
how $\M_{\ol p}$ relates to its context, namely the whole structure
$\A$. Thus for (\ref{see}) we will use  particular
properties of  $\A$. 

If $\beta$ is a sentence in the language of arithmetic, let $\tilde
\beta(\ol p)$ be the translation of $\beta$, namely the formula obtained
by replacing $=, +, \times$ by their definitions via $\phi_{\eqq}, \phi_+,
\phi_\times$ and relativizing the quantifiers to those $x$   satisfying
$\phi_{num}(x,\ol p  $. Then

\begin{equation}
\label{Ntrans}
 (\NN, + , \times) \models \beta \ \LR \AA \models \ex \ol p [
\alpha_{st}(\ol p ) \lland \tilde \beta(\ol p)].
\end{equation}
Since $\tilde \beta(\ol p)$ is obtained in an effective way,  $\ThN \le_m \ThA$.

In order to carry out (\ref{code}), it is often useful if one just has
to code a computable  directed graph\id{graph!directed} $(V,E)$ into  $\A$ with parameters (where
in fact $V=\NN$). Here we 
provide a parameterless coding of $\Nops$ in a particular such
graph\id{graph!coding $N$}

\begin{equation}
\label{VN}   (V_\NN, E_\NN),
\end{equation}

  which is
a recursive irreflexive \ird{partial order}. To construct this
 partial order, one starts with a countable antichain of minimal elements $p_n$ which will represent the numbers $n$.  Then, for each $n, m \in \NN$ one adds an element
$c_{n,m}$ to $P_A$ which represents the pair 
$(p_n,p_m)$.  Next, one adds ascending chains of lengths 2 and 3, respectively, from $p_n$ to $c_{n,m}$ and from $p_m$ to $c_{n,m}$.  Finally, to code addition, add a chain of length 4 from $p_{n+m}$ to
$c_{n,m}$ and for multiplication, add a chain of length 5 from
$p_{n\times m}$ to $c_{n,m}$.

We now discuss how to carry out 
 (\ref{see}),  assuming that some scheme $S_M$ as in  (\ref{code}) has
  been specified. We first assume that $\phi_{\eqq}(x,y; \ol p)$  defines
  the
  trivial equivalence relation $x=y$ (thus, numbers are represented by
  certain elements of $\A$). If we were allowed to  quantify over subsets of
 $\M_{\ol p}$, for any $\M_{\ol p}$, then we could simply  use \ird{Dedekind}'s second-order axiomatization of
 $\Nops$: we would require that each subset of $\M_{\ol p}$ which contains $0^{\M_{\ol p}}$ and
 is closed under successor equals $\M_{\ol p}$. Of course we cannot  quantify
 over all such subsets in the first-order language of $\A$, but we can
 try to quantify over sufficiently many, by using some uniform definability result. The
 following two facts specify which subsets must be included.

\begin{fact} \label{StSk} Suppose $\A$ is coded in $\Nops$. Then,
  for some fixed $k$, the standard
  part $S$ of each coded model $\M$ has a $\Sk$ index
  set. \end{fact}

\pf Note that $\gamma(i) \in S \ \LR \  \ex n\in \NN \ \ex y_0, \ldots, y_n\in \A$
  $$ [y_0 =  0^\M \lland y_n =
\gamma(i)  \lland (\fa i < n) \A \models \phi_\oplus(y_i, 1^\M,
y_{i+1})].$$

Since $\A$ is coded in $\Nops$, this is a $\Sk$ property of $i$, for some
fixed $k$ depending only on $\A$ and the scheme $S_M$.
\eop

We call a subset of $\M$ a $\Sk$-{\it subset}  if its index set is $\Sk$. 
Make sure not to  confuse $\Sk$-subsets of $\M$
with sets which can be defined by a  $\Sk$-formula  from the point of view of
$\M$.

\begin{fact} 
\label{DITA}
Suppose 
 that the collection of $\Sk$ subsets of any $\M_{\ol p}$ is weakly 
uniformly definable via a formula $\phi(x; \ol q)$, where $\ol q$ is a
parameter list containing $\ol p$. 
Then $\ThN$ can be interpreted in $\ThA$.
\end{fact}

\pf
 Let
$\alpha_{st}(\ol p)$ be the formula expressing

\begin{quote}  for all $\ol q$, if $\{x: \phi(x;\ol q) \}$ is a subset
  of $\M_{\ol p}$ which  which contains $0^{\M_{\ol p}}$ and
 is closed under successor then it equals $\M_{\ol p}$. 
\end{quote}

This is certainly satisfied if $\M_{\ol p}$ is standard. If $\M_{\ol p}$ is not
standard, then, the standard part is a $\Sk$--set which therefore can
be defined via some $\ol \q$. So the statement fails. \eop

If $\phi_{\eqq}$ defines a nontrivial equivalence relation, with an adjustment of the
terminology carried out in the following definition the previous considerations are still valid. 

\begin{definition}
\label{repr}
A class $\hat S
\sub \{x: \phi_{num}(x, \ol p)\}$ {\rm  represents}\id{represent} a subset $S \sub
\M_{\ol p}$ if $\hat S$ is $\equiv$-closed and $\hat S/_\equiv = S$. We call
$S$ a $\Sk$-subset of $\M_{\ol p}$ if $\hat S$ has a  $\Sk$ index set.
\end{definition}

 (Note that in the above we really mean the restriction of $\equiv$ to

$\{x:  \phi_{num} (x, \ol p)\}$.) To show that the $\Sk$-subsets of $\M_{\ol p}$ are weakly uniformly definable\id{uniformly definable!weakly}
usually involves proving a
sufficiently strong uniform definability result. Such a result  can be
derived for $\Rm$, as well as for $\EE$.
 The definability lemma\id{definability lemma} for $\Rm$ states that

\begin{eqnarray} \label{definability lemmaemmaRm}
  & \mbox{for each $k\ge 3,N \ge 1$, the class of $N$-ary $\Sk$
    relations which} \\
& \mbox{ are contained in some $[\zer,\c]$, $\c < \one$, is uniformly
    definable.} \nonumber 
\end{eqnarray}

It was first proved  in Nies \citelab{Nies:95}{DL} by induction on $k$ (see
Section \ref{definability lemmaRm}). Later, Harrington used the same general method to
prove a similar result  for $\EE$, which we call the ideal
definability lemma\id{definability lemma!ideal}: for an r.e.\ set $A$,
let $\BB (A)$ be the Boolean  algebra of components of c.e.\ splittings
of $A$, and let $\RR (A)$ be the ideal of $\BB (A)$ consisting of the
computable subsets of $A$.  An ideal $I$ of $\BB (A)$ is called {\it
k-acceptable} if $\RR (A)\subseteq I$ and $\{ e:W_e\in I\}$ is $\Sigma^0_k$.
Harrington's  
 ideal definability lemma  asserts  that, for any odd $k\ge 3$, 

 \begin{equation}
   \label{definability lemmaE}
  \begin{array}{c}
   \text{the class of} \  
k-\text{acceptable ideals of}  \\ 
\BB(A) \ttext{is uniformly definable.}
  \end{array}
 \end{equation}
  Again, it is proved by induction,
here over odd $k \ge 3$ (see Section \ref{ideal definability lemmaE}).

The definability lemma\id{definability lemma}  for $\Rm$ is in fact so strong that it can also
be used to  code copies of $\Nops$
with parameters, i.e., to carry out (\ref{code}). For $\EE$, some extra work
is required, but the ideal definability lemma is still a main ingredient. In this way, an
intermediate coding in 
cumbersome auxiliary structures, like  the recursive \ird{boolean pair}s in
Herrmann 
\cite{Herrmann:84}, can be avoided.
 
For the structure $\DTE$ of $T$-degrees of $\Delta^0_2$ sets (Shore \citelab{Shore:81}{1}), as well as for
$\Rtt$ (Nies and Shore \citelab{Nies.Shore:95}{1}), one also
satisfies (\ref{code}) and (\ref{see}), but no  general definability
lemmas\id{definability lemma}  are used to interpret $\ThN$ in the
theory. Instead, the coding of copies $\M$ of $\Nops$ is made more ``effective'' (in the
sense of  the arithmetical complexity of a function $g$ such that
$n^\M = \deg_r(W_{g(n)})$), so that, e.g.\
in the case of $\Rtt$, the standard part of any coded $\M$ is actually
$\SIII$. Now, a rather weak definability result  suffices: each $\SIII$-ideal of $T$-incomplete
c.e.\ $tt$-degrees has an \iid{exact pair}, namely it has the form $[\zer,\a] \cap [\zer,\b]$ for
appropriate $\a,\b\in \Rtt$.

In the case of $\RT$,  another approach yet has been
carried out in order to satisfy (\ref{see}): one considers not only coded
copies of $\Nops$, but also coded partial isomorphisms between them (called
comparison maps). Extending work by Slaman and Woodin, in 
Nies e.a.\ \cite{NSS:96},  schemes $S_M,S_g$ are determined  such
that for each coded copies $\M_0, \M_1$ of $\Nops$  there is a map $g$ which extends the
isomorphism between the standard parts of the coded models. Then,
$\Mp$ is standard iff for each $\M$, some $g: \Mp \rightarrow \M$ is
total. The latter  condition  can be expressed in the first-order language
of $\RT$. Thus, standard models are singled out as the ``shortest''
coded models.
The idea of  using ``comparison maps''
 is  essential  to obtain the definability results in
Nies e.a.\  \citelab{NSS:96}{def}, for instance the definability without parameters of $\Low_2$ and $\High_1$.

For $\Rwtt$ we will develop in Chapter  \ref{ChapterRwtt} a parameter free
coding of a copy of $\Nops$. We represent the number $n$ by all sets
of cardinality $n$ in  the  uniformly definable class of
``EN-sets''. A  scheme is needed to code  maps between
EN-sets in order to express for instance that two EN-sets have the
same cardinality.  A similar result can be obtained for $\Rm$ (Section
\ref{CharRm}). The first application of this variant was to the  upper
semilattice of c.e.\ equivalence relations\id{equivalence relation!c.e.}
modulo finite variants (Nies \citelab{Nies:94}{1}). 

\section{Undecidable fragments of  theories}
\label{fragments}

In the context of fragments of theories we only consider coding {\it up to
isomorphism}. A formula is $\SSk$\id{formula!$\Sigma_k$} if it has the form $$ (\ex \dots \ex) \  (\fa \dots \fa)
\ 
(\ex \dots \ex)  \ldots \  \psi,$$  with $k-1$ quantifier alternations and
$\psi$ quantifier free, and $\Pi_k$\id{formula!$\Pi_k$}  if it has the form

 $$(\fa \dots \fa) \ 
(\ex \dots \ex) \  (\fa \dots \fa)  \ldots \ \psi.$$ 

 Given an  undecidable theory $T$, an
interesting further question to ask is which {\it
  fragments}\id{theory!fragment of}  $T \cap \SSk$
and $T \cap \PPk$ are undecidable, for several reasons. Firstly, the
sentences which occur in mathematical practice usually have a low
number of quantifier alternations. So, even after undecidability of $T$ is known, the
question remains which feasible fragments are undecidable. Secondly,
a sharp classification at which fragment  an undecidable
theory $T$  becomes undecidable gives more precise  information about 
$T$ than  a plain undecidability proof. (Ershov gave an example of an undecidable
theory
where all fragments are decidable. However, an undecidability result 
obtained indirectly via Fact \ref{indundec} gives actually
undecidability of 
 some fragment.)
Finally, if $T=\Th(\CCC)$ for some class
of structures $\CCC$, sometimes one  can
 interpret  the sentences in a fragment algebraically. Then a
decision procedure for that fragment gives algebraic information about
$\CCC$. For instance, the $\Pi_1$-theory of a variety is closely
connected to the word problem of its finitely presented
members. Moreover, $\Pi_2$-sentences in the language of p.o. can be
interpreted as  statements about possible extensions of embeddings of
finite \ird{partial order}s.

In the following we will develop a version for fragments of the method to obtain undecidability
of theories of classes in an indirect way which was outlined in
Section \ref{Coding}. Given a first-order language $L$, a set of sentences  $U\sub L$ is \iid{hereditarily
  undecidable} (h.u.) if, for each $X$, 
$$ L-\text{valid}\cap U \sub X \sub U \RRA X \ \text{undecidable}.$$

This extends the previous definition given before Fact \ref{indundec}
for {\it theories} $U$.
Disjoint sets $A,B \sub \NN$ are called \iid{recursively inseparable}
if there is no computable set $R$ such that $A \sub R \lland  B \sub
\NN-R$. Then (provided we have chosen some G\"odel numbering of the
formulas in $L$), 
 $$U\sub L \ttext{is h.u.}
\LR \  L-\text{valid} \cap U, L-U \ \text{are recursively inseparable.}$$

In order to obtain an undecidability result for a low-level fragment
of $\Th(\AAA)$ from the
method in Fact \ref{indundec}, one has to
invent a coding of $\CCC$ in $\AAA$ of maximum economy. Therefore it
is useful to consider a class $\CCC$ in  a language $L_0$ without
equality. The following theorem was proved by Lavrov and is reproved
in Nies \citelab{Nies:96*1}{1}.

\begin{theorem}   \label{Lavrov} The $\Sigma_2$-theory of
  the class of finite undirected graphs\id{graph!undirected} in the language without
  equality
is hereditarily undecidable. \end{theorem}

As in Section \ref{Coding},  consider first-order languages $L_0, L_1$ over finite symbols sets,
and suppose that $L_0$ is relational and has no equality symbol.
First we have to clarify when an $L_0$-structure  is said to be coded in an
$L_1$-structure. Let $\phi_{eq}(x,y) \in L_0$ be the formula  expressing that $x,y$ behave in the
same way with respect to all elements of the structure. Thus
$\phi_{eq}(x,y)$ is the conjunction of formulas of the type

$$\fa x[ (Rxz \LLR Ryz) \lland (Rzx \LLR Rzy)],$$

for each relation symbol $R$ of $L_0$. Given an $L_0$-structure $\C$, let

$$eq(\C) = \{x,y: \C \models \phi_{eq}(x,y)\},$$

 and let
$\C/_{eq(\C)}$ be the structure on equivalence classes defined in the
canonical way. It is easy to verify that $$\C\models \psi \LLR
\C/_{eq(\C)}\models \psi,$$ for each $L_0$-sentence  $\psi$, by an
induction on $|\psi|$.   

A {\it $\Sigma_k$-scheme}\id{scheme!$\Sigma_k-$} is given by a list of
formulas
\begin{equation} 
  \label{Skscheme}
  S=\phi_{dom}(x; \wt p), \phi_R(x_1, \ldots, x_n; \wt p),  
\phi_{\ol R}(x_1, \dots,
  x_n; \wt p )_{R \ \text{relation symbol of} \ L_0},
\end{equation}
together with a \iid{correctness condition}  which
expresses  that 

$$D= \{x: \phi_{dom}(x; \wt
p)\}$$

 is nonempty,  and that the relations on $D$ defined by the
formulas $\phi_R, \phi_{\ol R}$ are complements of each other. 
These condition can be expressed by universally quantified boolean
combinations of $\SSk$-formulas, and therefore by a $\Pi_{k+1}$-formula $\alpha(\wt p)$.
 We say that the $L_0$-structure $\C$ is {\it
coded}\id{coding!$\SSk$-} in   $\A$ via the   $\SSk$-scheme $S$ and a list of parameters $\wt
a$ in $\A$ if $\A\models \alpha(\wt a)$ and $$\C/_{eq(\C)} \cong \D/_{eq(\D)},$$

 where $D=\{x: \phi_{dom}(x; \wt
p)\}$ and $\D$ is the $L_0$-structure on $D$ induced by the formulas $\phi_R$.  For uniform
$\SSk$-coding
of a class $\mathcal C$ of $L_0$-structures in a class $\mathcal A$
of $L_1$-structures one requires that a fixed $\SSk$- scheme of formulas with
parameters defines a copy of each structure from $\mathcal C$ in
some structure $\A$ from $\mathcal A$ if appropriate values in $\A$
are substituted for the parameters.
For example,  the class of finite undirected
graphs\id{graph!undirected}  from Theorem
\ref{Lavrov} is uniformly $\Sigma_1$-coded in the class of finite
p.o. as $L(\le)$-structures, via the following $\Sigma_1$-scheme without parameters:

\begin{eqnarray} \label{FPO}    \phi_{dom}(x) &\LR& \ex u \ex v [u <
  x < v] \\
   \phi_E(x,y) &\LR& x \not \le y \not \le x \lland \ex z \ x,y \le z
   \nonumber \\
  \phi_{\ol E}(x,y) &\LR& x \not \le y \not \le x \lland \ex z \ z \le
  x,y \nonumber 
\end{eqnarray}

(here  $x < y$ stands for $x \le y \lland y \not \le x$).

We are now ready to obtain a version of Fact \ref{indundec} for
fragments. 

\begin{transferlemma} 
\label{transfer} Let $r \ge 2, k \ge 1$.

\be
\item[(i)] If $\CCC$ can be uniformly $\Sk$-coded in $\DDD$ without
  parameters, then $$ \Sigma_r-\Th(\CCC) \ h.u. \RRA
  \Sigma_{r+k-1}-\Th(\DDD) \ h.u.$$

\item[(ii)] If $\CCC$ can be uniformly $\Sk$-coded in $\DDD$ with
  parameters, then $$ \Pi_{r+1}-\Th(\CCC) \ h.u. \RRA
  \Pi_{r+k}-\Th(\DDD) \ h.u.$$

\ee 
\end{transferlemma}

Thus, combining (\ref{FPO}) with Theorem \ref{Lavrov}, we obtain
from (i)
that the $\Sigma_2$-theory of the class of finite partial
orders\id{partial order!undecidability of $\Sigma_2$-theory}  is h.u. 

\pf The idea is to define an effective map $F$ from $L_0$-sentences
to  $L_1$ sentences which maps $L_0-\text{valid}$ into $L_1-\text{valid}$
and sentences $\phi \nii \Th(\CCC)$ to a sentence $F(\phi) \nii
\Th(\DDD)$. 

Given an $L_0$-sentence $\phi$ in normal form, the translation $\wt
\phi(\ol p)$ ($\wt \phi$ if no parameters are used in the $\SSk$-scheme) is
obtained by relativizing the quantifiers to $\{x: \phi_{dom}(x; \ol p)\}$ and
replacing the atomic formulas $R\ol x$ and $\neg Rx_1, \ldots,
x_n$ by $\phi_R (\ol x)$ and $\phi_{\ol R} (x_1, \ldots,
x_n)$ ($\ol x = x_1, \ldots, x_n$) in a way  to minimize the number of quantifier alternations. If the
innermost quantifier in $\phi$ is existential, replace $R\ol x$ by
$\phi_R (\ol x)$ and replace $\neg Rx_1, \ldots,
x_n$ by $\phi_{\ol R} (x_1, \ldots,
x_n)$. Otherwise, replace $R\ol x$ by the $\PPk$-formula 
$\neg \phi_R (\ol x)$ and replace $\neg Rx_1, \ldots,
x_n$ by the   $\PPk$-formula $\neg \phi_{R} (\ol x)$. For instance, if $\phi$ is $\ex x \fa y [Rxy \ \vee \ \neg
Ryx]$, then $\tilde \phi(\ol p)$ is 

$$\ex x[ \phi_{dom}(x; \ol p) \lland \fa y[\phi_{dom}(y;\ol p) \RRA
\neg \phi_{\ol R}(x,y; \ol p) \vee \neg \phi_R(y,x;\ol p)]].$$

Note that the translation of a $\Sigma_r$-sentence is a
$\Sigma_{r+k-1}$-formula and that the  translation of a $\Pi_{r+1}$-sentence is a
$\Pi_{r+k}$-formula.

Let  $$F(\phi) = \fa \ol p [ \alpha(\ol p) \RRA \wt \phi(\ol p)]$$

(and $F(\phi) =   \alpha \RRA \wt \phi$ if there are no
parameters). Clearly, $$\phi \in L_0-\text{valid} \RRA F(\phi) \in
L_1-\text{valid}.$$

Moreover, 

$$\phi \nii \Th(\CCC) \RRA F(\phi) \nii \Th(\DDD),$$

because, if $\phi$ fails in some structure $\C \in \CCC$, then $F(\phi)$
fails in a structure $\D \in \DDD$ coding $\C$, the counterexample for $\fa
\ol p[\ldots]$ being provided by the list of parameters used for the coding.

For the proof of (i), note that, if $\phi$ is a $\Sigma_r$-sentence, then
$F(\phi)$ is logically equivalent to a $\Sigma_{r+k-1}$-sentence, since $r+k-1
\ge k+1$ and $\alpha(\ol p)$ is a $\Pi_{k+1}$-formula. 

For (ii) we argue in a similar way, using the fact that for a
$\Pi_{r+1}$-sentence 
$\phi$, 
$F(\phi)$ is logically equivalent to a $\Pi_{r+k}$-sentence. \eop

\chapter{C.e.\ many-one degrees} 

\label{RRRmmm}
We first concentrate on the proof of the definability lemma for $\Rm$
 (\ref{definability lemmaemmaRm}) explained in Section \ref{sectionITA}, which leads to
an interpretation of $\ThN$ in $\Th(\Rm)$.
 Then we  proceed  to results  about 
$\Rm$ of  a model theoretic
nature. First we derive  a local definability
result for automorphisms. Next we strengthen the result that there is an interpretation
of $\ThN$ in $\Th(\Rm)$ by  developing a coding of $\Nops$ in $\Rm$ without
parameters. Let  $\Rme = \Rm-\{\one\}$. We show that  there is an incomplete $\e \in \Rm$ such that $[\zer,
\e)$ is an elementary submodel of $\Rme$ via the inclusion embedding. In particular, $\Rme$ (and
hence $\Rm$, since $[\zer,
\e) \cup \{\one\} \prec \Rm$) has a proper elementary
submodel\id{model!elementary submodel}, i.e.\
is not a  minimal model\id{model!minimal}  over the empty
set.  It follows from results of Slaman
and Woodin \citelab{Slaman.Woodin:nd}{2} that  $\DTE$ is both a minimal
model and a {prime model}\id{model!prime}. For the  c.e.\ degree structures except
$\Rm$, both
questions remain open. 

\section{Preliminaries}

\label{prelims}
To prove the definability lemma we interpret an elementary non-extensional set theory
in $\Rme$. A set $S \sub \Rme$ is called {\it uniformly
  computably enumerable (u.c.e.)}  if $S= \{\a_n: n \in \NN\}$ for
some u.c.e.\ sequence $(\a_n)$. We will first show the uniform definability of
the u.c.e.\  (and in particular of the nonempty finite) sets $S \sub
\Rme$ from one parameter. This leads to a definable ``element'' relation on
$\Rme$.

\begin{theorem} \label{uce}

The u.c.e.\ sets of incomplete $m$-degrees
are uniformly definable via a formula $\phi_\in(x;a)$ with one parameter. \end{theorem}

Notice that the parameter $a$ is not uniquely determined. Hence
the elementary set theory we interpret in $\Rm$ will not be extensional.

An {\it ideal}\id{upper semilattice!ideals in}\id{ideals!in upper
  semilattice} of an upper semilattice is a nonempty subset which is
closed downward and under supremum. 
If $I$ is an ideal in $\R_m$, we will say that $\b$ is a {\it strong
  minimal cover of }\id{strong minimal cover}\id{ideals!strong minimal
  cover of}   $I$ if $I = [\zer,\b)$.
 A strong minimal cover is necessarily {\it join irreducible}\id{join
   irreducible}, namely it is not the supremum of two smaller degrees. We first state a Lemma which is a
special case of  Theorem 3.1 in Ershov and Lavrov
\citelab{Ershov.Lavrov:73}{1} (where a
completely different notation is used). Inspection of their proof
shows that the  
strong minimal cover is obtained in an effective way.

\begin{lemma}[\cite{Ershov.Lavrov:73}] \label{Ershov} Suppose that $I \sub \Rme$ is a
$\SIII$-ideal and $D\sub \Rme$ is a u.c.e.\ set. Then one can effectively obtain a strong minimal cover $\b$ of $I$ such that 
$\b \not \le \d$ for each $\d \in D$. \eop \end{lemma}

{\it Proof of Theorem \ref{uce} } We have to determine a formula
$\phi_\in(x;a)$ with the following property: given a u.c.e.\ sequence $(\a_n)$ of incomplete $m$-degrees, there is a parameter $\a$ such that $\{\a_n: n \in \NN\}=\{\b: \RR_m \models \phi_\in(\b;\a)\}$.  Applying Lemma \ref{Ershov}
to each $\SIII$-ideal $[0,\a_n]$ and the u.c.e.\ set  
$D= \{\a_n: n \in \NN\}$, we obtain a u.c.e.\ sequence
$(\b_n)$ such that, for each $n$, $\b_n$ is a strong minimal cover
of $[0,\a_n]$ which is not below $\a_m$ for any $m $. For each $n,m$ ,
$\b_n = \b_m$ or $\b_n, \b_m$ are incomparable. 

 Recall that, by a result of Lachlan (see \cite[p.\ 45]{Soare:87}), the complete c.e. $m$-degree is join irreducible in
$\Rm$. So $I \sub \Rme$.  Now, let
$\a$ be a strong minimal cover of the $\SIII$-ideal $I$ generated by
$\{\b_n: n \in \NN\}$.
The set $\{\b_n: n \in \NN \}$ is definable in $[0,\a)$ as the set of 
maximally join irreducible elements:
if $\x < \a$, then $\x \le \b_0 \vee \ldots \vee \b_n$ for some
$n$. By the distributivity of $\RR_m$
(\ref{ldistr})\id{distributive!upper semilattice} , there exists $\c_i \le \b_i$ such that $\x = \c_0 \vee \ldots \vee \c_n$. If $\x$ is join irreducible, this implies that 
$x \le \b_i$ for some $i$, and if $\x$ is maximally join irreducible, then even 
$\x = \b_i$.

To show, conversely, that each $\b_m$ is maximally  join irreducible in $[0,\a)$,
suppose that $\b_m \le \x$ for some join irreducible $\x < \a$. In the same way as
above, we obtain $\x \le \b_i$ for some $i$. Hence $\b_m \le \b_i$ and
therefore $\x = \b_m$. 

From the definability of $\{\b_n: n \in \NN \}$ in $[0,\a)$ we obtain the definability
of $\{\a_n: n \in \NN \}$. The formula $\phi_\in(x;a)$ is given by

$$\ex b ( b \ \mbox{maximally join irreducible in }[0,a) \ \land \ [0,b)=[0,x]).$$

 \eopnospace

\begin{corollary}

\be

\item[(i)] Every finite set $F \sub \Rme$ is uniformly definable from an incomplete $m$-degree $\a$, which can be obtained 
effectively in $F$.

\item[(ii)] There exist definable projection functions $pr_1$, $pr_2$ so that

$$ \fa \x_1,\x_2 < \one \ex \y < \one [pr_1(\y)=\x_1 \ \land \ pr_2(\y)=\x_2].$$

Moreover, an index for such a degree $\y$ can be obtained 
effectively in indices for $\x_1,\x_2$.
\ee
\label{pairs}
 \end{corollary}

\pf (i)  Suppose that $F=\{deg_m(W_i): i \in D_z\}$, where $D_z$ is a strong index for a finite set. If $D_z=\ES$, let $\a=0$. Else in some 
effective way obtain a u.c.e.\ sequence $(\a_n)$ such that $F=\{\a_n: n \in \NN\}$ and apply Theorem \ref{uce}.

(ii) Recall that, in set theory, the ordered pair $\la x_1,x_2 \ra$ is
represented  by 
$\{\{x_1\},\{x_1,x_2\}\}$. Regardless of extensionality, the analog in
$\Rme$ of
such an
object will determine both of its
components. 
We let $pr_i(\y) = \x_i$ ($i=1,2$) in case $\y$ is a set  of this form, with respect to the element relation defined by the formula $\phi_\in$, and  $pr_i(\y) = \zer$, otherwise. The maps
$pr_i$ are clearly definable in $\Rm$. Moreover, given $\x_1,\x_2 \in \Rme$, by repeated applications of (i) we can effectively in indices for $\x_1,\x_2$ determine a $\y < \one$ such that 
$pr_i(\y) =\x_i$ for $i=1,2$. \eop

If $\x_1,\x_2$ are given by indices for c.e.\ sets, let $p(\x_1,\x_2)$ denote the degree $\y$ obtained at the end of the previous proof. If $\x_1,\x_2 < \one$, then $\x_1 \vee \x_2 \le p(\x_1,\x_2) < \one$. Keep in mind that $p(\x_1,\x_2)$ really depends on the indices
used to represent $\x_1, \x_2$.

\section{The definability lemma  for $\protect\Rm$}

\label{definability lemmaRm} 

This section and the following\id{definability lemma}  do not use any particular property of $\Rm$ beyond that the ordering is $\SIII$ as a relation on indices and Theorem \ref{uce}.

\begin{definability lemma} For each $k \ge 3, N \ge 1$, the class of $\Sk$ relations
  on intervals $[\zer, \c]$ of $\Rm$  such that $\c < \one$ is weakly uniformly
  definable\id{uniformly definable}. 
In fact, there exist  formulas $\phi_{k,N}(x_1, \ldots, x_N;c,a)$ 
with the following property: \begin{quote}  if $Z \sub [\zer,\c]^N$ is $\Sk$, then one can effectively 
in $\c$ and a $\Sk$-representation of $Z$ determine an $\a< \one$
such that $\phi_{k,N}$ defines $Z$ in $\Rm$ with the parameters
$\c,\a$. \end{quote}
                    
\end{definability lemma}

\pf   We first make the extra assumption that we also possess  a lower
bound $\d > \zer$ for $Z$, 
namely that $Z \sub [\d,\c]^N$. Thus we will construct formulas
$\phi_{k,N}$ as in the statement of the 
lemma, but with an additional parameter $d$, and we will obtain $\a$ effectively in $\c,\d$ and a $\Sk$-representation of $Z$.

The formulas $\phi_{k,N}$ are defined by recursion over $k$ ($N$
 is  fixed). For notational simplicity, we
will carry our the recursion for $N=2$. The fact that $\a$ can be
determined effectively is needed to make the recursion work.

First let $k=3$. We will define u.c.e.\ sequences of incomplete 
$m$-degrees $\a_n,\b_n$ such that

\begin{equation} \label{Z} Z=\{\la \a_n,\b_n \ra: n \in \NN\} \cap [\d,\one)^2. \end{equation}

Recall that $\Theta Z$ is the index relation associated with $Z$. 
Since $\Theta Z$ is $\SIII$, there is a u.c.e.\ sequence of sets
$(X_{\la i,j,k \ra})$ such that
 ${\la i,j \ra} \in \Theta Z  \RA  \ex k X_{\la i,j,k \ra}=\NN$  and
 $\la i,j \ra  \not\in \Theta Z  \RA  \fa k \ X_{\la i,j,k \ra}$
 finite. This follows e.g.\ from Soare \cite[p.\ 68]{Soare:87}. Now, for $n=\la i,j,k \ra$, 
let

 $$A_n=W_i \cap X_{\la i,j,k \ra} \ttext{and} B_n = W_j \cap X_{\la
  i,j,k \ra}. $$

 If $\a_n= \deg_m(A_n)$ and $\b_n= \deg_m(B_m)$,
then $\d \neq \zer$ implies (\ref{Z}).

Since $\{p(\a_n,\b_n): n \in \NN\}$ is a u.c.e.\ set of incomplete
degrees, by Theorem \ref{uce} this set can be defined from a parameter
$\a$. Then $Z$ is definable from the parameters $\d,\a$ via a fixed formula
$\phi_{3,2}(x_1,x_2;d,a)$ (the upper bound was not  needed yet):

\begin{eqnarray*}
  \label{nono}
  \phi_{3,2}(x_1,x_2;d,a) & \equiv & x_1 \ge d \lland x_2 \ge d \lland
  \\
              &  & \ex
y[\phi_\in(y,a) \lland pr_1(y)=x_1 \lland pr_2(y)=x_2].
\end{eqnarray*}

 Clearly $\a$ was obtained effectively in $\d$ and a $\SIII$-representation of $\Theta Z$.

Next suppose $\Theta Z$ is $\Sigma^0_{k+1} (k \ge 3)$. We will show the definability of $Y= [\d,\c]-Z$. Notice that $\Theta Y$ is $\Pi^0_{k+1}$. One might attempt the following:

\be

\item write $\Theta Y= \bigcap_{n\in \NN} R_n$, where each $R_n$ is an $N$-ary relation on indices, closed under $\equiv_m$ and $\Sk$ uniformly in $n$

\item use the inductive hypothesis to define from parameters $\d,\c $ and $\b_n$ the relation on $m$-degrees given by $R_n$ 

\item finally apply Theorem \ref{uce} to the u.c.e.\ sequence
$(\b_n)$ in order to define $Y$ from the new parameter $\b$ obtained,  as the
intersection of these relations.

\ee 
 But there is a flaw in this approach,
 the careless handling of index sets. For instance, if $k=4$ and
 $Y = [\zer,\v)$ for some $\zer < \v$, then such a representation

implies that $\Theta Y = R_n\cap \{e:  \deg_m(W_e) \le \v\}$,
for some $n$ and hence that $Y$ is $\SIII$. This  is not  the case
when $[\zer,\v]$ is effectively isomorphic to $\Rm$.
But the approach can be rescued if  carried out separately for 
each pair   of $m$-degrees. The relation $Y$ is the effective union of $\Pi^0_{k+1}$ relations $Y_{\la i,j \ra}$ of cardinality $\le 1$, where $$Y_{\la i,j \ra}  = Y \cap \{\la \deg_m(W_i), \deg_m(W_j)\ra\}.$$

We will find a u.c.e.\ sequence  $(\a_{\la i,j\ra})$ so that $a_{\la
  i,j\ra}) <\one$ and for some fixed formula $\psi(x,y;a)$, each
$Y_{\la i,j\ra})$ is definable from $\a_{\la i,j\ra})$ via
$\psi$. Then, if $\a <1$ defines the set $\{\a_{\la i,j\ra}: i,j \in
\NN\}$ via $\phi_\in$, 

\begin{equation} 
\label{lineY} Y=\{\la \x,\y \ra: \Rm \models \ex a' \ [\phi_\in(a',\a)
\ \land \ \psi(\x,\y;a')\}.
\end{equation}

Hence the complement $Z$ is definable from the parameters
$\d,\c$ and $\a$. Moreover 
$\a$ was obtained effectively in $\d,\c$ and a $\Sigma^0_{k+1}$ representation of $\Theta Z$.

To determine the sequence $(\a_{\la i,j\ra})$, let $A=W_i, B=W_j$ and $S=\Theta Y_{\la i,j\ra}$. We will construct $\a_{\la i,j\ra}$   uniformly in $i,j, \Theta Y$ and $\c,\d$. From a $\Pi^0_{k+1}$-representation of $\Theta Y$ one obtains a uniform sequence

$(R_n)$ of $\Sk$ relations such that

\begin{equation} \label{S} S= \bigcap_n R_n.\end{equation}

Since $Y_{\la i,j\ra}\sub \{\la \deg_m(A), \deg_m(B)\ra\}\cap [\d,\c]$
and $k \ge 3$, we may  suppose that

\begin{equation}
\label{relsn} R_n xy \RA \d \le deg_m(W_x),deg_m(W_y) \le \c \lland
W_x \equiv_m A \lland W_y \equiv_m B.
\end{equation}

We will define modified $\Sk$-relations $\tilde R_n$ such that
(\ref{S}) still holds, but also each $\tilde R_n$ is compatible with
$\equiv_m$. We cannot simply
 take the closure under $\equiv_m$, since it may be the case that $S=\ES$ because, for each $n$, $R_nxy$ holds with different indices $x,y$ for $A,B$. Instead, we reduce the relations $R_n$, making use of the fact that $S$ is $\equiv_m$-closed in order to maintain (\ref{S}). This takes three steps.

\be

\item First we find a uniform sequence $(R_n')$ of $\Sk$ relations

such that $R_n=R'_n$ if

\begin{equation}
\label{l5}
\fa u,v[ W_u \equiv_m A \lland W_v \equiv_m B \ \RA \ R_nuv],
\end{equation}

and otherwise $R_n'$ is finite. Note that (\ref{l5}) is a
$\Pi_1^{\ES^{(k)}}$  (= $\Pi_2^{\ES^{(k-1)}}$)-statement about $n$,
uniformly in indices for $A,B$.
 Thus we can effectively determine a set $X_n= W_u^{\ES^{(k-1)}}$ such that,
if (\ref{l5}) holds, then 
$X_n = \NN$, and $X_n$ is finite otherwise. Then $R_n'=R_n \cap
(X_n\times X_n)$ is a $\Sk$-relation as desired.

\item

We view each relation $R_n'$ as a relation c.e.\ in $\ES^{(k-1)}$.Thus
  we are effectively given  an enumeration $R'_n= \bigcup_s R'_{n,s}$,
 where the sequence $(R'_{n,s})_{s\in \NN}$ of strong indices for
 finite sets of pairs is recursive in $\ES^{(k-1)}$. We define
 $\ES^{(k-1)}$- recursive relations $R_{n}'' \sub R_{n}'$ as follows:
 at stage $s$, allow a pair $\la x,y \ra \in R'_{n,s}$
 into $R''_{n,s}$ only if $|R'_{m,s}| > n$ for all $m < n$. If (\ref{l5}) holds for each $R_n$ (and hence for each $R'_n$), then $\fa n R''_n = R'_n$, but otherwise $R''_n = \ES$ for almost all $n$.

\item Finally let $\tilde R_n$ be the closure of $R''_n$ under the equivalence relation $\equiv_m$.

\ee

To verify that $S= \bigcap_n \tilde R_n$, let $x,y \in \NN$ be arbitrary. First suppose that $Sxy$. Then, because $S$ is compatible with $\equiv_m$, (\ref{l5}) holds for each $n$. Hence
$\fa n\tilde R_nxy$.  Conversely, if  $\fa n\tilde R_nxy$, then 
$R''_n \neq \ES$ for each $n$. Then $R_n'$ is infinite, and therefore
(\ref{l5}) holds. Since $S= \bigcap_n R_n$, this implies that
$Sxy$.

Since (\ref{relsn}) holds for the sequence $(\tilde R_n)$ and $\tilde R_n$ is $\Sk$ uniformly in $n$, it is possible to apply the inductive hypothesis to each relation on $\Rm$ given by $\tilde R_n$. Thus we obtain a u.c.e.\ sequence $(\b_n)$ of incomplete $m$-degrees such that 

$$\la \x,\y \ra \in Y_{i,j} \ \LR \ \mbox{for each} \  n, \Rm \models \phi_{k,2}(\x,\y;\b_n).$$

Define $\{\b_n: n \in \NN\}$ from a parameter $\a_{i,j}< \one$.
Then

$$ Y_{i,j} = \{\la \x,\y \ra: \Rm \models \fa b [\phi_\in(b, \a_{\la i,j \ra }) \rightarrow \phi_{k,2}(\x,\y;b)]\}.$$

Since we have determined $\a_{\la i,j \ra} $ effectively, this
concludes  the proof of the lemma for an interval $[\d,\c]$, $\d \neq
\zer$.

 To reduce the general case to this, by Lemma \ref{Ershov}, (effectively in $\c$) obtain a minimal $m$-degree $\d$ such that
$\d \not \le \c$. Now we apply the above  with $\tilde \c= \c \vee \d$
instead of $\c$. Note that the intervals $[\zer,\c]$ and $[\d, \tilde
\c]$ are isomorphic: the isomorphism is $\x \ \rightarrow \ \x\vee \d$
and its inverse is $\y \ \rightarrow \ \y \wedge \c$. If $Z \sub
[\zer,\c]^N$ is $\Sk$ ($k \ge 3$), then so is $\tilde Z = \{\x\vee \d:\x \in Z\}$. Hence $\tilde Z$ is uniformly definable from the parameters $\d, \tilde \c$ and a parameter $\tilde \a$. Let $\a = p(\d, \tilde \a)$. Then $\a < \one$ and $Z = \{\x \wedge \c: \x \in \tilde Z\}$ is
uniformly definable from $\a$ and $\c$. (Note that, also in the
general case, $\a$ is obtained effectively.) \eop

The proof shows that in fact, for $k \ge 3$,  every $\Sk$ relation 
$Z\sub [\zer, \c)$ can be defined from parameters via a $\Sigma_{k+C}$
formula, for a fixed $C$. However, we only  obtain weak  uniform
definability, namely some extra  relations may be definable via our formulas.
 The following proposition shows that the  upper bound $\c< \one$ 
in the definability lemma\id{definability lemma}  is necessary.

\begin{proposition} There is a $\Sigma^0_4$-relation on $\Rm$ which is not definable from parameters.
\end{proposition}

\pf We use the fact that $\Rm $ has uncountably many automorphisms. First, by repeated applications of Lemma \ref{uce}, for each $m>0$ one can construct $\c$ such that
$|[\zer,\c]|=m$. Let

$$R= \{\la \a,\b \ra: \ex n >0( |[\zer,\a]| \ge n \lland \deg_m(W_n)=\b)\}.$$

Clearly $R$ is $\Sigma^0_4$. Assume that $R$ is definable from a parameter list. Because there are uncountably many automorphisms, there must be a non-identity automorphism $\Phi$ which fixes the parameter list. Then $\Phi$ 
respects $R$. We show that $\Phi(\b)=\b$ for each $\b$, a contradiction. Given $\b$, let $n > 0$ be minimal such that 
$\deg_m(W_n) = \b$. Then, for each $\a$, $R\a\b \LR |[\zer,\a]|\ge n$. Hence for each $\c$, $R\c\Phi(\b) \LR |[\zer,\c]|\ge n$. This implies that $\Phi(\b) = \deg_m(W_n) $. \eop

\section{Interpreting true arithmetic in $\protect\Th(\protect\Rm)$}

\begin{theorem} \label{ITARm}

$\ThN$ can be interpreted in $\Th(\Rm)$.

\end{theorem}

\pf We follow the framework of Section \ref{sectionITA}.
Fix any $\e <\one$ such that $[\zer,\e]$ is infinite. We carry out
(\ref{code}), representing numbers by 
 the degrees in $[\zer,\e]$ .

Let $h \le_T \ES^{(3)}$ be any map such that $\a_n=\deg_m(W_{h(n)})$
is a non-repeating list of all the degrees in $[\zer,\e]$ . Then addition and multiplication on $\{\a_n: n \in \NN\}$
(= $[\zer,\e]$), viewed as ternary relations, are $\Sigma^0_4$. By
the definability lemma, these relation can be defined from a list of parameters $\ol \p$, which includes $\e$.

To determine a \id{correctness condition} $\alpha_{st}(\ol p)$ in
(\ref{see}), consider an arbitrary parameter list $\ol \p$. First,
beyond the correctness condition $\alpha_0(\ol p)$ from Example \ref{codepa} we require
that $\e < \one$ and $\Mp \sub [\zer,\e)$. Let $k$ be the least number
such that, for each $\Mp$, the standard part is $\Sk$. Using the
definability lemma, we can now express in a first-order way that $\Mp$
is standard: we require that

\begin{quote} each subset
of  $\Mp$ defined from any parameters via the formula $\phi_{k,1}$
(and therefore, each $\Sk$-subset of $\Mp$) which contains $0^{\Mp}$
and  is closed under taking the successor function of $\Mp$ equals
$\Mp$. \eopnospace \end{quote}

\section{Model theoretic results on $\protect\Rm$}

\label{CharRm}

We survey several  results.

By the techniques of Denisov \citelab{Denisov:72}{2}, $\Rm$ possesses
continuum many automorphisms.
We apply the coding of copies of $\Nops$ to derive a uniform
definability result for the restrictions of automorphisms to proper
initial intervals. In particular, there are only countably many such
restrictions,  and the abundance of automorphisms stems  from the many
possibilities to put them together. 

\begin{theorem} \label{autos} The class of partial maps 
$$\{\Phi\lceil[\zer,\e]: \e < \one \lland \Phi \in \text{Aut}(\Rm)\}$$

is weakly uniformly definable.
\end{theorem}

\pf Suppose that $\e < \one$ and $\Phi \in \text{Aut}(\Rm)$. Then $\c=\e
\vee \Phi(\e) < \zer$ by the aforementioned  result of Lachlan
(see \cite[p.\ 45]{Soare:87}). Let  $$U=\{\la x,
\Phi(\x) \ra:
\x \le \e  \} .$$

Then $U \sub [\zer,\c]\times  [\zer,\c]$.   By the definability lemma, it is 
 therefore sufficient to show that, for some constant $k$ not
 depending on $U$, $U$ is $\Sk$. 

If $[\zer,\e]$ is finite then $U$ is $\SIII$. Now suppose
otherwise. As in  the proof  of Theorem \ref{ITARm}, fix a parameter list $\ol \p$ coding a copy $\M$ of $\Nops$ so that the
domain of $\M$ equals $[\zer,\e]$. If $\M'$ is the structure coded by
$\ol \p' = \Phi( \ol \p)$ via $S_M$, then the domain of $\M'$ is $[\zer,
\e']$ and $\M'$ is also a copy of $\Nops$. Clearly,

 $$ \la \x, \y \ra \in U \LLR \ex n \in \NN \  [ \x = n^\M \lland \y=
 n^{\M'}].$$

This shows that $U$ is $\Sk$ for some sufficiently large fixed $k$.   \eop

Next we give a coding without parameters of a copy of $\Nops$ and show
that the set of tops of finite initial intervals  is  definable. 

\begin{theorem} \label{NinRm} A copy of $\Nops$ can be coded in $\Rm$ without
  parameters. \end{theorem}

\pf The formula $\phi_\in$ from Theorem \ref{uce} determines a scheme $S_P$ to code subsets of 
$\Rm$  (with a vacuous correctness condition). We plan to
represent the number $n \in \NN$ by all sets $P$ such that $P
=n$. Thus let

$$\bfN = \{\a: |P_\a| < \infty \}.$$

To obtain a scheme as in (\ref{SMformulas}) but with an empty parameter
list, we have to give first-order definitions  without parameters of $\bfN$,
$\{\la \a, \b \ra:\a,\b \in \bfN \lland |P_\a = P_\b|\}$ and the
ternary relations on $\bfN$ corresponding to the arithmetical
operations $+,\times$.

The following formula determines an object scheme $S_C$ to code binary
relations:  $\phi(\x,\y;\b) \equiv$

$$ \ex z < 1 \ [ \phi_\in(z; \b) \lland \x = pr_1(z) \lland \y =
pr_2(z)].$$

(see Corollary \ref{pairs} for a definition of $pr_1, pr_2$.)
Clearly we can express in a first-order way that $C$ is a bijection
between sets defined coded via  some fixed schemes. 

For a first-order definition of $\bfN$, note that for $\a < \one$, $
\a \in \bfN$ iff there is a bijection between $P_\a$ and some initial
segment of a copy of $\Nops$ coded by the scheme of the preceding
section. By the results in Section \ref{prelims}, such a bijection can
be coded via $S_C$. 

Using elementary set theory in $\Rme$, we can also define in a
first-order way the  
other relations needed. Let $\phi_\equiv(x,y)$ be a formula expressing

$$ \ex C [C \ \text{is bijection} \ P_{x} \mapsto 
P_{y}]$$

and   $\phi_+( x,  y,  z)$ be a formula expressing
 $$\ex  u  \ex  v [\phi_\equiv( x,  u) \lland  \phi_\equiv( y,  v) 
\lland P_{ z}= P_{ u} \cup P_{ v} \lland P_{ u} \cap P_{ v} 
=\ES].$$

For  $\phi_\times(  x,   y,   z)$ we express in terms of definable
projection maps that $P_{  z}$ 
has the same size as the cartesian product $P_{  x} \times P_{  y}$. 
Thus $\phi_\times(  x,   y,   z)$ expresses

\begin{eqnarray*} \ex C_1 \ex C_2 
      && C_1: P_{  z} \mapsto P_{  x} \ttext{onto} \lland  C_2: P_{  z} \mapsto P_{  y} \ttext{onto} \lland \\
      && \fa a \in P_{  x} \ \fa b \in P_{  y} \ \ex!q \in P_{  z}
               [C_1(q) = a \lland C_2(q)=b].
\end{eqnarray*}   \eopnospace 

\begin{corollary}  $\{\b: [\zer, \b] \ttext{is finite} \}$ is
  definable in $\Rm$.\id{definability!in $R_m$}
\end{corollary}

\pf

 $[\zer, \b] \ttext{is finite} \LR\  \ex C \ \ex a \in \bfN [ C
\ttext{is bijection between} [\zer, \b] \ttext{and} P_a$]. \eop 

One can in fact obtain  a stronger result, using the proof  of Theorem
\ref{ITARm}: if $\CCC \sub \Rm$ is a class such that ``$\a \in \CCC$''
only depends on  the isomorphism type of  $[\zer, \a]$, then $\CCC$ is
definable iff $\CCC$ has an arithmetical index set. Thus a restricted
\ird{maximum definability property}\id{definability!maximum} holds (see Section \ref{Nonco}
below for a definition). The full {maximum
  definability property} in $\Rm$, which would state that a relation
is definable without parameters iff it is invariant under automorphisms and arithmetical,  is unknown. 
\vsp

Next we show the  existence of  an 
  incomplete $\e \in \Rm$ such that $[\zer,
\e)$ is an elementary submodel of $\Rme$ via inclusion. We use a
version of the 
 elementary chain principle. Write $\A \prec_{\Sigma_k} \B$ if $\A$
 is a submodel of $\B$ and the inclusion map is a
 $\Sigma_k$-elementary embedding. 
 
\begin{lemma}[\citelab{Chang.Keisler:73}{EC}] If $\AA_0
  \prec_{\Sigma_k} \AA_1\prec_{\Sigma_k} \ldots$ is a
  $\Sigma_k$-\ird{elementary chain} 
and  $\AA_\omega = \bigcup_{i \in \omega} \AA_i$, then $\AA_i \prec_{\Sigma_k}
\AA_\omega$ for each $i$. Moreover, if $\AA_i \prec_{\Sigma_k}
\B$ for each $i$, then $\AA_\omega \prec_{\Sigma_k}
 \B$. \eopnospace \end{lemma}

\begin{theorem}

$\fa  \a < \one \ex \e < \one \  [ \a\le \e \lland [\zer,\e) \ \prec \ \Rme]$.

\end{theorem}

\pf We use  the terminology and techniques of Denisov \citelab{Denisov:72}{EC} (see
also Odifreddi \cite{Odifreddi:nd}), which we review briefly. A main
concept is the notion
of an $L$-semilattice\id{L-semilattice} (called effective distributive upper semilattice
in \cite{Odifreddi:nd}), which is a type of distributive upper
semilattice with $0,1$.  Lachlan \citelab{Lachlan:70*1}{1} proved that up to isomorphism the
L-semilattices are the initial intervals of $\Rm$. 
 We also need  the following main tool for the
characterization of $\Rm$ from \cite{Denisov:72}. Enumerated
L-semilattices\id{L-semilattice!enumerated}  are L-semilattices with a presentation so that certain
effectivity conditions are satisfied. By the proof of  Lachlan's
characterization of initial intervals, each L-semilattice $[\zer, \x]$
is equipped with such an enumeration. Denisov's main technical result is
the following saturation property of $\Rm$.

\begin{quote} For enumerated L-semilattices $U_0, U$ and effective
  embeddings $g: U_0 \mapsto  \Rme, h: U_0 \mapsto U$ as initial
  intervals,  there is an effective embedding as an initial interval $f:
  U \mapsto \Rm$ such that $g = f \circ h$. 
\end{quote}

Moreover, the proof in \cite{Denisov:72} shows that an index for $f$ is obtained in
an effective way. 
Now, for each $\x$ we can effectively obtain $\y$ such that

\be

\item $\x <\one \RA \x<\y <\one$, and

\item $[\zer,\y] \cong \Rm $ via an (effective) isomorphism which acts
  as the identity on  $[\zer,\x]$.  

\ee 

 To see this, consider the (effective)  inclusion embedding $h$ of the enumerated
 L-semilattice $U_0=[\zer,\x]$ into the
enumerated  L-semilattice $L=\Rm \cup \{\t\}$, where $\t$ is a new largest
 element. 
By the above, obtain an effective $f: L \rightarrow \Rm$ which is the
 identity on $[\zer,\x]$, and obtain $\y$ as the image of $\one \in U$.
 Let us write $\y = F_0(\x)$. $F_0$ is an effective map on indices 
for c.e.\ $m$-degrees. Thus (like  the function $p$ introduced above)  $F_0(\x)$  really depends on the index via which $\x$ is given.
Iterating $F_0$  we obtain, by the effectivity of Denisov's construction, for any $\x < \one$ a u.c.e.\ chain

$$\x < F_0(\x) < F_0(F_0(\x)) < \ldots.$$

In a sense we will obtain $\e$ by iterating $F_0$ on $\a$
$\omega^\omega$ many times.  The construction bears some resemblance
to the \ird{reflection theorem}s from set theory.

Let $F_1(\x)$ be a degree $\y$ such that $[\zer,\y) = \bigcup_i[\zer, F_0^{(i)}(\x))$.
We can obtain $\y$ effectively in $\x$ by applying Theorem
\ref{Ershov}. Moreover, $\x< F_1(\x)$. More generally, if $F_k(\x)$
has been defined for all $\x$, $F_k$ is effective on indices and
$F_k(\x) > \x $ for $\x < \one$,
 let $F_{k+1}(\x) $ be a degree $\y$ such that $[\zer,\y) = \bigcup_i[\zer, F_k^{(i)}(\x))$. Then $F_{k+1}$ is a function on indices with the same properties.

\begin{claim} For $\x < \one$, $k \ge 0$, $[0,F_k(\x)) \prec_{\Sigma_k}\Rme$.\end{claim}

{\it Proof of the Claim.} By induction on $k$. For $k=0$, we assert  that $[\zer, F_0(\x))$ is embedded as an ordering into $\Rme$, which is correct. To prove the statement for $k+1$,
let $\z= F_{k+1}(\x), \z_j= F_k^{(j)}(\x) (j \ge 0)$. 
 By the inductive hypothesis, $[\zer, \z_j) \prec_{\Sigma_k} \Rme$, so the 
elementary chain principle implies that 

\begin{equation}
\label{zzz}  [\zer, \z) \prec_{\Sigma_k} \Rme \ttext{and} \fa j \
[\zer, \z_j) \prec_{\Sigma_k} [\zer, \z).
\end{equation}

Suppose $\b_0, \ldots, \b_{r-1} < \z$, and consider the formula 

$$\phi(\ol \b)=\ex \wt y \psi(\ol \b,\wt  y),$$

 where $\psi$ is a boolean combination of $\Sigma_k$- formulas and
 $\wt y$ is a tuple of variables of a certain length.
 We have to show that $$[\zer,\z) \models \phi(\ol \b) \LR \Rme \models \phi(\ol \b). $$

\be 

\item

First suppose that $[\zer,\z) \models \phi(\ol \b) $. Choose $j>0$ and
a tuple $\wt \c$ of elements in $[\zer, \z_j)$ such that $[\zer,\z)
\models \psi(\ol \b, \wt \c)$. 
 By (\ref{zzz}),  $[\zer,\z_j)  \models \psi(\ol \b, \wt \c)$. Then,
 because $[\zer, \z_j) \prec_{\Sigma_k} \Rme$, $\Rme \models \psi(\ol \b, \wt \c)$.

\item Now suppose that $\Rme \models \phi(\ol \b)$. Because $\Rme
  \cong [\zer, F_0(\z_j))$ via an isomorphism which acts as the
  identity on $[\zer, \z_j)$, there is a tuple of witnesses $\wt \c$
  in $[\zer, F_0(\z_j)) \sub [\zer,\z)$ such that $\Rme \models
  \psi(\ol \b, \wt \c)$. By
 (\ref{zzz}),  $[\zer,\z) \models \psi(\ol b, \ol c)$.

\ee   

Finally, let $\e > \a$ be such that $[\zer, \e) = \bigcup_{k\ge
  0}[\zer, F_k(\a))$.
 Since $[\zer, F_{l}(\a) ) \prec_{\Sigma_k} \Rme$  for all $l\ge k$,
 we conclude that  $[\zer, \e) \prec \Rme$ by the elementary chain principle. \eop

Notice that in fact $[\zer, \e) \cong \Rme$, because $[\zer, \e) $
satisfies the  characterization of $\Rme$ given in
\cite{Denisov:72}. However, the isomorphism cannot be $\Delta^0_3$
(let alone effective), because by construction of $\e$ we have a
u.c.e.\ chain $(F_k(\a))$ such that $\x < \e \LLR \ex k \ \x \le F_k(\a)$.
Such a chain converging to $\one$ cannot exist, because $\{i: W_i \equiv_m K\}$ is $\SIII$-complete.







\chapter{C.e.\ sets under inclusion}

\label{EEEEE}
\section{Outline}
\label{CEINC}

 We first  give a  proof  of Harrington's ideal definability
 lemma\id{definability lemma!lemma} explained in Section
 \ref{sectionITA}, (\ref{definability lemmaE}).  Based on this lemma 
we develop  a direct coding with parameters of a  standard model of arithmetic and
thereby give a  new   proof of Harrington's result  that true arithmetic can be interpreted in
$\hbox{Th}(\EE )$. 

Recall that $\EE^*$ is the lattice  of c.e.\ sets modulo finite differences. 
Both $\EE$ and $\EE^*$ are distributive lattices. The coding  methods can be used as well to give a
uniform coding of finite graphs\id{graph!directed} in $\EE^*$ via a $\Sigma_4$-scheme,
which proves the undecidability of $\Pi_6-\Th(\EE^*)$.
Furthermore they  yield 
  elementary differences between
relativized versions of $\EE$. A natural question due to E.\ Herrmann  is  if, for
$0<p<q$, the relativization\id{relativization!of $E$} of $\EE$ to $\emptyset^{(p-1)}$ (i.e.\ the
$\Sigma^0_p$-sets under inclusion)  and to $\emptyset^{(q-1)}$ are elementarily
equivalent\id{elementarily
equivalent}.  Evidence for an affirmative answer would come  from the fact
that constructions of c.e.\ sets which show that $\EE$ possesses certain
first-order properties, like the construction of a maximal set in 
Friedberg \citelab{Friedberg:58}{1},  relativize and therefore show that for each $Z\subseteq\NN$, $\EE^Z$, the
lattice of sets c.e.\ in $Z$, has the same
property. 
However, we answer the  question negatively.  Roughly speaking, an elementary
difference between the lattice of $\Sigma^0_p$- and the lattice of $\Sigma^0_q$-sets
$(0<p<q)$ is obtained by considering the ``coding power''  in the
structure of a 
scheme of formulas intended to code models of  $\mbox{PA}^-$ with an extra unary
predicate. This coding power  increases with the complexity of
the oracle $\EE$ is relativized to.  
 
Recall that $\LL^*(A)$\id{L\@$\LL^*(A)$} is the lattice of c.e.\ supersets of $A$ modulo
finite differences and that $A$ is  {\it
  quasimaximal}\id{quasimaximal} if $\LL^*(A)$ is
finite or, equivalently, if $A$ is the intersection of finitely many
maximal sets.  In Soare \cite{Soare:87} it is asked if the class of quasimaximal
sets is definable in $\EE$.  We answer this question affirmatively.
The definability of ``quasimaximal'' and of further classes 
of hh-simple sets can be obtained from the ideal definability
lemma\id{definability lemma!ideal}  and 
certain isomorphism properties of boolean algebras which are coded in 
$\EE$ with parameters.

The lattice $\EE$ is set apart from other structures studied in computability
theory by the fact that many results restricting coding and definability can be
obtained.  We show that no infinite linear order can be coded
(without parameters) even in the most general way, namely on
equivalence classes of $n$-tuples.  Moreover we give an example of a subclass
of $\EE$ which is nondefinable, but has an arithmetical index set and
is invariant under automorphisms.

For any class $\CCC\sub \EE$, $\CCC^*$ will denote the class
$\CCC/_{\ns =^*}$. We 
state our results for $\EE$ instead of $\EE^*$ mostly for notational
convenience.
 For 
definability and coding concerns, it does not matter whether the
setting of $\EE$ or of $\EE^*$ is 
used, unless we   study   fragments of the theory. The reason is
that from the methods in Lachlan \citelab{Lachlan:68*3}{1} one can derive that,
if $C\subseteq\EE^n$ is closed under finite variants, then

\begin{equation} \label{LachC} \CCC\ \mbox{
definable in}\  \EE\ \LLR  \CCC^*\ \mbox{definable in} \  \EE^*, \end{equation}

via a uniform translation between formulas,  and
similarly for definability with parameters.  Now our coding and
definability results do not refer to membership of particular elements. 
So one can easily transfer all the results from $\EE$ to $\EE^*$,
e.g.\ one can  prove that
$\{ A^*:\LL^*(A)\text{ finite}\}$ is definable in $\EE^*$ or that the
$\Sigma^0_2$-sets modulo finite variants are not elementarily
equivalent to $\EE^*$. 

Intervals play an important role in the study  of the lattice $\EE$. Several
interesting properties of a c.e.\ set can be given alternative 
definitions  in terms of  the structure of $\LL(A)$, the lattice of
c.e. supersets of $A$. For instance,  a coinfinite
c.e.\
set $A$ is \iid{hyperhypersimple} iff
 $\LL(A)$ is  a boolean algebra, and $A$ is {\it $r$-maximal}\id{r-maximal}
if and only if  $\LL(A)$ has no nontrivial complemented elements. 

Unlike to the case of $\RR_m$, the possible structure of intervals of
$\EE$ and $\EE^*$ is still not very well
understood. 
 Lachlan \citelab{Lachlan:68*3}{2} shows that the boolean algebras which can be
represented as $\LL^*$, $A$
hh-simple, are precisely the $\SIII$--boolean algebras (see Section
\ref{CEB} for a definition). 
The
class
of
$r$--maximal
sets is much more
elusive.
  Cholak
and Nies \citelab{Cholak.Nies:nd}{1} have  shown  that infinitely many non-isomorphic lattices
$\LL^*(A)$, $A$ $r$--maximal\id{r-maximal}, exist. 

We now review the notation and terminology used in this chapter. All
subsets of $\NN$ are c.e.\ unless otherwise mentioned. 

\begin{notation} \label{Enot}  {\rm
\begin{itemize}
\item Capital letters $A, B, C, X, Y$ range
over r.e.\ sets, letters $R,S,T$ over computable sets. 

\item  $X\sqsubset A\Leftrightarrow (\exists Y)[X\cap
  Y=\emptyset\lland X\cup Y=A]$, 

$\BB (A) =\{ X:X\sqsubset A\}$ and 

$\RR (A)=\{ R:R\sqsubset A\}$. 

\item An ideal
$I$ of $\BB (A)$ is $k$-acceptable\id{acceptable!$k$-acceptable} if $\RR (A)\subseteq I$ and $\{
e:W_e\in I\}$ is $\Sigma^0_k$.  If we say ``$I$ is
acceptable''\id{acceptable}  we mean
that $I$ is $k$-acceptable, where $k$ is a fixed number which depends only
on the context in which $I$ is defined (e.g.\ on formulas in some coding
scheme or on arithmetical constructions). 
\item
Given an r.e.\ set $A$ define a $\Delta^0_3$-enumeration
$(U_e )_{e\in\NN}$ of $\BB (A)$ as follows:

\noindent  if $e=\langle i,j\rangle$, $W_i\cap W_j=\emptyset$ and  $W_i\cup
W_j=A$ let $U_e =W_i$ and write $\overline U_e$ for
$W_j$.  Else let $U_e =\emptyset$ and $\overline U_e=A$.

\end{itemize}
}
\end{notation}
Recall that the major subset relation is defined as follows: for $A,B \in \EE$, 

$$  B\subset_m A \LR B\subset_\infty
A\wedge(\forall W\ \mbox{c.e.})[ A\cup W=\NN \Rightarrow B\cup
W=^*\NN].$$

\noindent A set  $B$ is a {\it small  subset}\id{major subset!small} of $A$ , denoted $B
\subset_s A$,
 if  $B \subseteq A$ and
\begin{equation} \label{small subset}
(\forall U,V)[  U \cap (A - B) \subseteq^* V \Rightarrow (U - A) \cup
V \hbox {r.e.}]. \end{equation}
We will make use of the  following well-known facts.

\begin{lemma} \label{sm}
   
\be
\item[(i)]
  If $B \subset_s A$, then each $Y
\sqsubset A$  such that  $Y \subseteq^* B$ must be computable.

\item[(ii)]  If  $B 
\subset_m A$, then for each computable  $R \sub A, R \sub^* B$.

\item[(iii)]  If $B \subset_s A$ and  $B 
\subset_m A$ (this is also denoted by $B\subset_{sm} A$)
and the set $X \sqsubset A$ is non-computable, then $X - B$ is
non-c.e.

\ee

\end{lemma}

\pf

\noindent (i).  Let $U = \NN, V = A - Y$.  Then $U \cap (A - B) =^* A - B 
\subseteq^* V$, so $\overline A \cup V = \overline Y$ is c.e.

(ii). Immediate because $A \cup (\NN-R) = \NN$ and $\NN-R$ is c.e. 

(iii).  If $X - B$ is c.e., then $Y: = X 
\cap B
\sqsubset A$, because $A - (X \cap B) = (A - X) \cup (X - B)$.  So by (i),
$X \cap B$ is computable.  Since $X$ is non-computable, $X - B$ is non-computable,
so we can
choose an infinite computable $R \subseteq X - B$.  This contradicts $B  \subset_m A$.
\eop

\section{The ideal definability lemma}

\label{ideal definability lemmaE} \label{IDLE}\id{definability lemma!ideal}
\begin{Idl}  For each $n\ge 1$ the class of
  $2n+1$-acceptable\id{acceptable} ideals is  uniformly
  definable\id{uniformly definable}. More precisely,   there is a formula with parameters
$\varphi_n(X:D,\overline C,A)$ $(|\overline C|=n)$ with the following property.

\begin{quote}  If $A$ is non-computable, for
 $D,\overline C$ ranging over tuples where the  correctness condition
 
\bc $D \sub
C_0  \sub C_1 \sub \ldots \sub C_{n-1} \sub A \lland D\subset_m C_1$ \ec

$(D \subset_m A \mmbox{in the case that} n=1)$ is satisfied,

 \bc $\{
X:\EE\models\varphi_n(X;D,\overline C,A)\}$ \ec

 ranges precisely
over the class of 
$2n+1$-acceptable ideals of $\BB(A)$. \end{quote}
\end{Idl}

\pf  The formulas $\varphi_n$ are defined recursively, by reducing the problem of defining
a $2n+3$-acceptable ideal to the problem to define a
$2n+1$-acceptable one.   

\vsp

{\it The Case  $n=1$.} Let
\begin{equation}
\label{phi1}
\varphi_1(X;D,C,A)\equiv X\sqsubset A \lland   X\cap  C\sub^* D].
\end{equation}

Clearly, the index set of any ideal $I$ defined via $\phi$ is a
$\SIII$-ideal of $\BB(A)$. Moreover, since $D\subset_m A$ is a
\id{correctness condition}, $\RR(A)
\sub I$ by Lemma \ref{sm} (ii). We now prove that, whenever
$D\subset_{sm}A$, then each 3-acceptable ideal of $\BB(A)$ has the
form $\{X: X\cap C \sub^* D\}$ for some $C$. To do so, we will in fact
prove a slightly more general fact about intervals $[D,A]$, where $D
\subset_m A$, which will be used again in Section \ref{IntE}. Consider
the set 

\begin{equation} \label{Splboolean algebra}\BB=\{X\cup  D^*: X \sq A \} \end{equation}

(we will write $X\cup Y^*$ instead of $(X \cup Y)^*$). By the reduction
principle (\ref{RedPr}),
 $\BB$ equals the set of
complemented elements in the lattice  $[D^*,A^*]$ and therefore is a
boolean algebra. Then $(U_e \cup D^*)_{e \in \NN}$ is a $\Delta^0_3$ listing of
$\BB$, (see Notation \ref{Enot} for the sequence $(U_e))$, and we obtain a notion of
index sets of subsets of 
$\BB$ with respect
to that listing, and
especially of $\SIII$ subsets of $\BB$. {\it In the following we will
identify  subsets of $\BB$ with their index sets.}

\begin{lemma} \label{Fidel} If $D \subset_m A$ 
and $I$ is a $\SIII$--ideal of $\BB$,
  then there is  $C$, $D \sub^* C \sub^* A$ such that

 \begin{equation}
   \label{fi}
   I = \{U_j \cup D^*: U_j \cap  C \sub^*D\}.
 \end{equation}

\end{lemma}

\pf First we give an effective representation of the
filter of complements of elements of $I$, using the following
uniformization
fact. 

\begin{fact} If $(W_{g(i)})_{i \in \NN}$  is a sequence of splits of $A$, $g
\le_T \ES''$, then there is a uniformly c.e.\ sequence of splits 
$(Z_i)$ of $A$  such that $\fa i \ W_{g(i)} \triangle Z_i \sub^* D$.
\end{fact}

\n To  prove this, choose a u.c.e\ sequence $(V_k)$ of initial segments of 
$\NN$ such that $W_p =  W_{g(i)} \LR \exists n V_{\langle i,p,n\rangle } =
\NN$ (this is possible since ``$W_p =  W_{g(i)}$'' is $\SIII$).
The desired u.c.e.\ sequence is 

$$\begin{array}{l}
Z_i=   \{a: \ex s \ex q=\langle i,p,n\rangle \\
 \ \mbox{max} \ \bigcup_{ \langle i,p',n' \rangle < q}
 V_{\langle i,p',n' \rangle,s}  <  a \le \ \mbox{max} \ V_{q,s} \lland
 a \in W_{p,s} \}.
\end{array}
$$

\n Given $i$, let $p=g(i)$ and let $q = \langle p,n\rangle $ be the least such that $V_{\langle i,p,n\rangle } =
\NN$. Then 
$Z_i=^* R \cup W_p$, where $R$ is the computable set $\{a: \ex s \ a \in Z_{i,s}
\lland a >  \mbox{max} \  V_{q,s} \}.$
Therefore 
$W_p \triangle Z_i \sub^* R \sub^* D$.
This proves the fact. 

Clearly, the indices of c.e.\ sets which are complements of elements
in $I$,

$$S = \{i: \ex k \in I \  W_i \cap U_k  \sub^* D \lland W_i \cup
U_k =^* A\}$$

 is $\SIII$ and therefore $S$ is the range of a function $g
\le_T \ES''$. Applying the preceding fact we obtain a sequence $(\tilde
Z_i)$. Let $Z_n = \bigcap_{i \le n} \tilde Z_i$. Then the  u.c.e.\ sequence $(Z_n
\cup D)^*_{n \in \NN}$ generates the filter of complements of elements in $I$. 

To build $C$, we 
meet for each $n$ the following requirement:

$$P_n : |W_e \cap Z_n \cap \ol{D}|= \infty \RA |W_e \cap C  \cap
\ol{D}| \ge k \ ( k = \langle e,n\rangle ).$$

\n The construction of $C$ is the following. Let $C_0 =\ES$. At a
stage $s+1$, for each $\langle e,k\rangle  =n <s$, act as follows. If $P_n$ is
currently unsatisfied, namely  $|W_{e,s} \cap C_s  \cap
\ol{D_s}| < k$, and there is an $x \in Z_{n,s} -D_s$ such that $x \in
W_{e,s}$, then enumerate the least such $x$ into $C$. 

We verify that $C$ satisfies (\ref{fi}). Notice that at
most $k+1$ elements which are permanently in $\ol{D}$ are enumerated into $C$ for the sake of
$P_{\langle e,k\rangle }$. Therefore $C \sub^* D \cup Z_m$ for each $m$. Now, if $j
\in I$, then choose an $m$ such that $Z_m - (A - U_j)
\sub^*D$, i.e.\  $Z_m \cap U_j \sub^*D$.
  Since $C \sub^* D \cup Z_m$, $C \cap U_j \sub^* D$.

If $j \not \in I$, then $D \cup (A - U_j)^*$ is not in the filter
  dual to $I$, so  $D \cup (A - U_j)$  does not $*$-include $Z_n$ for any
  $n$. Thus if $W_e=U_j$,  the hypothesis of all the requirements $P_{\langle e,k\rangle }$
  is satisfied. Hence $C \cap W_e \cap \ol{D}$ is infinite.  This
  proves Lemma \ref{Fidel}. \eop

Now assume that
$D\subset_{sm}A$ and let $I$ be a  3-acceptable ideal of $\BB(A)$. To
show that $I$  has the
form $\{X: X\cap C \sub^* D\}$ for some $C$, consider the ideal $\hat I=
\{(X \cup D)^*: X \in I\}$  of $\BB$.  Then  $D \subset_s
A$ implies $X \in I \LLR D \cup X^* \in
\hat I$: the direction ``$\RA$'' is immediate, and ``$\LA$'' follows
because $D \cup X=^* D\cup Y$ for $Y \in I$ implies that $X \trup
Y \sub^*D$. Hence $X \trup Y $ is computable and $X \in I$.
Since $\hat I$ is $\SIII$, we obtain $C$ such that $\hat I = \{ X
\cup D^*: X \cap D \sub^* C\}$. So $I = \{X: X \cap D \sub^* C \}$.

\vsp

{\it The Inductive Step.}  To complete the 
proof of the ideal definability lemma, 
  we will  show the following: if $m 
\geq 2$ and $I$ is
an $m + 3$-acceptable ideal of $\BB (A)$, then there is a
 non-computable $C \subseteq A$  and an $m + 1$-acceptable ideal $J$ of $\BB 
(C)$ such that, for each $X \sqsubset A$,

\begin{equation} 
\noindent   X \in I \ \Leftrightarrow \ \exists R
\subseteq A\ \forall S \subseteq A - R \ [ X \cap S \cap R \in J] . 
\label{S2} \end{equation}

Then,  if $\overline C=(C_0,\ldots,C_{n-1})$, let
\begin{eqnarray}
\label{phin}
&\varphi_{n+1}(X;D,\overline C, C_n,A)\equiv  & X\sqsubset
A\lland\\
& & \exists R\subseteq A \ \forall S\subseteq A-R\nonumber \\
&  & \varphi_n(X\cap S\cap C_n;D,\overline C, C_n).\nonumber
\end{eqnarray}
 (Recall that the variables $R$,$S$ range over computable sets. Notice
 that $C_n$ plays the role of $C$ in (\ref{S2}).)   For
 instance,

\begin{eqnarray*} 
\varphi_{2}(X;D,C_0, C_1,A) &\equiv&  X\sqsubset
A\lland\\
&& \exists R\sub A \  \fa S \sub A \  [X \cap S \cap C_1 \cap  D \sub^* C_0].
\end{eqnarray*}

We first check that this formula only defines $2n+3$-acceptable ideals.
Firstly,  if $X\sqsubset A$ is 
computable, then (\ref{phin}) 
holds via $R=X$.  Secondly, the class of $X$ satisfying $\varphi_n$
is downward closed, and if $X,Y$ satisfy $\varphi_{n+1}$ via
$R_X$ and $R_Y$ respectively, then $X\cup Y$ satisfies
$\varphi_{n+1}$ via $R_X\cup R_Y$, by inductive
hypothesis on $\varphi_n$.  Finally, to see that
the ideals defined by $\varphi_{n+1}$ are 
$\Sigma^0_{2n+3}$, we write  $\varphi_{n+1}(X;D,\overline C,
C_n,A)$ more
explicitly (for the moment, let $R,S$ range over arbitrary c.e.\ sets):

$$\aligned
&(\exists R\subseteq A)(\exists\tilde R)[R\cap\tilde
R=\emptyset\ \land R\cup\tilde R=\NN\ \land\\
&\quad (\forall S\subseteq A\cap\tilde R)[S
\ttext{noncomputable}(\Pi^0_3)\lor\\
&\qquad \varphi_n(X\cap S\cap C_n; D,\overline
C,C_n)(\Sigma^0_{2n+1})].
\endaligned
$$
Because $n\ge 1$, this shows that the corresponding index set is
$\Sigma^0_{2n+3}$.

To prove (\ref{S2}), we need some facts. First, we describe an
appropriate set $C$ for (\ref{S2}). Each noncomputable splitting $X$ of
$A$ effectively obtains a \iid{trace}  $T \cap C$ in $\BB(C)$, where
$T\sub X$ is computable.

\begin{specialtrace}  \label{Etrace} Let $A$ be non-computable.  Then there is $C \subseteq A$
such that  $(\forall X \sqsubset A$ non-computable) $(\exists T \subseteq X$ computable)
$$
  T \cap C \  \hbox {non-computable}.
$$
A strictly increasing (finite or infinite) computable sequence $b_0 < b_1 < \ldots$ such
that  $T = \{ b_0,b_1, \ldots  \}$  can be obtained effectively in (an index 
for) $X$.  We write $T = T_{X}$.

\end{specialtrace}

\pf   Let $B \subset_{sm} A$.  By an infinitary version of the 
proof of the
Friedberg splitting theorem in Soare \cite{Soare:87}\id{splitting!in
  computability theory}
\label{Friedberg:58.2}, obtain a u.c.e. partition $(B_k)$ of 
$B$ such that
$$
      (\forall W)(\forall k) [W - B \  \hbox {non-computable}\ 
\Rightarrow W - B_k \ \hbox {non-computable}]. \label{S1}
$$
Let $C = \bigcup \{ B_n : n\in K\}$.  We claim that $C$ is the desired set.
First we show that  for each $k$
and each non-computable $X \sqsubset A$, $X  \cap B_k$ is
infinite.
By Lemma \ref{sm} (iii), $X - B$ is non-c.e.  So, by (\ref{S1}), $X   - 
B_k$ is non-c.e., 
 thus  $X \cap B_k$ must be infinite.

Now define $T_{X} = \{ b_0, b_1 , \ldots \}$, where $(b_k)$ is an
effective strictly increasing sequence and $b_k \in X \cap
B_k$.  To do so, by induction over $k$, enumerate $X  \cap 
B_k$ until a new element
is found.  If $X$ is non-computable, then $T_{X}$ will be an infinite
computable subset of $X$.  Moreover, $T_{X} \cap C
\equiv_m K$, so $T_{X} \cap C$ is non-computable.
\eop

We now give a lemma on how to approximate $\Sigma^0_3$ sets. This
lemma will be relativized to $\ES^{(m)}$ in order to obtain (\ref{S2}).

\begin{lemma} \label{approxS3}  If $P$ is a $\Sigma^0_3$ set, then
  there is a u.c.e.\ 
sequence $(Z_i)$ such that $Z_i \subseteq \{ 0,\ldots ,i\}$
and
\bi

\item[a)]  $(\forall b \in P)$ (a.e.i~) $[b\in Z_i]$

\item[b)]  $(\exists^\infty i) \; [Z_i \subseteq P]$

\ei

\end{lemma}

 {\it Remark.}  If $b \not \in P$, then $b \not \in Z_i$ infinitely many
$i$, so 
$$b \in P \Leftrightarrow \ \mbox{for almost every}\  i \ b\in Z_i.$$  Note that the right
hand side is in $\Sigma^0_3$-form.

\pf  We first assume that $P$ is a $\Sigma^0_2$  and show that
there exist a sequence $(Y_i)$ of strong indices for finite sets with the
properties required in the lemma (this was first proved by Jockusch). 
 For  the general case, we will 
relativize  to $\emptyset '$.

  If $P$ is  $\Sigma^0_2$, there is a c.e. 
set $C$
such that $P = \{ (X)_0 : x \in \overline C\}$.  Suppose $C = \bigcup C_i$ ,
where $(C_i)$ is an effective  sequence of strong indices and $C_i \subseteq \{ 0,\ldots
,i\}$.  Let
$d (i) = \min(C_{i+1} - C_i) \  \hbox {if}\quad C_{i + 1} - C_i \not = 
\emptyset$ and 
$d(i)= i + 1 \  \hbox {else}$. 
Note that at most two arguments for the map $d$ can yield the same value.  Let 
$Z_i$ be a strong index for
$$
\{ c < d (i): c\not \in C_i \} .
$$
Then $a \not \in C \Rightarrow a \in Z_i$ for almost every $i$ and, if $j$ is a
non-deficiency state, i.e.
$$
d(j) = \mbox{min} \{ d(i) : i \geq j\},
$$
then $Z_j \subseteq \overline C$.  Now let
$$
Y_i = \{ (x)_0 : x\in Z_i\} .
$$
Then $Y_i \subseteq \{ 0 , \ldots i \}$ (because this  holds for $Z_i$).  For (a),
if $b \in P$, say $b = (c)_0$ for $c \in \overline C$, then $b \in Y_1$
whenever $d_i > c$, so almost every $i$ $b\in Y_i$.  For (b), note that $Y_i \subseteq
P$ whenever $Z_i \subseteq \overline C$.

Now suppose $P$ is $\Sigma^0_3$.  By relativization to $\emptyset '$, 
there is a
$\Delta^0_2$-sequence of strong indices for finite sets $Y_i \subseteq \{
0,\ldots ,i\}$ such that (a) and (b) hold.  By the Limit Lemma \cite{Soare:87}, there
is a computable array of strong indices $(Y_{i,k})$
such that for each $i$ and for almost every k, $Y_{i,k} = Y_i$.  Let
$$
Y^*_{\langle i,k\rangle} = \{ 0,\ldots ,i\} \cap \bigcup_{t\geq k} Y_{i,t} ,
$$
and let $f$ be a $1 - 1$ computable function such that
$$
rg(f) = \{ \langle i,k \rangle : k = 0 \vee Y_{i,k} \not = Y_{i,k-1} \} .
$$
Note that, for each $i$, there are only finitely many $j$ such that
$(f(j))_0 = i$.  We claim that $Z_j = Y^*_{f(j)} \cap \{ 0,\ldots ,j\}$ is the desired u.c.e.
sequence.  For (a), if $b\in P$,
then for almost every $i$, $b \in Y_i$. Since $Y_i = Y^*_{i,k}$ for almost every
$k$, by the above property of $f$, $b \in Z_j$ for almost every $j$.

For (b), if $i$ is such that $Y_i \subseteq B$ and $s$ is maximal such that
$s = 0$ or $Y_{i,s-1} \not = Y_{i,s}$, then, for $j$ such that $f(j) =
\langle i,s\rangle$, $Z_j = Y_i \subseteq B$.  So for infinitely many $j$
$Z_j \subseteq B$.   \eop

We are now ready to prove (\ref{S2}).
Let $P = \{ e : U_e \in I\}$ (recall $(U_e)$ is a listing  of the splits
of $A$).  By applying Lemma \ref{approxS3} relativized to $\emptyset^{(m)}$, we obtain a
sequence of sets $(Z_i)$ which are uniformly $\Sigma^0_{m+1}$ such that $Z_i
\subseteq \{ 0,\ldots, i\}$ and
$$
U_e \in I \Leftrightarrow (a.e. \  i) (e \in Z_i] \leqno (a')
$$
$$
 \exists^\infty i\ \forall e \in Z_i\  [U_e \in I] . \leqno (b')
$$
Let $C \subseteq A$ be the set obtained by the Trace Lemma \ref{Etrace}.  Moreover let $\BB
(A)_{\leq i}$ be the boolean algebra generated by $\{ U_e : e \leq i\}$
(assume $B (A)_{\leq 0} = \{ \emptyset ,A\}$).

\begin{claim}  There is a $\Delta^0_3$-sequence $(S_i)_{i\in\NN}$ of 
computable subsets of $A$ 
such that the sets $S_i$ are pairwise disjoint,

$$
 \forall R \subseteq A \ \exists i [R \subseteq S_0 \cup \ldots \cup S_i]
$$
and
$$
\forall i \ \forall V \in \BB (A)_{\leq i}  [ V \ \hbox {non
computable } \  \Rightarrow V \cap S_i \cap C\ \hbox {non-computable}].
$$
\end{claim}

Then we will define $J$ essentially as the ideal on $\BB (C)$ generated 
by the intersections $U_e \cap S_i \cap C$, where $e \in Z_i$.
Let $(R_i) $ be a $\Delta^0_3$ listing of $\RR(A)$.

 {\it Proof of the Claim.} 
Let $S_0 = \emptyset$ and, if $\hat S_{i} = S_0 \cup \ldots \cup S_i$,
$$
     S_{i+1}  =  (R_i - \hat S_i ) \cup  \bigcup \{ T_{V-\hat S_i} : V \in \BB (A)_{\leq i + 1} \}
$$
Then $R_i \subseteq S_0 \cup \ldots \cup S_{i+1}$.  Moreover, if $V \in \BB
(A)_{\leq i + 1}$ is non-computable, then, by the Trace Lemma \ref{Etrace}, $T_{V-\hat S_i} \cap C$
is non-computable (where  $T_{V-\hat S_i} \subseteq V$), so, since $S_{i+1}$
computably splits into  computable sets $T_{V-\hat S_i}$ and $S_{i+1} -  T_{V-\hat 
S_i}$, $V \cap S_{i+1} \cap C$ must be non-computable.

  Let $$J\ =  \text{ the
ideal of} \  \BB (C) \ttext{generated by} \RR(C) \ttext{and} 
\{ U_e \cap S_i \cap C: e \in Z_i\} .
$$

Since the relation ``$e \in Z_i$'' is $\Sigma^0_{m+1}$, $(S_i)$ is a
$\Delta^0_3$ sequence and $m\geq 2$, $J$ is a
$\Sigma^0_{m+1}$-acceptable ideal.  It remains to verify (\ref{S2}). 
  Suppose $U = U_{\tilde e}$.

{``$\Rightarrow$''}  If $U_{\tilde e} \in I$, choose $i_0$ such that
$e \in Z_i$ for all $i > i_0$.  We claim that $R = S_0 \cup \ldots \cup
S_{i_0}$ is a witness for the right hand side in \ref{S2}.  If
$S \subseteq A - R$, then $S \subseteq S_{i_0 + 1} \cup \ldots \cup S_j$ for
some $j > i_0$.  Now $U_{\tilde e} \cap S_i \cap C \in J$ for any
$i >i_0$, so, $U_{\tilde e} \cap S \cap C \in J $.

{``$\Leftarrow$''}  Suppose $U_{\tilde e} \not \in I$.  Given any
$R\subseteq A$, choose $k$ such that $R \subseteq S_0 \cup \ldots \cup S_k$.
By ($b'$), there is an $i > k$ such that $Z_i \subseteq \{ e : U_e \in 
I\}$,
and also $U_{\tilde e} \in \BB (A)_{\leq i}$.  We show that $S_i$ is a
counterexample to the right hand side in \ref{S2}, i.e.
$U_{\tilde e } \cap S_i \cap C \not\in J$.  Let $V = U_{\tilde e} -
\bigcup_{e\in Z_i} U_e$.  Since $U_{\tilde e} \not\in I$, $V$ is a non-computable
element of $\BB (A)_{\leq i}$.  So $V\cap S_i \cap C$ is not computable by the
claim above. But,  if $U_{\tilde e} \cap S_i \cap C \in J$, then, by
 the disjointness of the sets $(S_j)$,
$$
U_{\tilde e} \cap S_i \cap C
\subseteq R \cup (\bigcup_{ e\in Z_i} U_e \cap S_i
\cap C)
$$
for some computable subset $R$ of $C$.  So $V \cap S_i \cap C$  is 
computable as a split of C which is contained in a computable subset of C.

This concludes the proof of the ideal definability lemma. \eop

\section{Defining  classes of hh-simple sets}
 
\label{Defhhsim}
Recall that $A$ is \iid{hyperhypersimple} ({\it hh-simple}) if $\LL^*(A)$
forms a boolean algebra. In that case, $\LL(A) = \{A \cup R: R \
\text{computable}\}$.  We  consider parameterless \ird{definability}  in
$\EE^*$ of   classes of
hh-simple sets,  based 
on the ideal definability lemma.
For instance  the class of \ird{quasimaximal}
sets is definable in $\EE$. Recall that, by (\ref{LachC}), we can
disregard the difference between $\EE$ and $\EE^*$. We need two facts.
 
\begin{fact}\label{fact1}  If $\LL^*(A)$ is a boolean algebra, then there is a $\Delta^0_3$--
isomorphism $\Phi:\LL^*(A)\mapsto \BB_{rec}^*(A)$, where
$\BB_{rec}^*(A)=\{(R\cap A):R\hbox{ computable}\}$.\end{fact}

\pf  Let $\Phi(B^*)=(R\cap A)^*$, where
$B=A\cup R$.
 Note that it takes an oracle $\emptyset''$ to find $R$ from an input $e$ such
 that  $B=W_e$. \eop
 
Observe that $\BB_{rec}^*(A)$ is a subalgebra of $\BB^*(A)$ containing $\RR^*(A)$.
Thus we also obtain an isomorphism of the lattice of $\Sigma_k^0$-ideals 
$I$ of
$\LL^*(A)$ $(k\ge3)$ onto the lattice of $\Sigma_k^0$-ideals $\tilde{I}$ of
$\BB_{rec}^*(A)$ which contain $\RR^*(A)$.  The ideal definability lemma now 
implies that the
$\Sigma_k^0$-ideals of $\LL^*(A)$ ($k\ge3$ odd) are uniformly definable,
because $\tilde{I}=[\tilde{I}]_{\text{id}}\cap\BB_{rec}^*(A)*$, where
$[\tilde{I}]_{\text{id}}$ is the ($k$-acceptable) ideal of $\BB^*(A)$ generated
by $\tilde{I}$. 

Fix a $hh$-simple $A$ as a parameter. We consider definability 
of ideals of $\LL^*(A)$ with 
parameter $A^*$ in   $\EE^*$.  The  {\it derivative}\id{boolean
  algebra!derivative of} of a boolean algebra $\B$ is 
$\B/_{\ns U}$, where $U$ is the ideal generated by the atoms of $\B$.
 If $I$ is an ideal of $\LL^*(A)$, let $\AAA(I)$ be the
ideal
of $\LL^*(A)$ generated by the atoms of $\LL^*(A)/I$ (i.e\., 
$\LL^*(A)/\AAA(I)$ is the {derivative} of $\LL^*(A)/I$).

\begin{fact} \label{fact2} If $I$ is an ideal of $\LL^*(A)$ which  is
  definable in $\EE^*$ with parameter $\A^*$,
then so is $\AAA(I)$.  The formula defining $\AAA(I)$ only depends on the
formula defining $I$, not on $A$.\end{fact}

\pf  If $I$ in a $\Sigma_k^0$ ideal $(k\ge 3)$, then $\AAA(I)$ 
 is 
$\Sigma_{k+2}^0$.  So we can define $\AAA(I)$ as the least $\Sigma_{k+2}^0$
ideal of $\LL^*(A)$ which contains all the elements of $I$ and all $B^*\ge A^*$ such that
$B^*/I$ is an atom in $\LL^*(A)/I$. \eop

Note that we can also express that  $\AAA(I)$ contains infinitely many atoms of
$\LL^*(A)/I$:  this is the case iff $\AAA(I)$ describes a nonprincipal ideal
in $\LL^*(A)/I$, i.e.\   if there is no $B^*\ge A^*$ such that, for each 
$C^*\supseteq A^*$, $C^*\in\AAA(I)\Leftrightarrow(C-B)^*\in I$.

In the following Theorem, (i) for $n=1$ and $\B = \{0\}$ gives a first--order definition
of quasimaximality\id{quasimaximal}.  In (ii), we refer to Tarski's classification of the
completions $T$ of the theory of boolean algebras\id{theory!of boolean
  algebras}, in the form
presented in Chang and Keisler \citelab{Chang.Keisler:73}{BA}, 
Section 5.5. They assign invariants $m(\B), n(\B) \in \omega+1$ to Boolean algebras
and prove that two boolean algebras are elementarily equivalent iff they have the
same invariants. Thus if $T$ is a completion of the theory of boolean
algebras, we can also write $m(T),n(T)$ for $m(\B), n(\B)$, where
$\B$ is some model of $T$. 
\begin{theorem}\label{qmax}  The following classes of hh-simple\id{hyperhypersimple} sets are 
definable in $\EE^*$ without parameters. 

\be
\item[(i)] $\{A^*:$ the $n$-th derivative of $\LL^*(A)$ is $\B\}$, where
  $\B$ is a fixed finite boolean algebra or $\B=\{0\}$
\item[(ii)] $\{A^*:\LL^*(A)\models T\}$, where $T$ is any completion of the
theory of boolean algebra's except the one with the invariants
$m(T)=\infty, n(T)=0$.
\ee
\end{theorem}

Note that (ii) is non-trivial since some completions are not finitely
axiomatizable.

\pf (i).  Let $I_0^A$ be the least ideal of $\LL^*(A)$, and for each 
$n$ let
$I_{n+1}^A=\AAA(I_n^A)$.  Then, by Fact \ref{fact2}, there is a formula $\psi_n$ which
uniformly for each $A^*$ defines $I_n^A$ in $\LL^*(A)$. So we can express 
that the quotient boolean algebra of $\LL^*(A)$ through  $I_n^A$ is isomorphic to $\B$.
(ii). is left as an exercise to the reader. \eopnospace

\begin{corollary} \label{qqmax} The following classes of hypersimple
  sets are definable in $\EE^*$.
\be
\item[(i)] $\{A^*: A \ttext{is quasimaximal} \}$
\item[(ii)] $\{A^*: \LL^*(A)$ {is isomorphic to the boolean algebra of finite or cofinite
    subsets of}  $\NN\}$
\ee
\end{corollary}

\pf (i) and (ii) follow from (i) of the preceding theorem with $\B=
\{0\}$ and  $\B=
\{0,1\}$, respectively. \eopnospace

\section{Coding a recursive graph in $\EE$}
\label{CodeRG}
We will  develop a scheme  
\begin{equation}
 \label{graph-scheme}
  \phi_{dom}(x; \ol p), \phi_{\eqq}(x,y; \ol p), \phi_{E}(x,y; \ol p)
\end{equation}
 
for coding a recursive directed
graph\id{graph!recursive} $(\NN,E)$  into $\EE$. Applying this to the graph (\ref{VN}), we
obtain  a scheme as in Example \ref{codepa}. In particular,
$\phi_{num}(x;\ol p)$ expresses that $x$ is a minimal element in the
copy of the
partial order $(V_\NN, E_\NN)$ coded in $\EE$, and the equality
formulas $\phi_\eqq$ are the same.

 Let $A$ be any c.e.\ set such that  $\LL^*(A)$ is a
a boolean algebra with infinitely many atoms. Each atom of $\LL^*(A)$ has the form
$(A\cup H)^*$ for some computable set $H$. Now let
$$\mathcal{H}= \{H: H \ \mbox{computable} \lland (A\cup H)^* \mmbox{atom in}  \LL^*(A) \}.$$

Since the index set of $ \mathcal{H}$ is computable in $\ES^{(4)}$, there is a
$\ES^{(4)}$-sequence $(H_i)$ of computable sets such that $(A\cup
H_i)^*_{i \in \NN}$ lists the atoms of  $\LL^*(A)$ without
repetitions. {\it The variable  $H$ will range over sets in  $ \mathcal{H}$.}

\begin{remark}\label{Lachlan}
{\rm Lachlan \citelab{Lachlan:68*3}{3}  proved that each $\SIII$-boolean algebra $\BB$ is isomorphic to
some $\LL^*(A)$ (also see Soare \cite{Soare:87}). If $\BB$ is the boolean algebra of finite or cofinite subsets of
$\NN$, then  the set $A$ obtained by his construction has the
property that there is in fact a $\ES''$-list  $(H_i)$ as above.} \end{remark}

We will introduce  all-together six   acceptable ideals  of $\BB(A)$. The parameters
needed to define them via the ideal definability lemma will constitute the list of parameters
coding the graph $(\NN,E)$. For a set $\mathcal{C} \sub \BB(A)$, let

$$\ID{\mathcal{C}} \ = \mmbox{\it the ideal of } \BB(A) \mmbox{\it generated by} \RR(A) \cup
\mathcal{C}.$$

Hence  $\ID{\mathcal{C}}$ consists of the sets in  $\BB(A)$ which are
contained  in a set $R \cup X_1 \cup \ldots \cup X_n$, where $R \in
\RR(A)$ and $X_i \in \mathcal{C}$. 
Clearly, if $\mathcal{C}$ is $\Sk \ (k \ge 3)$ then
$\ID{\mathcal{C}}$ is $k$-acceptable.

We  obtain a u.c.e.\ partition
\begin{equation}
\label{Pk}
(A_k)_{k\in\NN} 
\end{equation}

 of $A$ by modifying the proof of the
Friedberg Splitting Theorem\id{splitting!in computability theory} in Soare \cite{Soare:87} so that a
splitting of $A$ into infinitely many sets is produced.
We intend to use $A_k/_\equiv$ to represent the vertex $k$, where
$\equiv$ is  the  equivalence relation defined via $\phi_{\eqq}$.  By
the argument in \cite{Soare:87}, for each c.e.\ $W$ and each
$k$,
$$
W-A \ttext{non-c.e.} \Rightarrow W-A_k
\text{ non-c.e.}
$$ 
In particular, $H-A_k$ is non-c.e.\ for each $H\in \mathcal{H}$,
and hence $H\cap A_k$ is not computable. 

For a better understanding, we first   consider a simplified version of the
coding scheme, ignoring the necessity of a nontrivial equivalence relation $\equiv$,
 at the cost of obtaining coding of the recursive graph\id{graph!recursive}
 only in the structure  in $\EE$ enriched by  an additional
unary predicate  for $\{
A_k:k\in\NN\}$. 
Think of $H_{\la k,h \ra }$ as representing the  pair $A_k, A_h$.
  Given a recursive graph  $(\NN,E)$, we define a copy of $E$ on $\{
A_k:k\in\NN\}$ by using two acceptable ideals $\wt J_0$ and $\wt J_1$.
Let

\begin{eqnarray}
\label{JJ}
 \wt J_0 &=  &\ID{ \{H_{\la k,h \ra}\cap A_j: 
 \neg Ekh \ \lor \  (Ekh \lland k\neq j) \}} \\
 \wt J_1 &=  &\ID{ \{H_{\la k,h \ra}\cap A_j: 
 \neg Ekh \ \lor \  (Ekh \lland h\neq j) \}}. \nonumber 
  \end{eqnarray}

Thus all the intersections $H_{\la k,h \ra} \cap A_j$ go into $\wt J_0$
[$\wt J_1$],
unless $Ekh $ holds  and $k=j$ [$h=j$].
Then one can recover $E$ from $\wt J_0, \wt J_1$
because
\begin{equation}
Eij\Leftrightarrow (\exists H\in \mathcal{H})[H\cap
A_i\not\in \wt J_0\ \land \ H\cap A_j\not\in \wt J_1
].\label{simplified}
\end{equation}
This can be verified using the facts that  for $H,N\in \mathcal{H}$, $H\cap A_k\not\in
\RR (A)$ for each $k$, and  either
$H\cap A_k$, $N\cap A_k$ are disjoint on the complement of  a set in
$\RR(A)$ or they are equal.

With an additional unary predicate for $\{
A_k:k\in\NN\}$, a copy of $E$ on this set can be defined
with parameters by (\ref{simplified}), since $ \mathcal{H}$, $J_0$ and
$J_1$ are definable with parameters.

We now describe how to obtain a definable equivalence relation
$\equiv$ such that
 $ A_k \not \equiv A_h$ for $k \neq h$ and $\{X: \ex k \ A_k \equiv X
 \}$ is definable.  After this we will come back to the coding of
 edges.
 The proof was inspired by Rabin's uniform coding
 of finite graphs in
 {\it boolean pairs}\id{boolean pair!coding in} (i.e.\ boolean algebras with a distinguished
 subalgebra, see 
 Burris and Sankappanavar \citelab{Burris.Sankappanavar:81}{BP}).
However, we don't make an explicit use of boolean pairs.

We picture the sets $A_k$ as {\it columns} and the sets
$H_i \cap A$ as {\it rows}.  The intersection of a column and a row 
is a nontrivial splitting of $A$.  Our goal is to find a parameter definable
collection of splittings of $A$ including the columns, and to define an equivalence relation such that each split
in the collection is equivalent to just one column. We call this
collection of splittings the {\it approximations to columns}.

Observe that, for each noncomputable $Q$, one can uniformly in an
index for $Q$ choose a maximal ideal $J$ of $\BB(Q)$\id{ideals!maximal}
(i.e.\ $|\BB(Q)/_{\ns J}|=2$) which is
4-acceptable,  as follows.  Let $(U_e )_{e\in\NN}$ be a
$\Delta^0_3$-listing of $\BB (Q)$ as described in Notation \ref{Enot}.
One builds an ascending sequence $(X_k)_{k\in\NN}$ of
elements of $\BB (Q)$ which generates  $J$, ensuring that
$(\forall e )[U_e\in J\lor\overline U_e\in J ]$
(to make $J$ maximal) and
$(\forall k) [Q-X_k\text{ noncomputable}]$
(to ensure $\RR(Q) \sub J \not=\BB (Q)$).    Let
$X_0=\emptyset$.  Inductively, for $k>0$, one has to make a
decision, computably in  $\emptyset^{(3)}$, if
$$
X_n=X_{n-1}\cup U_n\text{ or } X_n=X_{n-1}\cup\overline U_n.
$$
If one of these sets has a computable complement $R$ in $Q$,
one has to take the other (i.e.\ $R$ is added to $X_{n-1}$).
If both are non-computable, one can decide either way.  
This procedure guarantees $\RR(Q) \sub J$: if $\overline U_n$ is computable, then the first
set has the computable complement $\overline
X_{n-1}\cap\overline U_n$, so one  opts for $\overline U_n\in J$.  This
shows $\RR (Q)\subseteq J$, so $J$ is 4-acceptable.

Now choose such a maximal 4-acceptable  ideal $J_{i,k}$ for each $Q=
H_i \cap A_k \ (i,k \in \NN)$ in a uniform way, and
let 

\begin{equation}
\label{II}
I = \ID{\bigcup \{J_{i,k}: i,k \in \NN\} }.
\end{equation}

Also let 

\begin{equation}
\label{IA}
I_A =\ID{
\{A_k: k \in \NN \}}.
\end{equation}

 The approximations to columns will be in $I_A$. 
The two  ideals are 5-acceptable, and $I \sub I_A$.
To define a set of approximations to columns we let

\begin{equation}
  \label{att1} \phi_{dom}(X) \equiv X \in  I_A 
            \lland \fa H X\cap H/_{\ns I}
  \mmbox{is atom in} \BB(A)/_{\ns I}].
\end{equation}

This property can be expressed as a first-order property  of
 parameters  coding   the acceptable
ideals via the formulas of the ideal definability lemma. Moreover it is satisfied by
each $A_k$. 
Now, to express that $X,Y$ approximate the same column, one is tempted
to use the formula 

\begin{equation} \label{attt} \fa H \ (X\trup Y) \cap H \in I. \end{equation}

However,
some $X\sq A$ could satisfy (\ref{att1}) without approximating a
column,
 because it may happen that for two different $A_k$'s
there are infinitely many $H_i$ such that $X\cap H_i = A_k \cap
H_i$. Thus, $X$ seems to choose two different columns at the same
time. To avoid this, we restrict the set of $H$ considered in
(\ref{attt}). Let $S \sub \NN$ be any $\Sigma^0_5$-set which is
maximal in the lattice of $\EE^5$ of $\Sigma^0_5$ sets (i.e., $S^*$ is
a co-atom in $(\EE^5)^*$). The existence of $S$ follows by relativizing
Friedberg's maximal set construction (see Soare \cite{Soare:87}) to
$\ES^{(4)}$. Let

\begin{equation} \label{IS} I_S =\ID{\{H_i \cap A: i
\in S\}}. \end{equation}

Then $I_S$ is a 5-acceptable ideal representing the set of atoms of
$\LL^*(A)$ which are ``in $S$''. Finally let 

\begin{equation} \label{IM} I_H = \ID{\{H_i \cap A: i
\in \NN\}}.\end{equation}

In Table \ref{ideals} we summarize the definitions of acceptable
ideals which are needed for the coding of $(\NN, E)$. \vspace{12pt}

\begin{tabular}{llll}
\label{ideals}
Symbol  & Defined in & Function                          & $k$-acc.  \\
              &                   &              & for $k=$  \\              \hline
$I$         & \ref{II}        & $H_i \cap A_k$ atom  in $\BB(A)/_{\ns I}$ & 4 \\

$I_A$     &\ref{IA}        &auxiliary                & 3 \\
$I_S$     &\ref{IS}        & represents $S$, a maximal set in $\EE^5$      & 5 [3] \\
$I_{H}$ & \ref{IM}      &auxiliary                  & 5 [3] \\

$J_0, J_1$ &\ref{JJM} & Code edge relation on $\{A_k\beqq: k \in \NN\}$ & 6 [4]\\
\\ 

  \multicolumn{4}{c}{Table \ref{ideals}. Numbers [k]
    are for $A$ as in Remark \ref{Lachlan}}

\end{tabular}

 \vspace{12pt}

Modify (\ref{attt}) as follows:


\begin{eqnarray}
  \label{equal}
  \phi_{\eqq}(X,Y;\ol P) & \equiv & \ex U \in I_H \
   \\ 
 &&  \fa H  [H \cap A
  \not \in I_S \lland H \cap U \in
\RR(A) \RRA \nonumber \\
 && \ \ \ (X\triangle Y) \cap H \in I]. \nonumber
\end{eqnarray}

(Recall that the variable $H$ ranges over $ \mathcal{H}$.) This formula expresses that except for on finitely many relevant rows,
$X$ and $Y$ behave similarly. It   clearly defines an equivalence relation. We claim the
following (omitting the list of parameters).

\begin{claim}
\label{Pkh}
\be
\item[(i)] for each $k\neq h$, $\phi_{dom}(A_k) $ and $\neg
  \phi_{\eqq}(A_k, A_h)$.
\item[(ii)] $\fa X [ \phi_{dom}(X) \RRA \ex ! k \ \phi_{\eqq}(X, A_k)].$
\ee
\end{claim}

\pf (i). $\phi_{dom}(A_k) $ is obvious. If $k \neq h$, given  $U
\in I_H$ choose an $r \not \in S$ such that $H_r \cap A \not \sub^*
V$. Then $H_r \cap A_k$ and $H_r \cap A_h$ represent different atoms
in $\BB(A)/_{\ns I}$.

(ii). Clearly there can be at most one $k$ such that $\phi_{\eqq}(X,
A_k)]$. For the existence of $k$,
since $X \in I_A$, $X \sub A_0 \cup \ldots \cup A_n \cup R$ for some
$n \in \NN, R\in \RR(A)$.
So, in $\BB(A)/_{\ns I}$ the following relation between atoms holds for
each $H$: 

$$(X\cap H)/_{\ns I} \le \sup_{k \le n} (A_k\cap H)/_{\ns I}.$$

Thus  there is a $k \le n$ such that $X\cap
H_i/_{\ns I} = A_k\cap H_i/_{\ns I}$
for infinitely many $i, H_i \cap A \not \in I_S$. The set $\DDD= \{ i: X\cap
H_i/_{\ns I} = A_k\cap H_i/_{\ns I} \}$ is $\Sigma^0_5$ and $\DDD \not
\sub^* S$. Since $S$ is maximal, $\DDD \cup S =^* \NN$. If $U = \bigcup\{H_i
\cap A: i \not \in \DDD\cup S\}$, then $U \in I_H$ and, for each $H$ such
that $ H \cap U  \in \RR(A)$ and $H\cap A \not \in I_S$, $X\cap H/_{\ns I} = A_k\cap
H/_{\ns I}$.  \eop

Summarizing, the splittings $X$ satisfying $\phi_{dom}$ are those
which, for some $k$ on almost all rows $H \not \in I_S$ behave like
$A_k$. The finitely many exceptions must be taken into consideration
when determining $\phi_E$. We have to represent  an edge from $A_k\beqq$ to $A_h\beqq$
on infinitely many rows $\not \in I_S$. So pick a  sequence
$N_{\la k,h,i \ra }$ without repetitions of elements $H_j $ such that $H_j \cap
A \not \in I_S$. Such a  sequence can be chosen computably in
$\ES^{(5)}$ because $S$ is $\Sigma^0_5$. Modifying (\ref{JJ}), let

\begin{eqnarray}
 \label{JJM}
  J_0 &=  &\ID{ I \cup \{N_{\la k,h,i  \ra}\cap A_j: 
 \neg Ekh \ \lor \  (Ekh \lland k\neq j) \}} \\
  J_1 &=  &\ID{ I \cup \{N_{\la k,h,i \ra}\cap A_j: 
  \neg Ekh \ \lor \  (Ekh \lland h\neq j) \}}. \nonumber 
    \end{eqnarray}

Then $J_0, J_1$ are 6-acceptable. Let 

\begin{eqnarray*}
\phi_{E}(X,Y;\ol P)\ = \ \fa U \in I_H \ \ex H [H \cap A\not\in I_S
\lland H \cap U \cap
A \in \RR(A) \lland\\
H\cap
X\not\in  J_0\ \land \ H\cap Y \not\in  J_1 ].
\end{eqnarray*}

\begin{claim} 
\label{Ekh}
If $\phi_{\eqq}(X,A_k; \ol P) \lland \phi_{\eqq}(Y, A_h; \ol P)$,
  then $\phi_{E}(X,Y;\ol P) \LLR Ekh.$ \end{claim}

\pf  Pick $V \in I_H$ such that $\fa H [ H\cap A\in I_S \lland H\cap V
\in \RR(A) \RRA (X\triangle A_k) \cap H \in I$ and $(Y\triangle
A_h) \cap H \in I$. 

``$\LA$''.  Given $U \in I_H$, pick any $H, H \cap A  \nii I_S$ such that
$H\cap (U\cup V) \in \RR(A)$. Because $(X\triangle A_k) \cap H \in I$ and $(Y\triangle
A_h) \cap H \in I$,  $H\cap X \nii J_0$ and $H\cap Y \nii J_1$.

``$\RA$''.  If not $Ekh$, we show $\neg \phi_E(X,Y, \ol P)$: for each
$H,  H \cap A  \nii I_S$ such that $H\cap V \in \RR(A)$, $H\cap A_k \in J_0$ implies $H\cap
X \in J_0$, and  $H\cap A_h \in J_1$ implies $H\cap
Y \in J_1$. This conclude the construction of the scheme
(\ref{graph-scheme}) and the proof that $(\NN,E)$ can be coded via
this scheme.  \eop

\section{Interpreting true arithmetic in $\Th(\EE)$}

\label{ITAE}
\begin{theorem}[\cite{Harrington.Nies:nd}]

$\ThN$ can be interpreted in $\Th(\EE)$.

\end{theorem}

This theorem was first proved by L.\  Harrington. A simpler proof
appeared in Harrington and Nies \citelab{Harrington.Nies:nd}{2}. The
present proof   is a simplification once again because of an  improved coding
scheme $S_M$.

\pf  We will apply Fact \ref{DITA}. First we must   develop  a scheme
$S_M$  as in Example \ref{codepa} to satisfy (\ref{code}). Combining the coding of a copy of
$\Nops$ in a recursive graph $(V_\NN, E_\NN)$ in (\ref{VN}) with the
coding of any recursive graph in $\EE$ from the preceding Section, we
obtain formulas  (\ref{SMformulas}) with parameters $\ol P$. The list $\ol P$ consists
of a parameter $A$ (the base set) and parameters to define ideals as in Table
\ref{ideals}. As before let the variable $H$ range over computable
sets such  that $A\cup H^*$ is an atom in $\LL^*(A)$. 
Beyond the correctness condition $\alpha_0(\ol P)$ from Example \ref{codepa} we add some
more \id{correctness condition}s which enable us 
 to   quantify
over $\Sk$ subsets of $M_{\ol P}$  (in the sense of Definition \ref{repr}),
for a sufficiently large $k$, in order to satisfy (\ref{see}).

First, using Corollary \ref{qqmax} (ii)  we require in a first-order
way that $\LL^*(A)$ is {isomorphic to the boolean algebra of finite or cofinite
    subsets of} $\NN$ and that  $I^*_H$  is the ideal generated by
  the  set
  of  atoms in $\LL^*(A)$. As before let $(H_i)$ be a $\ES^{(4)}$-sequence  of computable sets such that $(A\cup
H_i)^*_{i \in \NN}$ lists the atoms of  $\LL^*(A)$ without
repetitions. Then the collection 
$\Sigma^0_5$ ideals of $\LL^*(A)$ contained in the  ideal $\wt I_H$ generated by
the atoms under inclusion is canonically isomorphic to
$\EE^5$. Moreover, as described in Section \ref{Defhhsim} we can
quantify over this collection of ideals. {\it In the
following we identify those ideals with $\Sigma^0_5$-sets.}  Let
$\wt I_S = \{(A\cup H)^*: A\cap H \in I_S\}$. As a  correctness
condition 
we require that 

\begin{equation} \label{Sismax} \wt I_S \ttext{(viewed as a}
  \Sigma^0_5-\text{set)  is maximal in} \ \EE^5. \end{equation}

This completes the description of the scheme $S_M$.

For $X,Y\in \BB(A)$, let

$$\DDD_{X,Y} = \{i: (H_i\cap X)/_{\ns I}= (H_i \cap Y)/_{\ns I}\}.$$

Since the sequence $(H_i)$ is $\ES^{(4)}$ and $I$ is 5-acceptable, 
$\DDD_{X,Y}$ is a $\Sigma^0_5$-set. So, under the identification we
make, 

\begin{equation}
  \label{alter}
    \DDD_{X,Y} \sub^* \wt I_S \ \vee \ \DDD_{X,Y} \cup \wt I_S =^* \NN.
\end{equation}

Moreover, $\DDD_{X,Y} \cup \wt I_S =^* \NN \LLR \phi_{\equiv}(X,Y;\ol
P)$.

 The following  lemma
will allow  us to quantify over $\Sk$-subsets of $M_{\ol P}$: if $\CCC
\sub M_{\ol P}$ and $\hat \CCC \sub \BB(A)$ represents $\CCC$ (in the
sense of Definition \ref{repr}),
then $\CCC$ can be recovered from $\ID{\hat \CCC}$, the ideal generated by
$\hat \CCC \cup \RR(A)$.

\begin{lemma} \label{Independence} Suppose that $\hat \CCC \sub \BB(A)$ represents $\CCC \sub
  M_{\ol P}$. Then for each $X \in \BB(A)$,

$$ \phi_{num}(X; \ol P) \lland X \in \ID{\hat \CCC}) \RRA X \in \hat \CCC.$$
\end{lemma}
\pf  Since $\hat \CCC$ represents $\CCC \sub M_{\ol P}$,
$\phi_{num}(X; \ol P)$ holds for each $X \in \hat \CCC$. We now use an
argument similar as the one to prove Claim  \ref{Pkh} (ii). Suppose $X
\in \ID{\hat \CCC}$, then $X \sub Y_0 \cup \ldots \cup Y_n \cup R$ for some
$n, R\in \RR(A)$. Then  $\fa i \ \ex k \le n (X\cap H_i)/_{\ns I} =
(Y_k\cap H_i)/_{\ns I}$.
If $\DDD_{X,Y_j} \sub ^* \wt I_S$ for each $j \le n$, then $\NN =
\bigcup_{j \le n} \DDD_{X,Y_j} \sub^*  \wt I_S$, contrary to the
correctness condition
(\ref{Sismax}). So, by (\ref{alter}), there is $j \le n$ such that
$\DDD_{X,Y_k} \cup \wt I_S =^* \NN$, hence $ \phi_{\equiv}(X,Y_j;\ol
P)$. Since $\hat \CCC$ is 
 is closed under 
$\equiv$, this implies that $X \in \hat \CCC$.  \eop

By the ideal definability lemma \ref{ideal definability lemmaE} and since $\ID{\CCC}$ is
$k$-acceptable
for any $\Sk$-set $\CCC \sub \BB(A)$ if $k \ge 3$, the set of
$\Sk$-subsets
 of $M_{\ol P}$ is weakly uniformly definable (in fact, uniformly
 definable  if $k$ is sufficiently large). By   Fact \ref{DITA}, this completes the
proof. \eop

\section{Fragments of $\Th(\EE^*)$}

\label{Efragments}

We will investigate decidability and undecidability for fragments of
$\Th(\EE^*)$ as a lattice. Lachlan \citelab{Lachlan:68*2}{1} proved
that $\Pi_2$-$\Th(\EE^*)$ is decidable. Here we will use the coding
methods developed in Section  \ref{CodeRG} in order to prove that 
$\Pi_6$-$\Th(\EE^*)$ is undecidable. While seemingly far from optimal,
this result improves the bound one obtains from the coding in
Harrington, Nies \citelab{Harrington.Nies:nd}{3} by two quantifier alternations. 
For  an optimal result, one would wish  to develop  a coding of a
sufficiently complex class, like the class of finite undirected graphs
(see Theorem \ref{Lavrov}), using
a $\Sigma_1$-scheme  with parameters. By the methods of 
Section \ref{fragments}, this would imply  the undecidability of the
$\Pi_3$-theory of $\EE^*$. However, such a proof is not possible,  since it
would show that the class of finite distributive lattices with the
reduction property (\ref{RedPr}) (also called {\it separated}
distributive lattices)\id{distributive!lattice!separated}
has an undecidable theory, contrary to a result of
Gurevich \citelab{Gurevich:83}{1}. The argument is as follows: Suppose that, via some
$\Sigma_1$-scheme in the sense of Section \ref{fragments} (or even a scheme (\ref{general scheme}) which
consists solely of
$\Sigma_1$ formulas) we can code each finite
undirected  graph\id{graph!undirected}  $(V,E)$ (say), using appropriate parameters $\ol p$.  Let $G$
be a   finite distributive sublattice of $\EE^*$ which contains $\ol p$, all the
elements of $\EE^*$ representing the vertices in $V$ and also witnesses
for all $\Sigma_1$--formulas involved to code $(V,E)$. Then $G$, and in
fact any distributive lattice $D$ such that $G \sub D \sub \EE^*$ codes
$(V,E)$ via the same scheme and parameters. Now by Lachlan
\citelab{Lachlan:68*2}{HH} let $D$ be such a
lattice which is also finite and satisfies the reduction principle.  In this
way, we have obtained a uniform coding of a complex class in  the class of finite distributive lattices with the
reduction property.

We  conclude that the best we can hope for  by the
standard coding methods is the undecidability of the $\Pi_4$--theory
of $\EE^*$,
which still would require a far more direct coding than the one
presented here.

\begin{theorem}
\label{Pi6}
The $\Pi_6$-theory of $\EE^*$ as a lattice is undecidable.
\end{theorem}

\pf We will show that the class of finite directed graphs $(V,E)$ can be  uniformly
$\Sigma_4$-coded in $\EE^*$. Then, by Theorem \ref{Lavrov} and (ii) of the Transfer Lemma
\ref{transfer}, $\Pi_6-\Th(\EE^*)$ is undecidable. The coding is based
on the formulas in  Section \ref{CodeRG}, but all the formulas will now be
interpreted in $\EE^*$.
We use same-type lower case letters to indicate this difference. For
instance, the formula (\ref{phi1}) becomes 

$$ x \sq a \lland x\wedge c \le d.$$

We use the abbreviations ``$ u = v \times w$'' for ``$u \wedge v = 0
\lland u \vee v = w$'' (so   $u \sq w \LLR \ex v \ u \times v=w$).

 The advantage of working in $\EE^*$ is that 
the formula (\ref{phi1}) is essentially quantifier free, since the
$\Sigma_1$-condition
``$x \sq a$'' can be stated independently of the rest. Moreover,
 (\ref{phin}) is
$\Sigma_{2n-2}$. For example,

$$\begin{array}{rl}

\varphi_{2}(x;d,c_0, c_1,a) &\LR  x\sqsubset
a\lland\\
& \exists r\le a \ex \wt r[ r \times \wt r =1 \lland \fa s \le a \ \fa \wt s( s
\times \wt s =1 \RRA \\
&  \ \ \ x \wedge s \wedge c_1  \wedge d \le c_0)].
\end{array}$$

 Thus, a $2n+1$-acceptable ideal of $\BB(A)$ is
defined in $\EE^*$ by a $\Sigma_{2n-2}$-formula with parameters.
 We will also work with the particular
hh-simple set\id{hyperhypersimple} obtained from Lachlan's
construction (see Remark \ref{Lachlan}) and  fix $D \sub_{sm}
A$. This is possible since  here we are
satisfied with a {\it particular} list of parameters in $\EE^*$ and  containing $A^*, D^*$  which
codes a  finite directed graph\id{graph!directed}.

Suppose we are given  a finite directed graph  $(V,E)$,  where $V= \{0,\dots,v\} \sub \NN$. We
follow the definitions in Section \ref{CodeRG}, but with a finite partition
$(A_k)_{k \le v}$ instead of the infinite one used in (\ref{Pk}). 
Since we work with the  particular set $A$, we can assume that $S$ is
a maximal set in $\EE^3$, and the sequence $(N_{\la k,h,i \ra})$
introduced before (\ref{JJM}) is computable in $\ES^{(3)}$. Then 
the acceptable ideals in Table \ref{ideals} have the (lower) complexities
indicated in cornered brackets. The formula $\phi_{dom}(X; \ol
P)$, rewritten for $\EE^*$, becomes:

\vsp

$$\begin{array}{rl}
 \phi_{dom}(x;\ol p) & \LR \   x \in  I_A^* \lland \\
                                 &\fa h \\
                                      &   \ \ \   (x\wedge h)/_{\ns I^*} \mmbox{is atom in}
                                 \BB^*(a)/_{\ns I^*}.
\end{array}$$

\vsp

Now $I_A$ is 3-acceptable, so $I_A^*$ is defined in $\EE^*$ by a $\Sigma_1$-
formula. 
Moreover,  $I$ is 4-acceptable, so $I^*$ is defined in $\EE^*$ by a $\Sigma_2$-
formula. Then  ``$(x\wedge h)/_{\ns I^*} \mmbox{is an atom in}
                                 \BB^*(a)/_{\ns I^*}$'' can be
                                 expressed by a $\Pi_3$-formula, and
                                 $\phi_{dom}(x;\ol p)$ is $\Sigma_4$.

Next we look at $\phi_{E}(X,Y;\ol P)$ in $\EE^*$ and obtain
\vsp

$$\begin{array}{lr}
 \phi_{E}(x,y;\ol p) &\LR \  \fa u \in I^*_M \ \ex h[h \wedge a
 \not\in I^*_S \lland  h
 \wedge u \wedge 
a \le d   \lland\\
&h\wedge
x\not\in  J^*_0\ \lland \ h \wedge y \not\in  J^*_1 ].

\end{array}$$

\vsp

This describes a  $\Pi_4$-formula, since $J_0, J_1$ are 4-acceptable. In order to obtain a
$\Sigma_4$-scheme, we use the $\Sigma_4$-formula $\neg  \phi_{E}(x,y;\ol
p)$ to code the complement $\ol E=\{0,\dots,v\}^2 -E$ of the edge
relation $E$ on $\{x:\phi_{dom}(x; \ol p)\} $. Moreover, from
(\ref{JJM}) for $\ol E$ we obtain two further 4-acceptable ideals $\ol
J_0, \ol J_1$ and a $\Pi_4$-formula $\phi_{\ol E}(x,y;\ol
p)$, which is like  $\phi_E$ but uses the parameters for 
$\ol J_0^*, \ol J_1^*$. Now use $\neg \phi_{\ol E}(x,y;\ol
p)$ to code $E$. This completes the description of the
$\Sigma_4$-scheme and thereby the proof. \eop

\section{The theories of relativized versions of $\EE$}

In this section we investigate and compare the theories of the lattices
$\EE^Z$ of sets c.e.\ in an oracle set $Z$, in particular for sets $Z
\sub \NN $ with the following property: $Z$ is called
 \iid{implicitly definable} in arithmetic if there is a formula $\psi_Z$ in
the language $L(+,\times)$ extended by a unary predicate $R$ such
that, for each $X\subseteq\NN$,
$$
(\NN ,+,\times)\models\psi_Z(X)\Leftrightarrow X=Z. \label{tag9}
$$
Note that a set which is  implicitly definable in arithmetic is 
$\Delta^1_1$, hence \ird{hyperarithmetical},   and
that implicit definability of $Z$ only depends on the arithmetical degree 
of $Z$. Hence each $Z$ which is in the same arithmetical
 degree\id{degree!arithmetical}  as
some $\emptyset^{(\alpha )}$,
 $\alpha$ a recursive ordinal, is  implicitly definable in arithmetic.
  However, "most" hyperarithmetical
sets  are not  implicitly definable in arithmetic, since both \ird{arithmetically 
generic} sets and \ird{arithmetically random} sets $Z$ cannot be  implicitly 
definable (this is described in more detail in  Nies\cite{Nies:96}).

We prove that, if 
$Z$ is implicitly definable in arithmetic, then $\hbox{Th}(\NN 
,+,\times,Z)$ can be interpreted in
$\hbox{Th}(\EE^Z)$\id{relativization!of $E$}.  Since an
interpretation in the other direction exists as well, the
two theories have the same $m$-degree.
To do so,  we exploit the coding power of a specific
collection of formulas in $\EE^Z$ to show that
 for some fixed
$c\in\NN$, if $Z$ is  implicitly definable in arithmetic and
$Z^{(c)}\neq W^{(c)}$, then $\EE^Z$ is not elementarily equivalent
to $\EE^W$. (In Shore \citelab{Shore:81}{rel}, similar ideas were
first applied to relativizations of 
the structure of $\Delta^0_2$ Turing-degrees.)
  In particular, if $Z=\emptyset^{(\alpha)}$,
$W=\emptyset^{(\beta)}$, where $\beta<\alpha$ are recursive ordinals, then
$\EE^Z\not\equiv\EE^W$. For finite $\alpha ,\beta$, this
 negatively  answers  the question
mentioned in the introduction to this chapter whether $(\Sigma^0_p,
\sub)$ is elementarily equivalent to $(\Sigma^0_q, \sub)$ for $p < q$.
As a further application, if $Z$ is sufficiently complex, namely $Z
\not \in Low_c$ ($c$ as above), then 
  $\EE^Z$ is not elementarily 
equivalent to $\EE$. This includes the case that $Z$ is arithmetically 
generic. We note that, for all arithmetically generic $Z$, the 
relativization $\EE^Z$ has the same theory. Similar remarks apply to 
arithmetically random sets (Nies \citelab{Nies:96}{1}).
  
We make some observations which will enable us to interpret
true arithmetic in $\hbox{Th}(\EE^Z)$ for
 each $Z$ and
$\hbox{Th}(\NN ,+,\times,Z)$ in $
\hbox{Th}(\EE^Z)$ if $Z$ is  implicitly definable in arithmetic.
In $\EE^Z$, define $\BB(A)$ as in \ref{Enot} and let $\RR(A)\sub \BB(A)$ be the
collection of subsets of $A$ which are computable in $Z$. An ideal $I$
of $\BB(A)$ is $k$-acceptable (relative to $Z$)\id{acceptable!relative
  to an oracle}  if $\RR(A) \sub I$ and $I$ has a $\Sigma_k^Z$ index
set. The proof of the ideal definability lemma relativizes to $Z$, so  in $\EE^Z$ the
class of $k$-acceptable ideals is uniformly definable for all odd  $k \ge 3$. Thus the
scheme $S_M$ also works in $\EE^Z$, and, relativizing the
considerations in Section \ref{ITAE}, we obtain:

\begin{corollary} \label{relTAZ} $\ThN$ can be   interpreted 
 in $\hbox{Th}(\EE^Z)$ for
 each $Z$. \eopnospace
\end{corollary}
 
Moreover, we observe
\begin{fact}
\label{observE}
\be
\item[(i)]  If   $\varphi (\overline X;\overline P)$  is a
$\Sigma^0_k$ formula with parameters in the language of $\EE$,
then for each $Z$, the index set with respect to the indexing
of $\EE^Z$, $(W_e^Z)_{e \in \NN}$, of the relation defined by $\varphi$ with a fixed
parameter list is computable in $Z^{(k+2)}$.

\item[(ii)] For some fixed number $h$ (which does not 
depend on Z), for each $M$, there
is $g\le_T Z^{(h)}$ such that
$$
(\forall i)[W^Z_{g(i)}/_\eqq=i^M].
$$
\ee
\end{fact}

\pf (i) is immediate since ``$W^Z_i \subseteq W^Z_j $'' 
is computable in $Z^{(2)}$. For (ii),
suppose that $M=M(\overline 
P)$.  Let $\varphi_S(X,Y;\overline P)$ be a formula defining the
successor function in (any) $M(\overline P)$.  By (i),
the corresponding binary relation on indices is computable
in $Z^{(h)}$ for some fixed number $h$, so there is a
partial "choice" map $f$ which can be computed with the oracle
$Z^{(h)}$ such that, in $\EE^Z$,
$$
\varphi_S(W_i^Z, W_j^Z;\overline P)\text{ for some }
j \RRA
\varphi_S(W_i^Z, W^Z_{f(i)}; \overline P).
$$
Fix $i_0$ such that $W^Z_{i_0}/_{\ns I}=0^M$.  Then, by iterating
$f$ with $i_0$ as an initial value, obtain $g$ as desired. \eop

From  (ii)  one immediately obtains  the relativization of 
Fact \ref{StSk}:
 for each  structure $M$ coded in $\EE^Z$ via $S_M$, $\{ e:W^Z_e/I\text{ is a
standard number of }M\}$ is $\Sigma^0_p(Z)$ for some fixed $p$.

\begin{theorem} \label{TAZ}  If $Z$ is  implicitly definable in arithmetic,
then there are interpretations of theories  which show
$\hbox{Th}(\EE^Z)\equiv_m\hbox{Th}(\NN ,+,\times,Z)$. 
\end{theorem}

\pf Suppose that $Z$ is  implicitly definable in arithmetic.  To interpret
$\hbox{Th}(\NN,+,\times,Z)$ in $\hbox{Th}(\EE^Z)$ we need an
extended scheme which enables us to encode structures
$(M,Z_M)$, where $M= M_{\ol P}$ is a  coded copy of $\NN$
 and $Z_M$ is $Z$, viewed
as a subset of $M$.  Let $\psi_Z$ be a formula describing
$Z$ as in (\ref{tag9}).  Given $M$ as above, suppose that $\hat{Z}_M
\sub \BB(A)$ represents $Z$ and let $I_Z = \ID{\hat Z_M}$ (the ideal
generated by $\hat Z_M$ and $\RR(A)$).
  Then, using the map $g$ from Fact \ref{observE} (ii),

    $$I_{ Z} = \ID{\{W^Z_{g(n)}:n\in Z\}^\eqq}.$$

  Since $g\le_T Z^{(h)}$
for some $h$, $I_{Z}$ is $q$-acceptable (in $\EE^Z$)
for some $q$. Suppose that $M$ is standard.  Since   Lemma \ref{Independence}
also holds in $\EE^Z$, $\phi_{num}(P;\ol P)$ implies that 

$$P\beqq \in Z_M \ \Leftrightarrow P\in I_{\hat Z_M}.$$

In the extended scheme, expand the list of parameters by
parameters defining a $q$-acceptable ideal $J$ of $\BB
(A)$.  As an additional  \id{correctness condition} for  the scheme we
require that $$[\phi_{num}(X;\ol P) \lland \phi_{num}(Y;\ol P) \lland X \in J \lland X\eqq Y]  \ \RA \  Y \in  J.$$ 

 Let $W$ be the
subset of $M$ represented by $J \cap \{X: \phi_{num}(X; \ol P)\}$
(the intended meaning is that $W=\hat Z_M$).

 The interpretation of $\hbox{Th}(\NN
,+,\times,Z)$ in $\Th(\EE^Z)$ is now given by 

$(\NN
,+,\times,Z)\models\varphi\ \Leftrightarrow$ 

\begin{quote} for some $(M,W)$
 coded via  the extended scheme,  $M$ is
standard,$W$ (as a subset of $M$) 
is represented by $J \cap \{X: \phi_{num}(X; \ol P) \}$, $M_{\ol P}
\models \psi_Z(W)$
 and 
 $(M,W)\models\varphi$. \end{quote}
The right-hand side can be expressed by a first -order sentence
effectively obtained from  $\phi$. \eop

Let $T \sub \NN$. Given an  $M$ 
  coded in $\EE^T$, let $g$ be the function from Fact \ref{observE}
  (ii).  For $Q \sub \NN$, let $J_Q = \ID{\{W^T_{g(k)}: k \in Q\}^\eqq}$.
The following is a key technical fact. 

\begin{lemma}
  \label{psuffbig} For a sufficiently big number $p$ and  any  $M$
  coded in $\EE^T$, the following holds: if $Q \sub \NN$ then 
  $$
Q\text{ is }\Sigma^0_p(T)\LLR
 J_Q \mmbox{is} p-\mbox{acceptable}. 
\label{tag14} $$
\end{lemma}

\pf Let $p\in \NN$ be such that all ideals needed for the coding of $(V_\NN,
E_\NN)$ in $\EE^T$ are $p-1$-acceptable,  the function $g$ in Fact
 \ref{observE} is computable in $\jump{p-1}$ and $X \equiv Y$ is
 recursive in $T^{(p-1)}$ as a relation between indices.

For the direction ``$\RA$'', note that 

\begin{eqnarray*}
   W^T_e\in J_Q &  \LR & \ex r \ \ex k_0, \ldots, k_r \ \ex X_0, \ldots,
  X_r \\         &  &    W_e^T \sub \bigcup_{i = 0, \ldots, r} X_i 
                         \lland \fa i \le r (k_i \in Q
                         \lland W^T_{g(k_i)} \eqq X_i). 
\end{eqnarray*}

It is easy to check that this can be expressed as a 
$\Sigma^0_p(T)$ property of $e$. 
For the direction ``$\LA$'',  if $J_Q$ has
a $\Sigma^0_p(X)$ index set, then, because

\begin{eqnarray*}
  n\in Q&\Leftrightarrow &W^T_{g(n)}\in J_Q \\
&\Leftrightarrow &(\exists e)[g(n)=e\lland W^T_e\in J_Q]
\end{eqnarray*}
and $g \le_T \jump{p-1}$, $Q$ is $\Sigma^0_p(T)$ (this uses Lemma \ref{Independence}). \eop

Now assume  in addition  that $p$ is odd, and let $c=p-1$.
 \begin{theorem}      
If
$Z^{(c)}\not\equiv_T W^{(c)}$ and  $Z$  or $W$ is  implicitly definable in 
arithmetic,
 then $\EE^Z\not\equiv\EE^W$.
\end{theorem}

\begin{corollary}  If $\alpha$ is a recursive ordinal and
$\beta <\alpha$, then
$\EE^{\emptyset^{(\alpha)}}\not\equiv\EE^{\emptyset^{(\beta
)}}$. \eopnospace
\end{corollary} 

{\it Proof of the theorem.}  Assume that $Z^{(c)}\nleq_T W^{(c)}$. 
Then, if $Z^{(p)}\in  \Sigma^0_p(Z)-\Sigma^0_p(W)$.
 Let $\varphi (Y;D,\overline C,A)$ be
the formula obtained from the  ideal 
definability lemma to define  uniformly in
$\EE^T$ for a set $A$ which is c.e., but not computable in $T$ all $p$-acceptable
ideals of $\BB (A)$.
 
 First suppose that
  $Z$ is  implicitly definable in 
arithmetic,  via the description $\psi_Z$. Then the
following is true in $\EE^Z$, but not in $\EE^W$.
\begin{quote}
 There is a structure $(M,Y)$ coded by the extended
scheme  such that $M$ is standard, $(M,Y)\models \psi_Z(Y)$
and, for some list $D,\overline C$, the ``intersection'' of M and the ideal
coded by $D,\overline C,A$ equals $Y^{(p)}$, i.e.

\begin{equation}
\label{PMY}
\forall P[\phi_{num}(P;\ol P) \RRA ((M,Y)\models P/_\eqq\in
Y^{(p)}\Leftrightarrow\varphi(P;D,\overline C,A))].
\end{equation}
\end{quote}

The statement holds in $\EE^Z$ via any standard $M$ and
$Y=Z_M$ (i.e.\ $Z$ viewed as a subset of $M$), for in
this case $J_{Z^{(p)}}$ is $p$-acceptable by Lemma \ref{psuffbig}.  In $\EE^W$, either  does
$\psi_Z(Y)$ hold in no
structure $(M,Y)$, $M$ standard,  defined by the extended scheme, or, if $(M,Y)$ is such a structure, then
(\ref{PMY}) fails.  For, in $\EE^W$, $\{ P\in\BB (A):\varphi
(P,D,\overline C,A)\}$ is an ideal with $\Sigma^0_p(W)$ index
set by the easy direction  of the  ideal 
definability lemma relativized to $W$.  So, if (\ref{PMY}) holds, by \ref{tag14},
$Z^{(p)}\in\Sigma^0_p(W)$, a contradiction to $Z^{(c)} \not \le_T W^{(c)}$.

Now suppose that $W$  is implicitly definable via $\psi_W$. The case that $W^{(c)}\nleq_T Z^{(c)}$
has  already been covered above. Otherwise there is an index $e$  such that 
$\{e\}^{Z^{(c)}}=W$. 
 Then the
following is true in $\EE^Z$, but not in $\EE^W$.

\begin{quote}
 There is a coded standard model M and a list $D,\overline C$
coding a $p$-acceptable  ideal of $K$ of $\BB(A)$ such that 

\bc $U_0 = K \cap \{X: \phi_{num}(X;\ol P)\}$ \ec

is  closed under $\eqq$ and  if $U= U_0\beqq$,
 then for some index $e \in M$, 

\bc $M \models
\psi_W(\{e\}^U) \lland  U \nii \Sigma^0_p(\{e\}^U)$. \ec

\end{quote}
This statement holds in $\EE^Z$ via the ideal $K=J_{Z^{(p)}}$, but fails in 
$\EE^W$, once again by the easy direction  of the ideal definability lemma. 
 \eop

\section{Non-coding and Non-definability Theorems}
 
\label{Nonco}
In the 
last section of this chapter we investigate the limits of definability and
coding in $\EE$. We show that no infinite linear order can be coded
(without parameters) even in the most general way, namely on
equivalence classes of $n$-tuples.  The proof makes use of the fact
that for each partition of $\NN$ into three infinite computable
sets $R,S,T$ there is a canonical isomorphism $\EE\to\EE^3$
given by $X\to (X\cap R,X\cap S,X\cap T)$, combined with a
model theoretic result due to Feferman and Vaught \citelab{Feferman.Vaught:59}{1}
that a first-order property of a tuple in a model of the
form $\A^n$ can be expressed as a certain boolean
combination of first-order properties of the components.
First we prove a noncoding theorem in the context
of uniform first--order definability  with
parameters, which can be considered as a weak version  of the
model--theoretic notion of stability for
$\EE$:  there is no uniform way to define, even  with parameters,
a linear order on arbitrarily large classes $\{ R_1,\ldots
,R_k\}$ of pairwise disjoint  computable sets.  This
implies that infinite no linear order can be coded in a first-order way on  atoms of
$\LL^*(A)$, if $\LL^*(A)$ is a boolean algebra with
infinitely many atoms. 

Hodges and Nies \citelab{Hodges.Nies:97}{1} have  shown that 
 in fact  no infinite linear order can be coded without parameters in  any 
structure isomorphic to a structure  $\A\times \A$ (as $\EE$ is one). However, the proof given here for
$\EE$ contains 
interesting insights into further self-similarity properties of $\EE$ 
and also puts an effective upper bound on the cardinality of a linear
order
which can be coded by a given formula.

If $\A$ can be coded in $\Nops$, then each relation on $\A$ which is
definable without parameters must  be invariant under automorphisms and has
an
arithmetical index relation. The questions arises if a \iid{maximum definability
property}  holds, namely if these two properties actually characterize the
definable relations.  The question has been answered affirmatively for the
structure of $\Delta^0_2$ $T$-degrees by Slaman and Woodin \citelab{Slaman.Woodin:nd}{1}. 
In Harrington and Nies \cite{Harrington.Nies:nd} it was  proved that
the maximum definability property fails for $\EE^*$ (and 
hence for $\EE$) by giving a binary relation as a counterexample.  The
counterexample provided  here is in fact  a subclass of the class of
quasimaximal\id{quasimaximal} sets.

\noindent  Let the variable $\bfR$ range over finite classes of pairwise disjoint
infinite computable sets.
We use the variable $\tilde X$ for  tuples of c.e.\  sets $(X_0,\dots
,X_{n-1})$.

\begin{theorem}  \label{Rfat} For each formula $\varphi (X,Y;\tilde P)$ 
 one can find in an effective way  a number $k$ such that for each $\bfR$,
   $|\bfR | \geq k$, and for each list of parameters $\tilde A$, the relation
$$
\{ (X,Y) : X,Y \in \bfR  \lland  \EE \models \varphi (X,Y;\tilde A)\}
$$
is not a linear ordering of $\bfR$.

\end{theorem}

\begin {corollary}   If $\LL^*(A)$ is a boolean algebra with infinitely many
atoms, then it is not possible to code, even with parameters, an
infinite  linear
ordering on  atoms of $\LL^*(A)$.  

\end{corollary}

\demo{Proof of the Corollary} If $F$ is a set of atoms and $|F|
= k$, then for some $\bfR$ such that  $ | \bfR | = k$, $F = \{ A \cup R^* : R\in
\bfR \}$.
Hence, if $\psi (X,Y;\tilde P)$  defines a linear order on the atoms, then
$\varphi (X,Y; \tilde P,A) \equiv \psi (X \cup A, Y \cup A; \tilde P)$
defines a linear order on sets $\bfR$ of arbitrarily large cardinality.
\eop

{\it Proof of the Theorem.}  Note that, if $R,S$ and $T = \overline {R \cup S}$ are 
infinite, then $\EE \cong \EE^3$ via the map
$$X \mapsto (X \cap R, X \cap 
S, X \cap
T).$$

\noindent By a result of Feferman and Vaught
\citelab{Feferman.Vaught:59}{NCE}, if $\A$ is a structure $k \ge 0$ and
$\varphi (X^0, \dots ,X^{n-1})$ is a formula in the language of $\A$, then

\begin{equation}\label{Fefi}
\A^{k+1}  \models \varphi
\left( \begin{array}{ccc}
a^0_0 & & a^{n-1}_0 \\
a^0_1 & ,\dots , & a^{n-1}_1\\
            & \dots  &            \\
a^0_k & & a^{n-1}_k
\end{array} \right)
\ \Leftrightarrow  \
 \bigvee_{\alpha = 1,\dots ,r} \bigwedge_{i=0,\dots,k}
\A \models \varphi^\alpha_i \left ( a^0_i ,\dots , a^{n-1}_i \right ) 
\end{equation}

for some formulas $\varphi^\alpha_i$ which only depend on $\varphi$ and can  be found effectively.
Thus, whether $\phi(\ol a^0, \ldots , \ol a^{n-1})$ holds in
$\A^{k+1}$ only depends on 
finitely many effectively determined first-order properties of the
components $a^0_i, \ldots, a^{n-1}_i \in \A$ ($i \le k$). This can
be proved by induction on $|\varphi|$.

Now suppose that $\varphi (X,Y; \tilde P)$ defines a linear order $<_L$ on
a set $\bfR$.  By the isomorphisms $\EE \Leftrightarrow \EE^3$ above, an 
element $A \in \EE$ corresponds to the vector

$$
\left( \begin{array}{c}
A \cap R \\
A \cap S \\
A \cap T
\end{array} \right) .
$$
Hence, if $R,S \in \bfR, R \not = S$, then

\begin{eqnarray*}
  R <_L S & \Leftrightarrow  \ \bigvee_{\alpha = 1,\dots,r} &
       (\EE (R) \models \varphi^\alpha_0 (R,\emptyset, \tilde P \cap
       R)  \\
 & & \lland      \EE (S) \models \varphi^\alpha_1 (\emptyset ,S,\tilde P \cap S) \\
   & & \lland  \EE (T) \models \varphi^\alpha_2 (\emptyset, \emptyset, \tilde P \cap T)) ,
\end{eqnarray*}

where $T = \overline {R \cup S}$ and
$\tilde P \cap X = (P_0 \cap X ,\dots ,P_{k-1} \cap X) $.  Note that ``$\EE
(T) \models \varphi^\alpha_2 (\emptyset ,\emptyset ,\tilde P \cap T)$'' does not
depend on the order of $R,S$.  We say that $R <_L S$ {\it via} $\alpha$ if the
disjunct corresponding to $\alpha$
 holds. Now, we can compute a number $M$ such that, for $| \bfR | \ge M$, there exist $\alpha$ and $A,B,C,D \in \bfR$ such that $A <_L B
<_L C <_L D$ and the ordering relations hold all via $\alpha$.  This is
verified by using Ramsey's Theorem:  assign one of $r$ possible colors to
$\{ X,Y\} \subseteq \bfR, X \not = Y$, according to the minimum $\alpha 
\leq r$ such that $X <_L Y$ or $Y <_L X$ holds via $\alpha$.
  For $ k = | \bfR |$ large enough, there
exists a homogeneous set $F$ for this coloring of cardinality 4.  Since 
either $X <_L Y$ or $Y
<_L X$ for each $X,Y \in \bfR, X \not = Y$, there must be $\alpha$ such 
that, for $X ,Y \in F$,

$$
X <_L Y \Leftrightarrow X <_L Y \quad \hbox {via}\quad \alpha .
$$
Now let $F = \{ A,B,C,D\}$, $A <_L B <_L C <_L D$.  We show $C <_L B$, a
contradiction.  $\EE (C) \models \varphi^\alpha_0 (C,\emptyset ,\tilde P \cap C)$
holds since $C <_L D$ via $\alpha$, and $\EE(B) \models \varphi^\alpha_1 
(\emptyset ,B,\tilde P \cap B)$ because $A <_L B$ via $\alpha$.  Finally
$\EE (\overline {B\cup C}) \models \varphi^\alpha_2 (\emptyset ,\emptyset ,\tilde P
\cap (\overline {B\cup C})) $ is true since $B <_L C$.  This shows $C <_L B$
via $\alpha$.  
\eopnospace

\begin {theorem} It is not possible to code an 
infinite linear ordering in $\EE$ without parameters.

\end{theorem}

{\it Proof.}  We write $\tilde X <_L \tilde Y$ for $\tilde X \le_L
\tilde Y\lland  \tilde Y \not \le_L \tilde X$.   Suppose for a contradiction that there is an $\EE$-definable
 $2n$--ary relation $\leq_L$  which is a linear preordering on $\EE^n$ 
such that the
equivalence relation $\tilde X \equiv_L \tilde Y \Leftrightarrow
\tilde X \leq_L \tilde Y \leq_L \tilde X$ has infinitely many equivalence
classes.  We say that a computable set $R$ {\it supports}\id{support}  $A$ if $A \subseteq R$ or
$\overline R \subseteq A$.  $R$ supports $(A_0 ,\dots ,A_{n-1})$ if $R$
supports each set
$A_i$.  Let $$\CCC = \{ R: |R| = |\overline R| = \infty \}.$$

\begin{lemma}  For each tuple $\tilde A = (A_0 ,\dots 
,A_{n-1})$ of sets there exists $R \in \CCC$ such that $R$ supports $\tilde A$.

\end{lemma}

\demo {Proof}  We say that $S$ {\it co-supports} $A$ if $\overline S$ supports
$A$, i.e. $S \subseteq A$ or $A \subseteq \overline S$.  This notion is
closed downwards in $S$.  We define inductively sets $S_k \in \CCC$  co--
supporting $A_0 , \dots ,A_k$.  Then $R = \overline S_{n-1}$ is as required.
 
Let $S_0$ be a set in $\CCC$ which is a subset of $A_0$ if $A_0$ is infinite
and of $\overline A_0$ else.  If $k < n - 1$ and $S_k \cap A_{k+1}$ is 
infinite
let $S_{k+1} \in \CCC$ be a computable subset of $S_k \cap A_{k+1}$.  Else 
let $S_{k+1} = S_k - A_{k+1}$. 
\eop

We now derive an effective bound on $ |\EE^n / \equiv_L|$ (depending on the defining
formula for $\leq_L$).  First we show that each equivalence class of
$\equiv_L$ is large in the following sense:  for each $\tilde A \in \EE^n$,

\begin{equation}
(\forall S \in \CCC) (\exists \tilde B \equiv_L \tilde A) [S \text
{ supports }  \tilde B] \label{tag18}
\end{equation}

Fix $R \in \CCC$ supporting $\tilde A$, and let $S \in \CCC$ be arbitrary.
First suppose that $R \cap S = \emptyset$, and let $\pi$ be a computable permutation of
order 2 which exchanges $R$ and $S$ and is the identity on $\overline {R
\cup S}$.  Let $B_i = \pi(A_i) (i < n)$.  Then $S$ supports $\tilde B$.  Now
$\tilde A \leq_L \tilde B$ is equivalent to $\tilde B = \pi (\tilde A)
\leq_L \pi (\tilde B ) = \tilde A$, since $\leq_L$ is definable.  So $\tilde
A \equiv_L \tilde B$.

\noindent If $R \cap S$ is finite, proceed as above, replacing $S$ by $S - R$.  Then
$\tilde B$ is supported by $S - R$ and hence by $S$.  If $R \cap S$ is 
infinite,
obtain first $\tilde B_0 \equiv_L \tilde A$ supported by $\overline R$ and
then $\tilde B \equiv_L \tilde B_0$  supported by $R \cap S$.  Then $\tilde
B \equiv_L \tilde A$ and $\tilde B$ is supported by $S$.

Suppose $|\EE^n / \equiv_L | \ge p$.  We derive a bound on
$p$.  By (\ref{tag18}), let $S_0, \ldots ,S_{p-1} \in \CCC$ be
pairwise disjoint sets and let $\tilde B^i$, $i < p$, be $n$-tuples of
 sets supported by $S_i$ such that
$$
\tilde B^0 <_L \dots <_L \tilde B^{p-1}.
$$
If a tuple $\tilde X = (X_0, \ldots ,X_{n-1})$ is supported by $S$, we
assign a signature $\beta \in \{ 0,1\}^n$ to $(\tilde X,S)$ by
$\beta (k) = 0 \Leftrightarrow X_k \subseteq S \;\;\; (k < n)$.  Fix
an arbitrary number $q$.
If $p \geq 2^nq$, then there is a subsequence
$(\tilde A^j, R_j)_{j<q}$ of $(\tilde B^i, S_i)_{i<p}$ such that all $(\tilde
A^i, R_i)$ have the same signature $\beta$.  Let
$$
A_k = \bigcup_{j<q} (A^j_k \cap R_j) \qquad (k < n) .
$$
We show that the parameters $A_0, \ldots ,A_{n-1}$ can be used to define in
a first-order way a linear order on $\{ R_0 ,\ldots, R_{q-1}\}$.  Clearly 
one can
decode each $\tilde A^j$ in a uniform first-order way from $R_j$ and this list of
parameters, because $A^j_k = A_k \cap R_j$ if $\beta (k) = 0$ and $A^j_k =
(A_k \cap R_j) \cup \overline R_j$ if $\beta (k) = 1$.  Thus for the formula
$\psi (R,S; A_0, \dots ,A_{n-1})$
expressing $\tilde C <_L \tilde D$, where
 $C_k$ is $A_k \cap R$ if $\beta (k) = 0$ and
$(A_k \cap R) \cup \overline R$ else, and
 $D_k$ is $A_k
\cap S$ if $\beta (k) = 0$ and
$(A_k \cap S) \cup \overline S$ else,

$$
\psi (R_i, R_j; A_0, \dots ,A_{n-1}) \Leftrightarrow \tilde A^i <_L \tilde
A^j ,
$$
so $\psi$ defines a linear order on $\{ R_0 ,\ldots ,R_{q-1}\}$ with the
parameters 

$A_0 , \ldots   , A_{n-1}$.  By Theorem \ref{Rfat}, this gives
an effective bound on $q$ depending on $\psi$ (where $\psi$ was obtained in an effective way from $\phi$ and $\beta$, but did not depend on $q$).  Hence $|\EE^n / \equiv_L |$ cannot
exceed $2^n$ times this bound. Since we can take the maximum over all possible $\beta$, we effectively obtain a bound  which only depends on $\phi$. \eopnospace

For {quasimaximal}\id{quasimaximal!$n(A)$} $A$, let

$$n(A) \ = \ttext{number of atoms in} \LL^*(A).$$
\begin{corollary}  The following relation, which is  arithmetical and
  invariant under automorphisms, is not definable in
$\EE^*$:

  $$\{ \langle A^*,B^*\rangle : n (A) \leq n (B) \}.$$

\end{corollary}

\demo {Proof}  Definability would enable us to code without parameters $(\NN , \leq )$ on
equivalence classes.  \eop

Let $A \approx B$ denote that $A$ is automorphic to $B$ in $\EE$. Soare \citelab{Soare:74}{2} proves that, for quasimaximal $A,B$,
$
A \approx B \Leftrightarrow n(A) = n(B)
$.

Therefore,
$
n(A) \leq n(B) \Leftrightarrow (\exists B') [B \subseteq B' \wedge B'\approx
A]
$.
In fact the automorphism obtained in \cite{Soare:74}  can be represented by a $\Sigma^0_3$ map on
indices.

\begin{corollary} The following relations (which are 
invariant under $=^*$) are non--definable in $\EE$:

\be

\item[(i)]  $ A \approx B$

\item[(ii)]  $A \approx B$ via a $\Sigma^0_3$ automorphism.

\ee

\end{corollary}

\demo {Proof}  Definability of either one of the relations, together with
(i) of  Corollary \ref{qqmax}, would imply the definability of
$$
\{ (A,B): A,B \ttext{quasimaximal} \lland  n(A) \leq n(B) \} ,
$$
so $(\NN ,\leq)$ could be coded in $\EE$ without parameters. 
\eop

 We conclude this Section with  an example of a {\it unary} relation on $\EE^*$ which is  arithmetical and
  invariant under automorphisms: the class
$$ \{A^*: n(A)\ge 2 \ \land \  n(A)  \ \mbox{is a  power of}  \ 2\} $$
is not definable. 

\begin{theorem} \label{qmaxnd} Let $X \sub \NN$ be an infinite  set of even numbers such that
  for each distinct $n,m \in X$,   $(n+m)/_{\ns 2}$ is not in $X$ (for example
let  $X= \{n \ge 2: n\ \mbox{is a  power of} \ 2 \}  $.
Then $\{A^*: n(A) \in X \}$ is not definable in $\EE^*$.
\end{theorem}

 Notice that  for each $X$, $\{A^*: n(A) \in X \}$
is invariant under automorphisms of $\EE^*$. Moreover, if $X$ is
arithmetical, then this class  has an arithmetical index
set. 

\pf  Let $P=\{A: n(A) \in X \}$. Since $P$ is closed under  finite differences, by a result of
Lachlan \citelab{Lachlan:68*3}{4} described in (\ref{LachC}), it suffices to prove nondefinability of $P$ in
$\EE$ (however, one could also perform some notational changes below
to give a direct proof for $\EE^*$). If $A$ is  quasimaximal and $R$ is an infinite coinfinite
computable set, then $A \cap R$ is quasimaximal in $\EE(R) = [ \ES,
R]_\EE$. Let $n_R(A)$ denote $n(A\cap R)$ (evaluated in $\EE(R)$). If
$B^*$ is an atom above $A^*$ in $\LL^*(A)$, then either $B \sub^* A
\cup R$, in which case  $(B \cap R)^*$ is an atom above $(A\cap R)^*$,
or  $B \sub^* A \cup \overline{R}$, in which case $(B \cap \overline{R})^*$ is an
atom above $(A \cap \overline{R})^*$. Conversely, each atom above $(A \cap
R)^*$ gives rise to one above $A^*$, and similarly for atoms above $(A
\cap \overline R )^*$.
Thus 
$$n(A) = n_R(A) + n_{\overline{R}} (A).$$

We use the result of Feferman and Vaught (\ref{Fefi}) for $k=1$.
\n If $R$ is an infinite coinfinite computable set, then
 $\EE \cong \EE \times \EE$ via the map
$$X \mapsto (X \cap R, X \cap 
\overline R). $$

\n Thus, if $P$ is definable in $\EE$ by a formula $\phi(x)$, 
then 

\begin{equation}
  \label{st}
  \EE \models  \phi(X)  \LR  \bigvee_{\alpha = 1,\dots,r}
      (\EE (R) \models \phi^\alpha_0 (X \cap R)
\lland 
       \EE (\overline{R}) \models \phi^\alpha_1 (X\cap \overline{R})).
\end{equation}

\n For each $C \in P$, choose some computable set $R_C$ such
that $n(R_C) = n(C)/2$. By the pigeonhole principle, there are sets $A,B
\in P$, $n(A) \neq n(B) $ so that (\ref{st}) holds via the same
$\alpha$, if $R$ is $R(A)$ ($R(B)$, respectively). After applying an
appropriate computable permutation, we can assume that $R=R_A=R_B$.
Let
$D= (A \cap R) \cup (B \cap \overline{R})$. Then $\EE \models \phi( D)$,
because 
 $$\EE (R) \models \phi^\alpha_0 (D \cap R)
\lland
       \EE (\overline{R}) \models \phi_1^\alpha (D \cap \overline{R}).$$

\n But $n(D)=(n(A)+n(B))/_{\ns 2} \not \in X$, contradiction. \eop

\chapter{Ideal lattices}

\label{CHIB}

We prove that if $\BB$
is an \ird{effectively dense} boolean algebra, then  the theory of
the ideal lattice $\IB$  is undecidable. The next chapter contains  applications of
this result: we present  a coding of a lattice $\IB$ in various
structures, in many cases even without parameters. Thus the theory of those
structures is undecidable. In a forthcoming paper \cite{Nies:eff}\label{gulliwulli}, the author proves
by a much harder argument that $\Th(\IB)$ actually interprets true
arithmetic. 

\section{Computably enumerable boolean algebras}
\label{CEB}
First we  define in detail  the concepts introduced in
Section \ref{EIB}. 
We specify the notion of a c.e.\ boolean algebra\id{boolean algebra!c.e.} as follows. A c.e.\
boolean algebra is represented by a model

\begin{equation} \label{REboolean algebra}(\NN, \preceq, \vee, \wedge) \end{equation}

\n such that $ \preceq $ is a c.e. relation which is a preordering, $
\vee, \wedge$ are total computable binary 
functions, and the quotient structure

\begin{equation} \label{REboolean algebra2}  \BB =(\NN, \preceq, \vee,
  \wedge)/_{\ns  \approx}
\end{equation}

\n is  a boolean algebra (where $n \approx m \LR n \preceq m \lland m \preceq
n$.) We  require that 0 is an index for the least element of $\BB$, and
1 is an index for the greatest element. Then $ 0 \not \approx 1$ by the
definition of boolean algebras. Note that, in an effective way, for each $n$ we can find an index
for a 
complement of $n/_{\ns  \approx}$ in $\BB$, denoted by 

\begin{equation}
\label{cpln}
\mbox{Cpl}(n).
\end{equation}

  At stage $s$ of the algorithm, see if there is $b \le s$
such that $n \wedge b \approx 0$ and $n \vee b \approx 1$, and these
equivalences can be verified in $\le s$ steps. If so, return $b$ as an
output. We write $b-c$ for $b \wedge  \mbox{Cpl}(c)$ and $b \prec c $
if  $b \preceq c \lland c \not \preceq b$.   In general, ``$b \prec c
$'' is not decidable.

We will often relativize our results to $\jump{k-1}$. To define
 the notion of a $\Sigma^0_k$-boolean algebras\id{boolean algebra!$\Sk$},
  one requires that $\pr$ be
 $\Sigma^0_k$ and 
that $\wedge, \vee$ 
be computable in $\ES^{(k-1)}$. 

For a $\Sk$--boolean algebra  $\BB$, let 

\begin{equation} \label{IBBB} \IB \ =  \  \mbox{the lattice of
$\Sk$--ideals of}\id{ideals!$\Sigma^0_k$}\id{ideals!lattice of} \ \BB. \end{equation}
In the following we will mostly use the terminology of c.e.\ boolean
algebras. It will    be clear how to relativize the notions to the
$\Sk$-cases for $k>1$. 

A c.e.\  boolean algebra $\BB$ is called  {\it effectively
  dense}\id{effectively dense!definition of}
if there is a computable $F$ such that $\fa x\ [F(x) \preceq x] $ and 
\begin{equation} \label{eedd} \fa x \not \approx 0 \ \ [0 \prec F(x)
  \prec x]. \end{equation}

More generally,  a $\Sk$-boolean algebra $\BB$ is effectively
dense\id{effectively dense!$\Sk$-boolean algebra} if  (\ref{eedd}) holds with some  $F \le_T \ES^{(k-1)}$.
All effectively dense boolean algebras  are  dense and hence isomorphic,
 but not necessarily effectively isomorphic. Thus our study of boolean algebras is in the
 spirit of recursive model theory, and not along the lines of
Feiner \citelab{Feiner:74}{1}, where (classical) isomorphism types of
c.e.\ boolean algebras are investigated. Feiner proves that there is
a c.e.\ boolean algebra which is not isomorphic to a recursive one.

\begin{example} \label{Rosser}  Let $T$ be a consistent recursively
  axiomatizable theory, and let  $\BB_T$ be 
                    Lindenbaum algebra of sentences over $T$. If $T$
contains Robinson's $Q$, then 
 $ \BB_T$ is effectively dense. \end{example}

\pf We use   Rosser's Theorem
 which asserts that, from an index of a c.e.\ theory $S
\supseteq Q$ one  can effectively obtain a sentence $\alpha$ such that

$$ S \ \mbox{consistent} \ \RA \ S \not \vdash \alpha \ \mbox{and} \
S \not \vdash \neg \alpha.$$

Given $\phi \in \BB_T$, to determine $F(\phi)$ let $S = T \cup
\{\phi\}$. If $ \phi \not \approx 0 $ (in $\BB_T$), then $S$ is
consistent, so $\phi \not \pr \alpha$ and $ \phi \not \pr \neg
\alpha$. Thus let $F(\phi)= \phi \wedge  \alpha_S$. \eop

Notice that, by a result of  Montagna and Sorbi
\citelab{Montagna.Sorbi:85}{1}, the boolean algebras for all such theories are
effectively isomorphic. In that paper the notion of effective
density\id{effectively dense}
(for general c.e.\ lattices) is apparently mentioned   for the first time.  

\section{The theory of ideal lattices}

This section  will be devoted to the following result.

\begin{specialmain} \label{ASS} \label{ITAIB} Suppose $\BB$ is an \ird{effectively dense}
  $\Sk$-boolean algebra. Then 
  $\Th (\IB)$ is hereditarily undecidable. \end{specialmain}

The  main component of the proof is an uniform
definability lemma for the $\SIII$--ideals of $\BB$ which contain a
certain ``separating'' c.e.\ ideal $I_0$, where $|\BB/I_0| =
\infty$. This proof uses some ideas from Section \ref{IDLE}  in the
 context of c.e.\ boolean algebras.

In what follows, for notational reasons we will
actually give codings  in the two sorted structure $(\BB,
\IB)$. This structure can be interpreted in the lattice $\IB$ in a
natural way: represent $b \in \BB$ by the principal ideal $\hat b = [0,b]_\BB$.
 Since the principal ideals are just the complemented elements in 
$\IB$, the set  of ideals in $\IB$ representing elements of $\BB$ 
is definable in $\IB$ without parameters. Moreover, the membership
relation ``$b \in I$'' can be  translated into ``$ \hat b \sub
I$''.

We will find a formula with parameters
 $\psi(x; L,I_0)$ 
such that, as $L$ varies over c.e.\ ideals, $\{x: (\BB, \IB) \models \phi(x;L,I_0)\}$ 
ranges over the $\SIII$--ideals of $\BB$ containing $I_0$. 
 Then, intuitively speaking because $\SIII$ is far
beyond the level of complexity of the c.e.\ structure $\BB$ itself, it
will be possible to give an  interpretation of $\EE^3$ in $\IB$, using
$I_0$ as a parameter.  

\n We say that a c.e.\ ideal $I_0$ of $\BB$ is \emph{separating} 
 if the following holds in $\BB$:
 \begin{equation}
   \label{xysep} \forall x \exists y \preceq x \ y \in I_0 \lland (x \not \in
   I_0 \RA y \not \approx 0)
 \end{equation}

\n and, moreover,   $y$  can be determined
effectively in  $x$. 
The intuition is that a separating ideal nontrivially meets all the 
principal ideals $\neq \{0\}$, in an effective way. 

\begin{lemma} \label{exsep}$\BB$ possesses a separating c.e.\ ideal $I_0$
such that the boolean algebra $\BB/_{\ns I_0}$ is infinite.
\end{lemma}

\pf  We write $b_n$ instead of $n$ if we think of the number $n$ as
determining an element of the boolean algebra under consideration, and
call $b_n$ an \emph{index} for the element $b_n/_\approx$.   First we consider the easier problem how to build  a separating ideal $I_0$
such that $\BB/_{\ns I_0}$ has at least two elements. Recall that $F$ is the
function from (\ref{eedd}). Let $y_0=F(b_0)$ (so $y_0 \not \approx 1$). If $y_0,
\ldots, y_n$ have already been defined, then let $y_{n+1}= y_n \vee F(b_{n+1}
- y_n)$. Let $I_0$ be the ideal generated by $\{y_i: i
\in \NN \}$. Then $I_0$ is c.e.\ and separating, because $b_{n+1} - y_n
\not \approx  0$ if $b_{n+1} \not \in I_0$.  Also $I_0 \neq \BB$: otherwise
suppose that $n$ is the least number such that $y_{n+1} \approx 1$.
Then  $F(b_{n+1} - y_n) \ge \mbox{Cpl}(y_n)$, which is impossible by
our hypothesis on $F$ and since $\mbox{Cpl}(y_n) \not \approx 0$.

We now refine the construction in order to make  
$\BB /_{\ns I_0}$
infinite. To this end, we also define elements $z_0 <  z_1
< \ldots$  of $\BB$ such that $(z_n/_{\ns I_0})_{n \in \NN}$ is a
strictly ascending sequence  in
$\BB/_{\ns I_0}$. As above, let $y_0=F(b_0)$, and let $z_0=0$.
Now,  if $y_0,
\ldots, y_n$ and $z_0 < \ldots < z_n$ have already
been defined, then consider the ``partition''
$$p_0=z_1-z_0, \ldots, p_{n-1} = z_n - z_{n-1}, p_n=
\mbox{Cpl}(z_n).$$

\n Our intention is never to put so much into $I_0$ that one of the
components of the partition goes completely into $I_0$. Let $c_{n+1} =
b_{n+1}- y_n$. Note that, if $c_{n+1} \not \approx 0$, then the same
must hold for  $c_{n+1} \wedge p_i$ for some $i$. Thus if we let

$$y_{n+1} = y_n \vee \bigvee_{i \le n} F(c_{n+1} \wedge p_i),$$
\n we make sure that (\ref{xysep}) is satisfied for $x=b_{n+1}$ via $y_{n+1}$.
\n To make progress on  the ascending sequence, also let
 
\begin{equation}
  \label{z}
  z_{n+1} = F(\mbox{Cpl}(y_{n+1})) \vee z_n.
\end{equation}

\n Again, let $I_0$ be the ideal generated by $\{y_i: i\in \NN \}$.
 We verify that $I_0$ has the required properties.
Since the sequence $(y_n)$ is effective, $I_0$ is c.e.  Moreover,
$I_0$ is separating because, if $b_{n+1} \not \in I_0$, then
$c_{n+1} \not \approx 0$, and therefore $y_{n+1} \not \approx
0$. Furthermore,  $y_{n+1}$ was 
determined effectively from $b_{n+1}$. 
If  $n$ is least such that $y_{n+1} \approx 1$,  then 

$$\bigvee_{i \le n} (F(b_{n+1}
\wedge \mbox{Cpl}(y_n) \wedge p_i)) \succeq \mbox{Cpl}(y_n),$$

\n contrary to the fact that $F(b_{n+1}
\wedge \mbox{Cpl}(y_n) \wedge p_i) < b_{n+1}
\wedge  \mbox{Cpl}(y_n) \wedge  p_i$ for some $i$. 
Thus $y_n \not \approx 1$ for each $n$.

Fix $n$. We now show  that $d= z_{n+1} - z_n \not \in I_0$. Since $
y_{n+1} \not \approx 1$,   $d \not \approx 0$, so $ d \not 
\preceq y_0$. Suppose $k$ is least such that $d \preceq y_{k+1}$. Then $k >n$,
because $d \preceq \mbox{Cpl}(y_{n+1})$, but $y_{k+1} \preceq y_{n+1} $ for $k
\le n$. We now argue as above, but restrict ourselves  to the interval
$[0,d]$. By the minimality of $k$, $d \wedge  y_k \prec  d$, so 

$$ 0 \prec  d - y_k \preceq  \bigvee_{i \le k} F(c_{k+1}
\wedge  p_i).$$

\n Since the $(p_i)_{i \le k}$ form a partition and $d=p_n$, in the
supremum above only the term $F(c_{k+1} \wedge p_n)$ matters. Thus
(recall that $c_{k+1} =
b_{k+1}- y_k$) 
$$d-y_k \preceq F(c_{k+1} \wedge d) = F((d-y_k) \wedge b_{k+1}).$$

\n Since $d-y_k \not \approx  0$,
this contradicts
  the properties of
$F$. \hfill $\square$

\begin{lemma}
Suppose that $I_0$ is a c.e.\ separating ideal. 
For each $\SIII$--ideal $J$, $I_0 \sub J$, there 
is a c.e.\ ideal $L \sub I_0$ such that 
\begin{equation}
  \label{J}
   x \in J \LR \exists r \in I_0 \forall s \in I_0 \ 
        [s \wedge r \approx 0 \RA x \wedge s \in L].
\end{equation}

\end{lemma}

We write $\psi(x,L;I_0)$ for the right hand side in (\ref{J}).

\pf  Choose a  computable function $x \mapsto y(x)$
such that, given $x$, $y(x)$ is a witness for  (\ref{xysep}).
We first define a computable sequence $(s_n)$ which generates $I_0$ as an
ideal and has further useful properties. To start with, since $I_0$ is
c.e.\  there is some 
computable sequence $(y_i)$ generating  $I_0$.  Let $\BB
_{\le e}$ be a finite set of indices for the boolean algebra generated
by $\{ 0, \ldots, e \}$ ($\BB_{ \le e}$ can be obtained effectively
from $e$). Moreover let $s_0=y_0$ and
\begin{equation}
  \label{s}
  s_{n+1}= (y_{n+1}- \hat s_n) \vee \bigvee
\{y(z- \hat s_n): z \in \BB_{ \le n} \},
\end{equation}

\n where $\hat s_n = s_0 \vee \ldots \vee s_n$.
Clearly $ s_i \wedge s_j  \approx 0$ for $i \neq j$.

 Applying Lemma \ref{approxS3} to $P=J$ (viewed as an
index set), we obtain a u.c.e.
sequence $(Z_i)$ such that $Z_i \subseteq \{ 0,\ldots ,i\}$
and

\begin{enumerate}
\item[$\bullet$] $e \in J \RA \mbox{a.e.} \ i \ [e \in Z_i]$
\item[$\bullet$] $\exists^\infty i \ [Z_i \sub J ]$. 
\end{enumerate}

\n Let $L$ be the ideal generated by

$$ \{ e \wedge s_i : e  \in Z_i \}.$$

\n Clearly $L \sub I_0$  and $L$ is c.e.\ We now verify (\ref{J}). 

\smallskip

\n ``$\RA$'' Suppose that $x \in J$. Choose $ \tilde i$ such that
$\forall i > \tilde i \ (x \in Z_i)$ and let $r= s_0 \vee \ldots \vee
s_{\tilde i}$. If $s \in I_0$ and $s \wedge r =0$, then, for some $j >
\tilde i$, 
$s \preceq s_{\tilde i +1} \vee \ldots \vee s_j$. But $x \wedge s_k \in L$
for all $k > \tilde i$. Therefore $x \wedge s \in L$.

\n ``$\LA$'' Now suppose that $ x \not \in J$. Given $r \in I_0$,
choose $k$ such that $r \preceq \hat s_k$.  Choose $ i > k$ such that $Z_i
\sub J$ and also $i >x$. We show that the witness $s_i$ is a
counterexample to (\ref{J}), i.e.\ $x \wedge s_i \not \in L$. 

\n Let $v = x- \bigvee_{e \in Z_i} e - \hat s_{i-1}$. Then $v \not \in
I_0$: else, since $ \hat s_{i-1} \in I_0 \sub J$ and 
$\bigvee_{e \in Z_i} e \in J$, we could infer that $x \in
J$. Therefore $y(v) \not \approx  0$. Also $z=x - \bigvee_{e \in Z_i} e \in
\BB_{\le i-1}$, so $v= z - \hat s_{i-1}$ occurs in the disjunction
(\ref{s}) where  $s_n$ is defined. Hence
$y(v) \preceq s_i \wedge v$ and therefore $ s_i \wedge  x - \bigvee_{e
  \in Z_i} e \not \approx 0$. But this implies that $s_i$ is a counterexample:  if $x \wedge s_i \in L$, then by the fact
that the $(s_k)$ are pairwise disjoint,  $x \wedge s_i \preceq \bigvee_{e
  \in Z_i} e \wedge s_i$. This means that $s_i \wedge (x - \bigvee_{e
  \in Z_i} e) \approx 0$, a contradiction. \hfill $\square$

\begin{lemma} 
$\EE^3 = (\SIII, \sub)$ can be coded  in
$\IB$.
\end{lemma}
\pf By Lemma \ref{exsep}, fix a separating ideal $I_0$ of $\BB$ such that $\BB /_{\ns I_0}$
is infinite. By the previous lemma, the lattice $\mathbf L$ of
$\SIII$--ideals of $\BB$ which contain $I_0$ can be coded in
$(\BB, \IB)$, using $I_0$ as a parameter. We represent a ideal $J$, $I_0 \sub J$ by any $L\sub I_0$ such
that (\ref{J}) is satisfied.

 For completeness' sake we include the
coding scheme in the language of the one-sorted structure $\IB$.  Let
lower case letters  range over principal (i.e.,
complemented) ideals. The scheme is 

\begin{eqnarray*} 
\phi_{dom}(L;I_0) & \equiv & L \le I_0 \\
\phi_\sub(L,H; I_0) & \equiv &  \fa x [\psi(x, L;I_0) \RRA
  \psi(x,H; I_0)] \\
\phi_\equiv(L,H; I_0) & \equiv &
  \phi_\sub(L,H; I_0) \lland \phi_\sub(H,L; I_0).
\end{eqnarray*}

It is now sufficient to show that 
$(\SIII, \sub) \simeq  [C,D]_{\mathbf L}$ for some $C,D \in \mathbf
L$, since $\Th(\SIII, \sub)$ is h.u.\ by Section \ref{CodeRG} and Corollary
\ref{indundecN}. We distinguish two cases. 

\n \emph{Case  A: \   $\BB/_{\ns I_0}$ has infinitely many atoms.} \ Let $C=I_0$
and let $D$ be the ideal generated by $I_0$ and the preimages in $\BB$ of atoms 
of $\BB/_{\ns I_0}$. Notice that ``$x/_{\ns I_0}$ is an atom of $\BB/_{\ns I_0}$" is a $\Pi^0_2$--
property of indices, so there is a function $f \le_T \ES''$ such that
$(f(n)/_{\ns I_0})_{n \in \NN}$ is an enumeration of the atoms  of
$\BB/_{\ns I_0}$ without
repetition. This implies that $D$ is a $\SIII$--ideal 
and, moreover, 
$$ J \mapsto \{n \in \NN: f(n) \in J \}$$

\n is an isomorphism between $[C,D]_{\mathbf L}$ and $(\SIII, \sub)$.

\n \emph{Case  B: \  $\BB/_{\ns I_0}$ has  only finitely many atoms.} \ If   $\BB/_{\ns I_0}$ has only finitely
many atoms, let $b$ be a preimage in $\BB$ of their supremum.
Replacing $I_0$ by the separating ideal $I_0 \vee [0,b]$ if necessary, we can in
fact assume that  $\BB/_{\ns I_0}$ is dense and hence free. Note that  $\BB/_{\ns I_0}$ is c.e.  The standard step--by step
construction of a free generating sequence for a dense countable boolean
algebra produces in the case of  $\BB/_{\ns I_0}$ 
a $\ES'$--sequence $(a_i)$ such that $(a_i/_{\ns I_0})$ is a free generating
sequence for  $\BB/_{\ns I_0}$. Now let $\FF$ be the boolean algebra of
finite or cofinite subsets of $\NN$, and consider the natural map
$g: \BB/_{\ns I_0} \mapsto \FF$ defined by $g(a_i/_{\ns I_0}) = \{i\}$. Clearly $g$
is computable in $\ES'$ if viewed as a map from indices 
for $\BB$ into an effective representation for $\FF$. Let $C$ be the
ideal   $\{x: g(x) =0 \}$ and let $D$ be the ideal generated by the $a_i$'s and
$I_0$. Then $C,D$ are $\SIII$--ideals of $\BB$, contain $I_0$,  and
the 
$\SIII$--ideals $X$ of $\BB$ such that  $C \sub X \sub D$ correspond
to the $\SIII$--ideals of $\FF$ which are contained in the ideal
generated by the atoms. So again $(\SIII, \sub) \simeq  [C,D]_{\mathbf
  L}$. This concludes the proof.
\hfill $\square$

\chapter{Coding Ideal lattices}

\label{CDIL}
\setcounter{theorem}{0} 

Lattices $\IB$ can be coded without parameters in a natural way into three
interesting types of structures.  For the first type, $\BB$ is a c.e.\ boolean algebra, next a
$\SII$, and finally a $\SIII$-boolean algebra. Since $\Th(\IB)$  interprets true
arithmetic (Nies \citelab{Nies:eff}{2}), so do their theories. Here we contend
ourselves with proving undecidability.

\be
 \item Lattices of c.e.\ theories $\{T': T \sub T'\}$ under inclusion,
   where $T$ is a c.e.\ consistent theory containing Robinson 
   arithmetic $Q$

 \item Initial  intervals $[\zer, \a]$ of $\Rec^p_r$, where $\le_r$ is
   a polynomial time reducibility and $\a$ is the $r$-degree of a
   ``super sparse'' set
 \item All intervals of $\EE^*$  which are not boolean algebras.

\ee

Our  first result follows directly from the Main Theorem
\ref{ITAIB}. Let $\LL_T$ be the lattice of c.e.\ extensions of $T$ closed
 under logical inference.

\begin{theorem}  \label{Lin}  $\Th(\LL_T)$ is undecidable.
\end{theorem}
\pf Let $\BB = \BB_T$. In Example \ref{Rosser} it was proved that
$\BB_T$ is \ird{effectively dense}.  Notice that elements of $\LL_T$ are  the c.e.\ filters in $\BB=\BB_T$.
So $\LL_T \cong \IB$ via negation.  \eopnospace

\section{Intervals of $\EE^*$ and $\EE$}
\label{IntE}

In Section \ref{CEINC} we discussed intervals of $\EE$ and $\EE^*$ and gave
several examples. 
A further type
of intervals is obtained by considering the  major subset
relation\id{major subset}. 
\n Maass and Stob \citelab{Maass.Stob:83}{3} proved that  for each  pair $A,B$ such
that   $A\subset_m B$, up to an  (effective)  isomorphism one obtains the same lattices
 $[A,B]_{\EE}$ and $[A^*,B^*]_{\EE^*}$.
These structures  are  denoted by $\mathcal{M}$ and $\mathcal{M^*}$.
From the Maass--Stob result, it follows that $\mathcal{M}^*$ (say) is a distributive lattice with strong homogeneity
properties: all nontrivial closed intervals are isomorphic to the
whole structure, and all nontrivial complemented elements are
automorphic within $\mathcal{M}^*$. However, $\MM^*$ is not a boolean algebra.

A natural question to ask is which  intervals $[A,B]_\EE$  have
an undecidable theory. For instance, Maass and Stob pose this question
for $ \mathcal{M}$, as a part of a programme      to analyze the structure
of $\mathcal{M}$. A complete answer is given by the following result.

\begin{theorem} \label{intervals} Suppose $D \sub A$, where $D,A  \in \EE$. If
  $[D^*,A^*)]_{\EE^*}$ is not a boolean algebra, then 
  $\Th([D^*,A^*]_{\EE^*})$ and in $\Th([D,A]_{\EE})$ are undecidable.
\end{theorem}

Thus $\EE^*$ differs considerably from $\Rm$ and also  the
$\Delta^0_2$- Turing degrees, where intervals of a very different type
with a decidable theory
exist. For instance, in both degree structures there are initial intervals
which form linear orders of order type $\omega+1$ (Lachlan
\citelab{Lachlan:70*1}{lin}  and Lerman \cite{Lerman:83}).

\pf We will reduce the problem to the case of $\MM^*$.   First, we can assume
that $A= \NN$ since each interval of $\EE$ is isomorphic to an end
interval. Moreover we use the
following fact, due to Lachlan.

\begin{fact}
\label{MLachlan}
If $\LL(D)$ is not a boolean algebra, then there exist sets $\wt D, \wt A$ such
that $D \sub \wt D \sub \wt A$ and $\wt D \subset_m \wt A$.
\end{fact}

\pf  Since $\LL(D)=[D, \NN]$ is not a boolean algebra, we can choose
$\wt A \supset D$ such that $\wt A$ is not complemented in $\LL(D)$,
i.e.,  ${\NN-\wt A} \cup  D$ is not c.e.  Pick  $E
\subset_{sm} \wt A$. Then $\wt D:= E \cup D \subset_\infty \wt A$, because,
by the definition of small subsets\id{major subset!small} (\ref{small subset}), 

$$ \wt A =^* E \cup D \RRA \NN \cap (A-E) \sub^*  D \RRA (\NN-\wt
A)\cup D \ \text{c.e.}$$

It  is sufficient to prove the following.

\begin{claim} If $D \subset_m A$, then for some \ird{effectively dense}
  $\SIII$ -boolean algebra $\BB$, $\IB$ can be coded in 
 $ [\tilde D^*, \tilde A^*]$
\end{claim}

(Of course, by \citelab{Maass.Stob:83}{2}, all these intervals are isomorphic. But we
don't make use of the Maass-Stob result, since we directly see that the interpretation is
independent of the particular choices of $D,A$.)
For the case of $\EE^*$, to see the Claim suffices recall  that  $\Th(\IB)$ is h.u.\ by the Main Theorem \ref{ITAIB}, and $\IB$ can
be coded also in  $[D^*,A^*]$ if we use $\tilde D^*, \tilde A^*$ as
parameters. 
By  Fact \ref{indundec}, $\Th([\tilde D^*, \tilde A^*])$ is
undecidable. 

To obtain the result for 
$\Th([\tilde D, \tilde A])$ (and hence for $\Th([D, A])$), note that 
  $\Th([\tilde D^*, \tilde A^*])$ can be interpreted in $\Th([\tilde D, \tilde A])$:
 since $\wt D \subset_m \wt A$,  for $\tilde D \sub P,Q  \sub \tilde A$,
$P =^* Q \LR [P \cap Q, P \cup Q]$ {is a boolean algebra}.

To prove the Claim for $\EE^*$, we will code without parameters the lattice
$\IB$, for the $\SIII$-boolean algebra 

$$\BB=\{X\cup  \wt D^*: X \sq \wt A \} $$

which was already introduced  in (\ref{Splboolean algebra}) (with sets $D \subset_m A$
instead of $\wt D \subset_m \wt A$). Recall that $(U_e
\cup \wt D^*)_{e \in \NN}$ is a $\Delta^0_3$-listing of  $\BB$ (see
Notation \ref{Enot}), via
which $\BB$   becomes a
$\SIII$-boolean algebra. More precisely,  the induced ordering on indices

$$ e \preceq i \LR U_e \sub^* U_i \cup \wt D$$

  is $\SIII$ and
$\ES''$-computable functions $\vee, \wedge$ as in (\ref{REboolean algebra}) can
be defined in the appropriate way. Moreover $\BB$
is
$\ES''$--\ird{effectively dense}, by the Owings
Splitting Theorem (see Soare \cite{Soare:87}): given $e$, the Theorem
provides $V_0,V_1$ such that $\wt D \cup U_e = V_0 \cup V_1$ and 

$$U_e-\wt D \ \text{not co-c.e.} \RRA V_i -D \ \text{not co-c.e.} \
(i=0,1).$$

Let $F(e)$ be such that $U_{F(e)}= V_0$. If $U_e \not \sub^* \wt D$,
then  $\wt D \subset_m \wt D \cup U_e$, so $U_e-\wt D \ \text{not
  co-c.e}$.  Thus, in $\BB$, $0 \prec F(e) \prec e$.
 In fact,
the Owings Splitting Theorem is effective, but
it takes $\ES''$ to determine  $U_k \cup \wt D
^*$ from $k$. 

 By
 Lemma \ref{Fidel}, if  $I$ is a $\SIII$--ideal of $\BB$,
  then there is  $C_I$, $\wt D \sub^* C_I \sub^* \wt A$ such that

$$
   I = \{j: U_{j} \cap  C_I \sub^*\wt D\}.
$$

Conversely, an ideal $I$ satisfying this  for some $C$ must be
$\SIII$.  Now,
for  the desired 
coding of $\IB$, one  represents $\SIII$-ideals $I$ of $\BB$
ambiguously by elements  $c=C_I^*$. Thus, to specify the scheme for
this coding, vacuously  let
$$\phi_{dom}(c) \LR c=c.$$ 
Inclusion of  $\SIII$--ideals  can be defined within
$[\wt D^*,\wt A^*]$ using the formula

\begin{equation}
\phi_\le(c_1,c_2) \equiv \forall x (x \ \mbox{complemented} \
 \RA
(x\wedge c_1= 0
 \RA x \wedge c_2= 0)),
  \end{equation}

\n where $\wt d=\wt D^*$ etc., and $\phi_\eqq(c_1,c_2) \LLR \phi_\le
(c_1,c_2) \lland \phi_\le(c_2,c_1)$. Here, as usual $0$ stands for the
least element in the p.o.\ under consideration, namely $\wt D^*$.   \eop


\section{Complexity theory}
\label{EIBC}

We  proceed to  applications to  
 complexity theory of the method to  code lattices $\IB$ . Since  polynomial time reducibilities are $\SII$ on effectively
presented collections of computable sets, the \ird{effectively dense} boolean algebras we
deal with will be $\SII$. 

\begin{definition}
\label{defssp}
$A$ is {\rm super sparse}\id{super sparse!definition of} via $f$ if  
\be
\item $f$  is a  strictly increasing,
time constructible 
function  $\NN \mapsto \NN$ 

\item $A \sub \{ 0^{f(k)}: k \in \NN \}$ 
and    ``$ 0^{f(k)} \in A$ ?'' 
can be determined in time
$O(f(k+1))$ 
\ee
(Ambos-Spies \citelab{Ambos:86}{EIBC}).
Moreover we require that   

\bi \item[3.] $(\fa r \in \NN) 
 (\mbox{a.e.}\  n)   \ [  f(n)^r < f(n+1)]. $ \ei
A string $w$ is {\rm relevant}\id{relevant}  if $w=0^{f(k)}$ for some $k$.

\end{definition}

Because of the time-constructibility of $f$, we obtain 
\begin{fact}\label{relevant} The set of relevant strings is in $\PP$. \end{fact}

Given a reducibility $\le^p_r$, we denote the degree of a set $X$ by
$\x$ and also write $deg^p_r(X)$ for $\x$.
$\Rec^p_r( \le \a)$ denotes the initial interval of $r$-degrees $\le \a$.

A polynomial time 1-$tt$ reduction of $X$ to $Y$ is a polynomial  time Turing reduction
where in a computation at most one oracle question is asked. 
Thus

\begin{definition} \label{1tt} $X \le^1_{tt}Y$ if there are polynomial time
  computable functions $g: \Sigma^{<\omega} \times \{0,1\} \mapsto
  \{0,1\}$ and $h:\Sigma^{<\omega} \mapsto \Sigma^{<\omega}$ such that

$$ \fa w \in \Sigma^{<\omega} \ [X(w) = g(w, Y(h(x)))].$$
\end{definition}
Polynomial time 1-$tt$ reducibility is a reducibility of more
technical interest. Here is one application of the notion, due to Ambos-Spies.

\begin{theorem}[\citelab{Ambos:86}{1}] \label{Ambos} Suppose $A$ is
  super sparse. Then the polynomial time T-degree
of any set $B \pT A$ consists of  a single 1-$tt$--degree.
\eopnospace \end{theorem} 

Supersparse sets exist in the time classes we are interested in here.

\begin{lemma}[\cite{Ambos:86}] \label{ssplemma} Suppose that $h: \NN
\mapsto \NN$ is
an increasing time constructible function
with $\PTIME \subset \DTIME(h)$, so that $h(n) \ge n+1$ and $h$
eventually dominates all polynomials. Then there is a \ird{super sparse}
computable $A\in \DTIME(h)- \PP$. 
\end{lemma}

{\it Sketch of Proof.} Let $f(n) = h^{(n)}(0)$. Since $h$ eventually
dominates all polynomials, we can construct $A \sub \{0^{f(k)}: k \in
\NN\}$ such that $A \in \DTIME(h)$, but still diagonalize against all
polynomial time machines. \eop

\begin{theorem} \label{ssp} If $A \sub \{0\}^*$ is  super sparse, $A \not \in
  \PTIME$ and $\a= \deg^p_r(A)$, then 

\bc  $\Th(\Rec^p_r( \le \a))$ is undecidable. \ec
\end{theorem}

\pf In a sequence of lemmas, we will code  $\IB$ into $[\zer, \a]$
without parameters,  for
an appropriate \ird{effectively dense} $\SII$--boolean algebra $\BB$. 
We make $\BB$ a very easy, well controlled
part of $[\zer, \a]$, but  use all of $[\zer, \a]$ to sort out $\SII$--ideals of
$\BB$. We begin by introducing  $\BB$.
For an $r$-degree  $\c$, we
let  $\BB(\c)$ be the set   of complemented elements in
$[\zer,\c]$, i.e.

\begin{equation} \label{BB}  \BB(\c) = \{ \x \le \c: \exists \y \ 
   \x \wedge \y = \zer  \ \land \ \x \vee
\y = \c \}.
\end{equation}
A {\it splitting}\id{splitting!in complexity theory}  (or {\it split})
of a set  $B$
is a  set $X$
such that for some $R \in \PP$,
$X=B\cap R$.
We denote this by  $X \sqsubset B$  (via $R$). The advantage of taking a
super sparse $\a$  is that not only is $\BB(\a)$
 a boolean algebra, but in fact it is effectively isomorphic to the boolean algebra
of splittings of $A$, modulo the equivalence relation  under which  two splittings are identified if their
symmetric difference is in $\PP$.  The isomorphism is obtained by
mapping a split $A\cap R$ (represented by an index of the  $\PTIME$ set $R$) to
its degree. In this way, 
 $\BB$ is  well 
controlled as desired. (We could in fact easily ensure that
$A$ has no infinite $\PP$ subsets. In that case $\BB$ is isomorphic to
the boolean algebra of splits modulo finite sets.)


We first show that decomposing a super sparse set $A$ into splits  gives 
complements in the degrees.

\begin{lemma}
\label{l2}
Suppose that  $A$ is super sparse and via $f$ and 
 $A_1= A\cap R, A_2 = A \cap \ol R$ for  $R \in \PP$.
Then  $A_1$ and $A_2$
form a $T$-minimal pair, in the sense that if $Q\leq_T^P A_1,A_2$, 
then $Q\in P$.
\end{lemma}

\pf By Theorem \ref{Ambos}, it is sufficient to prove that 

$$ Q \leo A_1,A_2 \RRA Q \in \PP .$$

Suppose that $Q \leo A_i$ via $g_i,h_i$ ($i=1,2$). The idea to show $Q\in \PP$
is that, if both $h_1(w),h_2(w)$ are relevant oracle queries, then one
of them must be much shorter than the other, so that membership of the
shorter one in the appropriate oracle set can be determined in time polynomial
in the input. The procedure is as follows. Given $w$, compute $h_1(w)$
and $h_2(w)$. If for some $i$ $h_i(w)$ is not relevant,
then $Q(w)=g_i(w,0)$. Else,
\be
\item  if $k =|h_1(w)|=|h_2(w)|$, then see whether $0^k \in
R$. If so, then $Q(w)= g_2(w,0)$, else $Q(w)= g_1(w,0)$.

\item Otherwise, say $|h_1(w)| < |h_2(w)|$. Evaluate $Q(w) =
  h_1(w, A_1(v))$, where $v = h_1(w)$. This is possible in polynomial time, because, by the
  definition of super sparseness, the computation for $A(v)$
 takes  at most $O(|h_2(w)|)$ many steps.  \eopnospace
\ee

Next  we show that, conversely, each pair of complements is represented
by a decomposition into splits. 

\begin{lemma}
\label{l4} Suppose that $\a_1\vee \a_2 =\a$ and  $\a_1\wedge \a_2=\zer$.
Then there exists   a split  $A_1\sqsubset  A$ 
such that $A_1 \in \a_1$ and $A_2=A-A_1 \in \a_2$.
 \end{lemma}

\pf
 It is
sufficient to consider the case that $ r \in \{m, 1\mbox{-}tt\}$.  It is well
known that $\le^p_m$ and $\leo$ induce \ird{distributive} uppersemilattices  on the
computable sets. This is because,  if $X \ppr Y \oplus Z$, then there is
$R \in \PP$ such that $X\cap R \ppr Y$ and $X \cap \overline R \ppr
Z$ (provided that $ r \in \{m, 1\mbox{--}tt\}$). Now, pick sets $B_i \in
\aa_i$ and apply this to $A \ppr B_1 \oplus B_2$ in order to obtain
$R$. It is sufficient to show that in fact $A_1 = A\cap R \equiv^p_r
B_1$ and $A_2 = A\cap \overline R \equiv^p_r
B_2$. Notice that since $B_1 \ppr A_1 \oplus A_2$,
there is $Q \in \PP$ such that  $B_1\cap Q \ppr A_1$ and $B_1 \cap \overline Q \ppr
A_2$. But $B_1, A_2$ form an $r$-- minimal pair, so  $B_1 \cap \overline Q
\in \PP$ and therefore $B_1 \equiv^p_r B_1 \cap Q \ppr A_1$. \eop

 Finally , we show that the order is preserved when passing from
splits modulo $\PP$-subsets of $A$ to degrees.

\begin{lemma} Let $P,Q \in \PP$. Then 
 
$$A \cap P \le^p_r A \cap Q \LR A \cap (P-Q) \in \PP. $$

\end{lemma}

\pf The implication from right to  left is immediate. For the other
implication, notice that $A\cap P$ splits into $A \cap P \cap Q$ and
$A \cap (P-Q)$. But $A \cap (P-Q)$ and $A \cap Q$ form a $T$-minimal pair
by Lemma \ref{l2}. Therefore if $A \cap P \le^p_r A \cap Q$, then
 $A \cap (P-Q) \in \PP$. \eop

Let $(P_e)_{e \in \NN})$ be an effective listing of the polynomial time
sets. 
We  have obtained a representation of $\BB$ in the sense of Section \ref{CEB}: let $e \in \NN$
represent $deg_r(A \cap P_e)$. The  computable functions $\vee,
\wedge$ on $\NN$ are obtained by taking unions and intersections of polynomial
time sets. Clearly,
  ``$A \cap P_e \sub A \cap P_i$''  is $\SII$ in $e,i$.  

\begin{lemma} \label{Bised} $\BB$ is an \ird{effectively dense} $\SII$-boolean algebra. \end{lemma}
\pf By the  uniform diagonalization technique from Landweber e.a.\
\citelab{Landweber:81}{3},
 given a splitting $A
\cap P_e$, we can effectively obtain $Q=P_{F(e)} \sub P_e$  such that 
$A
\cap P_e \not \in \PP$ implies that $A \cap Q, A \cap (P-Q) \not \in
\PP$. For details, see the proof of Theorem 7.3 in 
Balcazar e.a.\ \cite{BDG:88}.
\eop

This concludes our analysis of $\BB$. Next we show how to obtain a
coding  of $\IB$ in $[\zer, \a]$. The
idea is to represent a $\SII$-- ideal  $I$ by a degree $\c_I$ such
that
$$ I = \{\x \in \BB: \x \le \c_I \}.$$
\n Clearly any ideal defined in this way must be $\SII$ (even if $\c_I$
is just the degree of any computable set, not necessarily in
$[\zer, \a]$). The final lemma will show that, 
conversely, each $\SII$ ideal can be represented in that way by a degree
$\c_I \le \a$. Then one obtains the desired parameter free coding of $\IB$
in $[\zer, \a]$,  using the same framework  as  we did for intervals
$[\tilde D^*, \tilde A^*]$ of
$\EE^*$ where $\tilde D \subset_m \tilde A$ in Section \ref{IntE}:   the scheme is given by  the formulas $\phi_{dom}(c)\LR c=c$,

$$
\phi_\le(c_1,c_2) \equiv \forall x \ \mbox{complemented} \
(x\le c_1
 \RA x \le c_2)
$$

 and $\phi_\eqq(c_1,c_2) \equiv \phi_\le
(c_1,c_2) \lland \phi_\le(c_2,c_1)$.

\begin{lemma} 
\label{sortout}
Suppose that $A$  is  super sparse via $f$. Then
  for each $\Sigma^0_2$ ideal $I \tria \BB(\a)$
there is $\c_I \le \a$ such that
$ \fa  \x \in \BB(\a) \ (\x \in I \LR \x \le \c_I). $

\end{lemma}

{\it Proof of Lemma \ref{sortout}}  Recall  that $w$ is {relevant}
if $w=0^k$  for some  $k \in \ \mbox{range}(f)$. We will build $C_I \le^p_m A$ via a $g$ which is computable in
polynomial time. By Theorem \ref{Ambos} it is sufficient to consider
the cases of $\le^p_m$ and $\le^p_{1-tt}$ reducibility. 
Since $I$ is in $\Sigma^0_2$, there is a function $q \le_T \ES'$ such
that $\mbox{range}(q) = \{e: \mbox{deg}^p_r(P_e \cap A) \in I\}$. By
the Limit Lemma in Soare \cite{Soare:87}, there is a computable function $q(e,t)$
such that $q(e) = \lim_t q(e,t)$. 
Let $(h_j)$ be a list of polynomial time $m$-reductions if we consider
$m$-reducibility,
and of polynomial time 1-$tt$ reductions if we consider 1-$tt$-reducibility. 
We meet the coding requirements

$$K_e: A \cap P_{q(e)} \le^p_m C_I$$

by specifying  polynomial time  $m$-reductions to $C_I$. To do so, we
assign $K_e$-{\it coding locations} to  certain relevant $0^s$. If
$s = f(m)$, a $K_e$-coding location for $0^s$ will have the form
$0^n$, $n = {\langle e,r \rangle}$, 
where $r \ge e$ and $f(m) \le n   < f(m+1)$. We will
ensure that $K_e$-coding locations exists for all sufficiently long
relevant $0^s$. We
require that in $n$ steps one  can determine that  $0^s \in P_u$,
where $u$ is the current guess at $q(e)=\lim_t q(e,t)$. We define $C_I$ by specifying
a polynomial time computable $g$ such that $C_I= g^{-1}(A)$, mapping coding
 locations for relevant $0^s$ to $0^s$.
Thus, eventually just the relevant $0^s \in  P_{q(e)}$ are assigned a  $K_e$-coding
location, which  is in  $C_I$ just if  $0^s$ is in $A$.
An
appropriate choice of the $K_e$-coding locations  will ensure that the requirements

$$ H_{\langle i, j \rangle} : A \cap P_i \le^p_r C_I \ \mbox{via} \
[g_j,] \  h_j
\RA A \cap P_i \le^p_r \oplus_{m \le  k} A \cap P_{q(m)} \ (k = \la
i,j \ra )$$

are met. We can suppose that computing $h_j(x)$ takes at most
$p_j(|x|)$ steps, where $$p_j(n)= (n+2)^j .$$  The main idea of the proof
is how to ensure that the coding of $K_e$ does not interfere with the 
requirements $H_i$, $i < e$. We make the length  of any 
$K_e$-coding location for $0^s $ exceed $p_{e-1}(s)$. 

\vsp

{\it The algorithm for $g$.} 

Given an input $x$,  $n =|x|$, first determine in quadratic  time the
maximal  $s \le n$ such that $0^s$ is relevant. This is possible by
the time constructibility of $f$. Now proceed as follows.

\be
\item See if there are  $e,r$ such that 
   $x= 0^{\la e,r \ra}$

\item perform computations $q(e,0), q(e,1), \ldots $ till $n$
steps have passed and let $u$ be  the last value (or $u=0$ if there was no
value so far).

\item  see if   $0^s \in P_u$ in $n$ steps

\item check if  $ p_{e-1}(s) \le n$.

\ee

If (1) and  (3)   are  answered affirmatively and the computation in (4)
stops,
then let $g(x) = 0^s$ (so $x$ is a $K_e$-coding location for $0^s$). Else let  
$g(x)$ be the string $(1) \not \in A$. This completes the algorithm.
Clearly the algorithm takes at most $O(n^2) $ steps.

\vsp
 
Let $C_I = g^{-1}(A)$. We verify that $C_I$ has the required
properties. 

\vspsmall 

{\it Claim 1.} Let $q(e)= \lim_t q(e,t)$. Then $A \cap P_{q(e)}
\le^p_m C_I$.

\pf Let  $p(s)$ be a polynomial which dominates $p_{e-1}(s)$ and  the
number of steps  it takes to
compute $P_{q(e)}$ on the input $0^s$.  Pick an $s_0=f(m)$ such that the
value returned in (2.)  of the algorithm is $q(e)$ for all $s \ge s_0$
and also that, by super sparseness, $\la e , p(f(k))\ra  < f(k+1) $ for all
$k \ge m$. Then for all  $s \ge s_0$, $0^s$ relevant,

$$ 0^s \in A \cap P_{q(e)} \ \LR \ 0^{\la e, p(s)\ra } \in C_I.$$

\vspsmall 

{\it Claim 2.} The requirements $H_{\la i, j \ra}$ are met.

\pf  We first consider the case of $m$-reducibility. Suppose
 that  $A \cap P_i \le^p_r C_I \ \mbox{via} \  h_j$. We
obtain  an $m$-reduction of 
$ A \cap P_i $  to  $\bigoplus_{m < k} A \cap P_{q(m)} \ (k = \la
i,j \ra )$ as follows.   Given a relevant  string  $0^s$,  first
compute  $x=h_j(0^s)$. Since 
$0^s \in A\cap P_i \LLR x \in C_I$, it is sufficient 
 to determine if $x\in C_I$.   Run the algorithm for $g$ on
input $x$. If $g(x) = (1)$ then $x \not \in C_I$. Otherwise $x$ is a
coding location.

{\it Case 1: $|x| < s$.}
 Then
give  $A(g(x))$ as an answer. Since $A$ is super sparse and
$|g(x)| < s$ ,
this answer can be found in time $O(s)$.

{\it Case 2: $n=|x| \ge s$. }

We can suppose that $s \ge s_0$ where $s_0$ is so large that for all
relevant $t \ge s_0$ $|h_j(0^t)|  $ is less than 
the least relevant number bigger than $t$ (by the last  condition in 
Definition \ref{defssp}), and also the computation 
in Step 2 of the algorithm for $g$ with input $0^t$ gives the
final value $ q(e)$ for each $e \le k$.
By the main idea , if $x \in C_I$, then $x$ must be a coding location
for a requirement $K_e$, $ e \le  k$. Since $s \ge s_0$,  $x \in C_I \LR
g(x) \in A \cap P_{q(e)}$.

\vsp

To prove Claim 2 for 1-tt reducibility, suppose that $A \cap P_i \le^p_{1-tt}
C_I$ via $g_j, h_j$. In Case 2, as before obtain an answer $b \in
\{0,1\}$ to ``$x \in C_I$ ?'', depending on a query to the oracle set. Now
give as an output $g_j(x,b)$.
  \eop

\begin{corollary} 
\label{DNC} Suppose that $h$ is time constructible and
  hyperpolynomial.
Then the degrees of (1) all sets  and (2)  of all tally sets in $\DTIME(h)$ have an undecidable theory.
\end{corollary}

\pf  Choose a  super sparse $A \in \DTIME(h) - \PP$ and let $\BB$ be
as before. 
 Because  
$h$ is hyperpolynomial, all the sets
$A \cap P_e$, as well as the sets $C_I$ are in $\DTIME(h)$. By the
preceding result, we obtain a coding of 
 $\IB$ in $(\DTIME(h), \le^p_r)$ with parameter $\a$. Because of
 Fact \ref{inundecN}  and the Main Theorem \ref{ITAIB} this implies that 
$\Th(\DTIME(h), \le^p_r)$ is undecidable. For (2), observe that all sets
involved are tally sets.  \eop

{\it Note}. If $\PTIME=\NP$, then the polynomial time honest degrees 
below any super sparse set form a boolean algebra (Ambos-Spies and
Yang \citelab{Ambos.Yang:90}{1}). So the dishonesty of the
 reduction of $C_I $ to $A$ in the proof of Lemma \ref{sortout}
 appears  to be 
inevitable.

\vsp

One can relativize\id{relativization!in complexity theory} a polynomial time reducibility $\le^p_r$ to a computable oracle
$U$ by replacing the underlying Turing machine model 
by an oracle Turing machine. 
We denote this relativized reducibility by $\le^U_r$.  The
relativization  process is most 
natural
for $\le^p_T$, since

$$ X\le^U_TY \LLR X \oplus U \le^p_T Y \oplus U.$$

Thus, if $\a=\deg^p_T(A)$, then the $\le^U_T$-degrees of the 
computable sets are isomorphic to the end interval  $\{\x \in 
\Rec^p_T: \x \ge \a\}$.

An interesting question arising from Corollary \ref{DNC} is:

\begin{equation} 
\label{PNP} \PTIME \neq \NP \RRA \Th(\NP, \le^p_T) \ \text{is 
undecidable} ?
\end{equation}  

Let $\EXP = \bigcup_{k \in \NN} \DTIME(2^{(n^k)})$. We show that the conclusion holds when relativized to any
computable oracle $U$ such that $\NP^U= \EXP^U$. Such $U$ exist 
by a result of Heller \citelab{Heller:84}{1}.  Clearly $\EXP^U$ is closed downwards
under $\le^U_T$.

\begin{theorem} \label{Dodo}   
$$\NP^U= \EXP^U\RRA \Th(\NP^U, \le^U_T) \ttext{is 
undecidable}.$$
\end{theorem}

\pf To relativize the notion of a super sparse sets to $U$, we change
the second condition in Definition \ref{defssp}: we now require
that 
  ``$ 0^{f(k)} \in A$ ?'' 
can be determined in time
$O(f(k+1))$ with the help of the oracle $U$. All the arguments used in order to prove Theorem \ref{ssp}
are relativizable, including Ambos-Spies' Theorem \ref{Ambos}.  For instance, Lemma \ref{ssplemma} relativized to $U$ states the existence of a $U$-super sparse
$A \nii \PTIME^U$ such that $A$ can be computed in time $h(n)$ with 
 oracle $U$.  We apply this with $h(n)=2^n$.

Notice that the boolean algebra  $\BB$ remains $\SII$ because $U$ is computable.
So we obtain a coding of $\IB$  in the  structure $R^U_A$ of $\le^U_T$-degrees below $A$. (Of course, $R^U_A$ is isomorphic
to the interval $[\u,\a]$ of polynomial time T-degrees, where $\a=\deg^p_T(A \oplus U)$).
Since $\NP^U= \EXP^U$,
$R^U_A$ is an initial interval of the $\le^U_T$-degrees of $\NP^U$-sets. So we obtain a coding of $\IB$ in $(\NP^U, \le^U_T)$. \eop

Next we consider relativizations of the  lattice of $\NP$ sets under
inclusion (which to some extent can be considered as a complexity
theoretic analog of $\EE$). It is not known if $\NP=\coNP$,
i.e. whether this lattice is a boolean algebra. The strongest possible analog to
the question (\ref{PNP}) would thus be: 

\begin{equation}
\label{PNPINC} \NP \neq \coNP \RRA \Th(\NP, \sub) \ttext{is 
undecidable} ?
\end{equation}

One can construct  oracles $X,U$ such that
$\NP^X=\coNP^X$ and $\NP^U\neq \coNP^U$. Here we extend the second oracle result:

\begin{theorem}[\cite{Downey.Nies:97}] There is a computable oracle
  $U$ such that   
$\Th(\NP^U, \sub)  \ \text{is 
undecidable}$. \end{theorem}

\pf  We develop a coding  with parameters of a lattice $\IB$, where $\BB$ is an \ird{effectively dense} 
$\SII$-boolean algebra. 
The proof necessarily produces an oracle $U$ such that
$\NP^U \neq \coNP^U$. In fact we make $\BB$ a boolean algebra which is closely related to

$$\CCC^U := \NP^U \cap \coNP^U,$$

and use the rest of $\NP^U$ to represent $\IB$. A similar idea was used in 
the proof of Theorem \ref{ssp}. {\it Let the variables $R,S$ range over $\CCC^U$.} 
We use the concept of oracle nondeterministic Turing machine (oracle
NTM) which is described in Balcazar e.a.\ \cite{BDG:88}. 

\vspsmall

{\it Outline of the proof.} The construction of $U$ extends Baker e.a.\ \citelab{BGS:75}{1}.
As a parameter, we determine a  set $Q \in \NP^U-\CCC^U$, where for some polynomial
time $S \sub \{0\}^*$, 

\begin{equation}\label{QfromS}
  Q= \{w \in S: \ \ex v\in U \ |v|=|w| \}.
\end{equation}

Then we  let $\BB=\BB(Q)/_{\ns \RR(Q)}$, where

\begin{eqnarray} \label{BRQ}
  \BB(Q) & = & \{Q\cap R: R \in \CCC^U\}, \\
  \RR(Q) & = & \{R \in \CCC^U: R \sub Q \}, \nonumber \\
  \Co\RR(Q)& = & \{Q-R: R \in \RR(Q). \nonumber
\end{eqnarray}

Clearly $\RR(Q)$ is an ideal of $\BB(Q)$. 
With an appropriate numbering of $\NP^U$, $\BB$ is  an \ird{effectively dense}
$\SII$-boolean algebra. 

The general frame for  the coding of $\IB$ follows Lemma \ref{Fidel}, the
$n=1$ case of the ideal definability lemma for $\EE$. However, here  we prefer the language
 of filters. A filter $\FF$ of $\BB(Q)$ is {\it 2-
acceptable}\id{acceptable!in complexity theory}
 if $\Co\RR(Q) \sub \FF$ and $\FF $ has a $\SII$-index set. The construction
 of $U$ will ensure that $\FF$ is 2-acceptable iff  for some $D\sub Q$ in $\NP^U$, 

\begin{equation} \label{FF}
\FF = \{ X \in \BB(Q): \ex R \in \RR(Q) [D-X \sub R] \}.
\end{equation}

Hence the class of 2-acceptable filters is uniformly definable in $\NP^U$. Moreover it
is in 1-1 correspondence with
 the
 class of $\SII$-filters of  $\BB=\BB(Q)/_{\ns \RR(Q)}$, and hence
to $\IB$. In this way we code $\IB$ into $\NP^U$ with a parameter $Q$.

\vspsmall

{\it The  details.} First we need an appropriate listing
of $\CCC^U$. We rely on the fact that $U$, and therefore $Q$, 
 is  given   by a construction which at stage $s$ determines $U^{=s} =
 U \cap \Sigma^s$.

\begin{lemma} \label{CoNP} There is a uniformly computable pair of sequences $(C_e), (\ol C_e)$ such that \be

\item[(i)] for each $e$ we are effectively given   oracle NTMs computing $C_e, 
\ol C_e$ with time bound $(n+2)^{e}$  

\item[(ii)] $C_e \cap \ol C_e =^*\ES$ and $C_e \cup \ol C_e =^*\Sigma^{<
    \omega}$.
\ee
\end{lemma}

\pf Fix some listing  of all oracle NTM $(N_k)$ such that $N_k$ has time
bound $ (n+2)^k$. We write $N_i^U$ for the set accepted by $N_i$ when
the oracle is $U$.
  To determine $C_e$, $e= \la i,j\ra$,  we assume  that $N^U_i$ is the complement of $N^U_j$ until, if ever, this can be refuted in real  time based on oracle queries whose answer has been already determined. Given input $w$, to obtain $C_e(w), \ol C_e(w)$, run $s=|w|$ steps of the following:

\begin{quote}  in lexicographical order, 
for  strings $x$ such that $(|x|+2)^e < s$,
 see whether $x \in N_i^U \LR x \in N_j^U$. 
If so, stop.
\end{quote}
If we stop in $\le s$  steps, then our assumption was wrong, so arbitrarily let $C_e(w)=0, \ol C_e(1)=1$. Else let 
$C_e(w)=N_i^U(w), \ol C_e(w)=N_j^U(w)$.

Clearly (i) and (ii) are satisfied. Moreover, if $N_i^U$ actually is the complement of $N_j^U$, then $C_e= N_i^U$ and $\ol C_e=N_j^U$. \eop

Notice that $\BB(Q)= \{Q \cap C_e: e \in \NN\}$,
so we obtain a presentation in the sense of (\ref{REboolean algebra2})
for  $\BB(Q)$, and hence for  $\BB$. Moreover, $\BB$ with this
presentation is a $\SII$-boolean algebra, because $U$ is computable:

  $$e \pr i \LLR Q\cap (C_e - C_i) \in \RR(Q) \LLR \ex j \ Q\cap
  (C_e-C_i) =  C_j \sub Q,$$

and the matrix  of the last expression is $\Pi^0_1$.

It remains to be proved that $\BB$ is \ird{effectively dense}. This is implied
 by the following relativizable lemma.

\begin{lemma} If $B $ is decidable and $B \nii  \coNP$, then one can in an effective way 
 from a decision procedure for $B$ determine
 a set $R \in \PTIME$ such that $B\cap R, B-R \nii \coNP$. 
\end{lemma}

\pf An easy  application of the delay diagonalization technique,
similar to the proof of Lemma \ref{Bised}.  \eop

Effective density of $\BB$ is  obtained as follows:  given $e$, consider
$B=Q \cap C_e$. Applying the preceding  lemma relativized to $U$ yields
$R \in \PTIME^U$ such that $B \nii \coNP \RRA 
B\cap R, B-R \nii \coNP^U$. Using $\ES'$ as an oracle one can compute
 $i=F(e)$ such that  $B\cap R = Q \cap C_i$.  So $\BB$ is effectively dense
via $F$.

We next describe how to ensure  $Q\not \in \CCC^U$ and introduce a 
first version of the set $S$ needed for (\ref{QfromS}). Using the
technique  of    Baker e.a.\ \citelab{BGS:75}{2},  for
each $e$, we produce a 
witness  $w$ such that $Q(w)=N_e^U(w)$. Thus we meet
the requirements

  $$R_e: Q\neq \Sigma^{<\omega} -N_e^U.$$

If $w$ is our witness  and we see an accepting computation
$N_e^U(w)=1$, we have to put a string $u$ of the same length as $w$
into $U$ which is not an oracle query asked in that computation (or in
accepting computations for requirements which have already been satisfied). 
Let $S = \{0^{s_0}, 0^{s_1},   \ldots \}$, where $s_0=0$ and, for $k > 0$

\begin{equation} \label{ssk}
  s_k= \min\{s> s_{k-1}: s>(s_{k-1}+2)^{k-1}\lland 2^s > G_k(s)\}.
\end{equation}

Here  $G_k(s) =(s+2)^k$, but this definition of $G_k(s)$ will be modified when we
add further requirements. Clearly $S \in \PTIME$ (apply the logarithm
with base 2  to ``$2^s > G_k(s)$''). 

\vspsmall

{\it Construction of $U$, Part 1.} 

For each string $w$, $U(w)=0$ unless otherwise specified. 
\begin{quote} To determine $U^{=s}$ for $s=s_k$, 
check whether there is an $e < k$ such that $R_e$ is not yet met,
namely 

 $$ \fa w \in S [ |w| < s \RRA N_e^U(w)\neq Q(w)].$$

If not, $U^{=s}=\ES$. If so for $e$ minimal, we meet requirement
$R_e$: see whether
$N_e^U(0^s)=1$ via some accepting computation $\Gamma$ based on  the current
oracle. Let $v \in \Sigma^s$ be the lexicographically  first
string which is not an oracle query in $\Gamma$,  and define $U(v)=1$,
thereby causing $Q(0^s)=1$.

\end{quote}

Next we describe how  we obtain, for each 2-acceptable $\FF$ a set
$D\sub Q$ in $\NP^U$ satisfying (\ref{FF}). We identify subsets of
$\BB$ and their preimages  under the canonical map associated with the
presentation (\ref{REboolean algebra2}). Note that there is an  effective
listing
 $(\FF_e)_{e>0}$ of $\SII$-indices for  2-acceptable filters: let $\FF_e$ be the filter generated by $\Co\RR(Q)$ and the $e-1$-th $\SII$-set. (We need $e>0$ for notational reasons.)

Since each $\FF_e$ is infinite (when viewed as a subset of $\NN$), there is a binary function
 $\alpha \le_T \ES'$ such that, for all $e > 0$,  
  $$\FF_e = \{\alpha_e(n): n \in \NN \}.$$

By the Limit Lemma in Soare \cite{Soare:87}, there is a computable
$\beta$ such that, for each $n, e>0$,
$\alpha(e,n) = \lim_k \beta(e,n,k)$. We can assume that

\begin{equation} \label{betasmall} \beta(e,n,k)< k. \end{equation}

To obtain a good representation of $\FF_e$, let

\begin{equation}
\label{Fchain}
F^e_{n,k} = Q\cap \bigcap_{m\le n} C_{\beta(e,m,k)}.
\end{equation}

Then, for each $n$, $F^e_n=\lim_k \ F^e_{n,k}$ exists in the sense that
 an index for an oracle NTM obtained from (\ref{Fchain}) stabilizes. Moreover, the sequence $F^e_0 \supset F^e_1 \supset \ldots$ generates $\FF_e$.

For $e >0$, let

\begin{equation}\label{De} D_e= \{0^{s_k}: e<k \lland \ex w \in U \ |w|=s_k+e\}. \end{equation}

For the inclusion ``$\sub$'' in (\ref{FF}), we ensure that 

\begin{equation} \label{global}  \fa m\ D_e \sub^* F^e_m. \end{equation}
 
Then $X \in \FF_e \RRA \ex m F^e_m \sub X \RRA D_e -X \ \text{finite}.$

For the converse inclusion,
we meet the requirements

$$P_{\la e,m \ra}:  |F^e_m \cap \ol C_m|= \infty \RRA D_e \cap \ol C_m \neq \ES.$$

Then, if $X = Q\cap C_i\nii \FF_e$, we can deduce that $D_e-X \not
\sub R$ for each $R\in \RR(Q)$. Observe that  $X\cup R \nii \FF_e$ because
$\Co\RR(Q) \sub \FF_e$. Choose an $m$ such that $X\cup R =Q\cap C_m$,
and also that  $\ol C_m$ is the complement of $C_m$. Then the hypothesis of $P_{\la e,m \ra} $ is satisfied, thus
$D_e \cap \ol C_m \neq \ES,$ which means that $D_e-X \not \sub R$.

We extend the construction by putting at most one element of length
$s_k+e$, $0<e < k$  into $U$ in order to meet the P-type requirements: according to 
  (\ref{De}) this will determine the sets $D_e$. After presenting  the construction we will determine an appropriate choice of the function $G_k(s)$ needed in (\ref{ssk}).

\vspsmall
{\it Construction of $U$, Part 2.} 

\begin{quote} For $s=s_k$, after determining $U^{=s}$, if we  placed 
 some string of length $s$ into $U$, we also do the following: search for 
 a minimal $\la e,m \ra < k, e> 0$ such that $P_{\la e,m \ra}$ is not
 yet {\it satisfied}, namely 
 \begin{equation}
   \label{before}
  D_e \cap \ol C_m \cap \Sigma^{<s}= \ES, 
 \end{equation}

and also
(based on  the current oracle)
 
\begin{equation} \label{sin}  0^s \in F^e_{m,k} \cap \ol
  C_m. \end{equation} 

 If  $\la e,m \ra < k$ exists, find a $w \in \Sigma^{s+e}$ which does
 not occur as an oracle query in an accepting computation in
 (\ref{sin}), and also not in the  accepting computation $\Gamma$ from
 Part 1, stage $s$. Define $U(w)=1$. We say that $P_{\la e,m \ra}$
{\it receives attention}.
\end{quote}

Now to make sure we can find $w$, we have to count relevant  accepting computations and define $G_k(s)$ appropriately. For a $Q$-type requirement there is at most one, and to determine
$0^s \in F^e_{m,k}$ we need at most $k+1$ many, see (\ref{Fchain}). Notice that these computations have a time bound
 $(s+2)^k$, by the property (\ref{betasmall}). There is one more accepting computation for $0^s \in \ol C_m$. So the definition

  $$G_k(s) = (k+3)(s+2)^k$$

is as desired.

 Clearly $U$ is computable and 
$Q\in \NP^U$. The R-type requirements are met for  the same reasons as
before.
 No requirement is ever injured by a ``later'' $U$-change by the fact that $s_k > (s_{k-1}+2)^{k-1}$ and the construction. 
So by the condition (\ref{before}),  each requirement receives attention at most once. We conclude that (\ref{global}) holds: 
given $e > 0$ and $m$,  choose a $k$ such that for $n <m$,  $\beta(e,n,k)$ has reached its limit and $P_{\la e,n\ra}$ does not receive attention from $s_k$ on. If a requirement  causes $v \in D_e$ at a stage
 $s \ge s_k$, then $s=s_h+e$ for some $h \ge k$ and the requirement is $P_{\la e,n\ra}$ for some $n > m$. Hence $v \in F^e_{n,h} \sub F^e_m$.

To prove that $P_{\la e,m\ra}$ is met, suppose that  $|F^e_m \cap \ol
C_m|= \infty$.  Choose a $k$ such that $\beta(e,m,k)$
has reached its limit and no requirement $P_{u}$, $u < \la e,m \ra$
receives attention at a stage $\ge s_k$. Since $F^e_m\sub Q \sub \{0^{s_i}: i \in \NN\}$,
there is an  $s=s_h \ge s_k$ such that  $0^s \in F^e_m \cap \ol C_m$.  Since $P_{\la
  e,m \ra}$ has the highest priority at $s$, we cause $0^s \in
D_e$. So $P_{\la e,m \ra}$ is met.   \eop

\chapter{C.e. weak truth-table degrees}

\setcounter{theorem}{0}

\label{ChapterRwtt}

We give a coding without parameters of a copy of $\Nops$ in $\Rwtt$. 
This implies that $\ThN$ can be interpreted in $\Th(\Rwtt)$. As a tool
we develop a theory of two sorts  of parameter definable subsets, using 
the distributivity\id{distributive} of $\Rwtt$ in an essential way.  One of them is
 the uniformly definable class of EN-sets (``EN'' stands for end segment). These are 
relatively definable without parameters in an \iid{end segment}, i.e.\ an
upward closed subset $E$ of $\Rwtt$, 
while $E$ is definable from two parameters $\c,\d$. The number $n \in \NN$ 
is represented by (parameters defining) any EN-set of size $n$ (but there may 
also be infinite EN-sets).  Using the combinatorics of $EN$-sets, we  give 
first-order definitions in terms of 
parameters of whether two EN-sets have the same size, and of the operations 
$+ $ and $\times$. For instance, for $+$, we express that an EN-set is the 
disjoint union of two others.

The second type of uniformly definable set, called ID-set (``ID'' stands for ideal) is needed to 
single out the finite EN-sets. We will compare EN-sets to ID-sets, using 
uniformly definable maps between the first and the second. We need to
introduce various schemes.  To understand the formulas related to these
schemes, it is vital to keep in mind the convention in \ref{letter
  convention}: if a scheme $S_X$ is given, then $X, X_0, \ldots$ are
objects coded via $S_X$.

\begin{notation} \label{wttnotation} {\rm
As in Soare \cite[p.\ 49]{Soare:87},  we assume that the use of the
computation $\{e\}^A_s(x)$,
$u(A;e,x,s) \le s$. For $e = \la e_0,e_1\ra$ let

 $$[e](x) \simeq \max_{y \le
  x} \phi_{e_1}(y).$$

 Let $[e]^A(x)$ be $\{e_0\}^A(x)$ if $[e](x)$ and
$\{e_0\}^A(x)$
 are  defined,  and the computation has use $\le  [e](x)$. Otherwise
$[e]^A(x)$ is undefined
. In a similar way define the approximations at stage 
$s$, namely $[e]_s(x)$ and $[e]^A_s(x)$. } \end{notation}

Note that $A \le_{wtt}B \LLR A= [e]^B \ttext{for some} e.$ This
implies that

\begin{equation} 
\label{wttS3}  \{ \la e,i \ra: W_e \le_{wtt} W_i \} \ttext{is} \SIII.
 \end{equation}

\section{Uniformly definable classes in $\Rwtt$}

We prove some  facts which lead to the concepts of EN- and ID-sets. 
Most
of the facts are algebraic. We outline the duality between the two concepts,
as far as the non-symmetric framework of an upper semilattice which may not be 
a lattice allows this. In the following let $\DU$ be a distributive upper semilattice with
least and greatest elements $0,1$.

\begin{lemma} \label{prell}  Suppose that
$b, y_0, \ldots, y_n \in D$. 

\be \item If $b \wedge y_i = 0$ for each $i$, then 
$b \wedge \sup_i y_i = 0$
\item If $b \vee y_i = 1$ for each $i$, then there is $t \in D$ such that
$b \vee t =1$ and $t \le y_i$ for each $i$.
\ee
\end{lemma}
\pf  (i) If $0 < x \le b, \sup_i y_i$, then by distributivity, there is an $i$ 
and $r \in D$ such that $0< r \le x, y_i$. But then $r \le b, y_i$, 
contrary to $b \wedge y_i =0$.

(ii) If $n=0$ let $t = y_0$. Else, since $y_1 \le b \vee y_0$, we can 
choose a $t_1 \le y_0$ and $b_1 \le b$ such that $y_1 = b_1 \vee t_1$. 
Then, 
$1 =b \vee b_1 \vee t_1= b \vee t_1$, so if $n \neq 1$, $y_2  \le b\vee 
t_1$ implies that we can pick $t_2 \le t_1$ and $ b_2 \le b$ such that $y_2=
b_2 \vee t_2$. Continuing in this way we obtain $t=t_n \le y_0, \ldots 
y_n$
such that $b \vee t=1$. \eop

For $d_0, \ldots, d_n \in D$, let

\begin{equation} \label{ENS}
E(d_0, \ldots, d_n) \ = \ \{x \in D: \fa y [ {\fa i \le n} (y \le d_i)
             \RRA y \le x]\}.
\end{equation}

Thus $E(d_0, \ldots, d_n)$ is the set of upper bounds of the ideal
$[0,d_0] \cap  \ldots \cap  [0,d_n]$. Note that $d_i \in E(d_0, \ldots, 
d_n)$ for each $i$. Finite EN-sets $\{p_0, \ldots, p_n\}$ will be sets
which are relatively definable in $E(p_0, \ldots, p_n)$. First we need a 
characterization of the elements in such an end segment.

\begin{lemma} \label{E3} For $d_0, \ldots, d_n \in \DU$,

$$ x \in E(d_0, \ldots, d_n) \LLR x= \inf_{i \le n} (x\vee d_i).$$

\end{lemma}
\pf For the direction from right to left, clearly $x \vee d_i \in E(d_0, 
\ldots, d_n)$ for each $i$. Hence, if the infimum exists, it is also an 
upper bound for the ideal $[0,d_0] \cap  \ldots \cap  [0,d_n]$.

For the other direction, we the  argument is similar to the one used 
in the proof of Lemma \ref{prell} (ii). If $y \le x \vee d_i$ for each $i$, 
then by distributivity we can choose $x_0 \le x$ and $q_0  \le p_0$ such 
that
$y = x_0 \vee q_0$. If $n \ge 1$, choose $x_1 \le x, q_1 \le p_1$ such 
that
$q_0=x_1 \vee q_1$. Continuing in this way we obtain $q_n \le p_n$ 
such that 
$q_{n-1} = x_n \vee q_n$. Moreover $q_n \le p_0, \ldots, p_n$, so $q_n 
\le x$.
Hence $q_{n-1} \le x, \ldots, q_0 \le x $ and finally $y \le x$. \eop

For $x, y \in D$, we write 
$$nd[x,y]$$
if $x < y$ and the interval $[x,y]$ does not embed the 4-element boolean
algebra preserving least and greatest element. Clearly $nd[x,y]$ can be 
expressed
in the language of p.o.
In the next lemma, (i) leads to the definition of EN-sets, and (ii) to the 
definition of $ID$-sets. 

\begin{lemma}
\be
\item[(i)] Let $p_0, \ldots, p_n$ be a finite sequence of elements of $D$ 
such 
that for each $i$,  $nd[p_i,1]$ (in particular, $p_i < 1$) and for $i \neq j$,
$p_i \vee p_j =1$. Then $\{p_i: i \le n \}$ is the set of minimal elements
$x$ in $E= E(p_0, \ldots, p_n)$ such that $nd[x,1]$.

\item[(ii)] Let $(a_i)$ be a finite or infinite sequence of elements of $D$ 
such 
that for each $i$ $nd[0,a_i]$  and for $i \neq j$,
$a_i \wedge a_j =0$. Then $\{a_i\}$ is the set of maximal elements
$x$ in $I$ such that $nd[0,x]$, where $I$ is the ideal of $D$ generated by 
$\{a_i\}$.
\ee
\end{lemma}

\pf (i) It is sufficient to prove that

 $$x \in E \lland nd[x,1] \RRA \ex j \ p_j \le x.$$

Since $x<1$, by Lemma \ref{E3} there is $j$ such that $x \vee p_j < 1$. 
Moreover, by Lemma \ref{prell} there is $t$ such that, for all $i \neq j$,
$t \le  x    \vee p_i$ and $x\vee p_j \vee t=1$. We can suppose that $x \le 
t$. By Lemma \ref{E3},  $x = \inf_{k \le n} x \vee p_k$, so $(x\vee p_j)
\wedge t =x$. By $nd[x,1]$, this implies $t=1$, so $x\vee p_i=1$ for
$i \neq j$ and $x =\inf_{k \le n} x \vee p_k =x \vee p_j$.

(ii) It is sufficient to prove that

 $$x \in I\lland nd[0,x] \RRA \ex j \ x \le a_j.$$

Since $x \in I$, $x \le \sup_{i \le n} a_i$ for some $n$. By distributivity,
$x=\sup_{i \le n} \wt a_i$ for some $\wt a_i \le a_i$ ($i \le n$). Since $0<x$, 
some $\wt a_j$ does not equal $0$. By Lemma \ref{prell}, 
$\wt a_j \wedge \sup_{i\le n, i\neq j} \wt a_i =0$, so $nd[0,x]$ implies
that $\wt a_i =0$ for ${i\le n, i\neq j}$, hence $x= \wt a_j \le a_j$. \eop

In the context of $\Rwtt$, we are able to give first-order definitions 
with
parameters of the set $E$ in (i) of the preceding  Lemma, and also of $I$ in
(ii) if $(\a_i)$ is 
a finite or an infinite u.c.e.\ sequence. We use the following 
theorem of Ambos-Spies, Nies and Shore.

\begin{theorem}[\citelab{ANS:92}{1}] \label{EPTnew} Let $I$ be a $\SIII$-ideal of $\Rwtt$. 
Then 
there exists $\a,\b \in \Rwtt$ such that $I=[\zer,\a] \cap [\zer,\b]$. \eop
\end{theorem}

Degrees $\a,\b$ as above are called an \iid{exact pair} for $I$. Note
that, conversely, 
 each ideal which has an exact pair is  $\SIII$, so that the 
 theorem 
constitutes a uniform definability  result for the class of $\SIII$-ideals.

\begin{lemma} \label{uniformly definableL} 
\be
\item[(i)] Suppose $\{\p_0, \ldots, \p_n\}$ 
is a subset of  $\Rwtt$ such that 

$nd[\p_i, \one]$ for each $i$ and 
$\p_i\vee \p_j = \one$ for $i \neq j$. Then  $\{\p_0, \ldots, \p_n\}$  is 
definable from two parameters $\c,\d$ via a formula $\phi_P(x;c,d)$.

\item[(ii)] Suppose $(\a_i)$ is a finite or infinite u.c.e.\ 
sequence in $\Rwtt$ such that  $nd[\zer, \a_i]$ for each $i$ and 
$\a_i\vee \a_j = \zer$ for $i \neq j$. 
 Then $\{\a_i\}$ is definable from two parameters $\c,\d$ via a formula 
$\phi_A(x; c,d)$.
\ee
\end{lemma}

\pf (i) Observe that $I=[\zer,\p_0] \cap \ldots \cap [\zer, \p_n]$ is a 
$\SIII$-ideal by (\ref{wttS3}), so $I=[\zer,\c] \cap [\zer, \d]$ for some 
$\c,\d$. Thus $E(\p_0, \ldots, \p_n)= E(\c,\d)$ is definable from $\c,\d$ 
via the  formula $\psi(x;c,d) \ = \ \fa y[y \le c,d \RRA y \le x]$. Let  
$\phi_P(x;c,d)$ be the formula expressing that $x$ is a minimal element 
in
$\{z: \psi(z;c,d)\}$ such that $nd[x,1]$.

(ii) Let $I$ be the ideal generated by $\{\a_i\}$. It follows from 
(\ref{wttS3}) that $I$ is $\SIII$. So, once again, $I=[\zer,\c] \cap [\zer, 
\d]$ for some 
$\c,\d$. Let $\phi_A(x;c,d)$ be the formula expressing that $x$ is a 
maximal element $\le c,d$ such that $nd[0,x]$.  \eop

We are now ready to specify the notions of EN-sets and ID-sets by 
appropriate schemes of the same type as in Example \ref{rels}.

\begin{definition} \label{EN} 
\be
\item[(i)] Let $S_P$ the scheme given by the formula  $\phi_P(z;c,d)$ 
and the \id{correctness condition} $\alpha(c,d)$ expressing that whenever $x,y$ satisfy the 
formula and $x\neq y$, then $x \vee y=1$. Subsets of $\Rwtt$ coded
via $S_P$ are called {\rm EN-sets}.

\item[(ii)]  Let $S_Z$ the scheme given by the formula  $\phi_Z(z;c,d)$ 
and the correctness condition $\beta(c,d)$ expressing that whenever $x,y$ satisfy the 
formula and $x\neq y$, then $x \wedge y=0$. Subsets of $\Rwtt$ coded
via $S_Z$ are called {\rm ID-sets}.
\ee
\end{definition}

Notice that subsets of finite EN-sets are EN-sets themselves.

\section{The undecidability of $\Th(\Rwtt)$}

Undecidability of $\Th(\Rwtt)$ was first proved in Ambos-Spies
e.a. \citelab{ANS:92}{2}. We use  the fact that there is an easy way to produce finite
EN-sets in order give a quite elementary new proof. The methods will also be used to obtain 
a coding of a copy of $\Nops$. Along the lines of Theorem \ref{NinRm}
we develop a scheme, also denoted by $S_C$, to code arbitrary
relations between finite EN-sets.

The abundance of EN-sets stems from the fact\footnote{The author
  would like to thank  Klaus  Ambos-Spies for suggesting this.}  that each low $\p \in
\Rwtt$ satisfies $nd[\p,\one]$. 
Thus, 
whenever 
$\p_0, \ldots, \p_n$ are low and $\p_i\vee \p_j = \one$ for $i \neq j$, 
then $\{\p_0, \ldots, \p_n\}$ is an EN-set. For each $n$, such
$wtt$-degrees  $\p_0, \ldots, \p_n$ can be obtained by the method of 
the Sacks splitting theorem (see Soare \cite{Soare:87}). In view of later 
applications, we will prove a more general version of this in Proposition 
\ref{Sacks} below.

\begin{theorem} \label{ASP}  If $\p \in \Rwtt$ is low, then  
$nd[\p,\one]$. \end{theorem}
 
\pf We slightly modify the proof of an extension of the Lachlan Non- 
Diamond Theorem in Ambos-Spies \citelab{Ambos:84}{1}. He proves that, if 
$\a_0, \a_1, \b_0,\b_1$ are c.e.\ Turing degrees such that 
$\a_0 \vee \a_1 = \deg_T(\ES')$ and $\b_0\vee \b_1$ is low, then, for 
some $i \le 1$, $\a_i$ is not $\b_i$-\iid{cappable}. Here $\a$ is $\b$-
cappable if there is a $\c \not \le \b$ such that $\b = \a \wedge \c$.
An inspection of the proof reveals that it can be adapted to
$wtt$-reducibility. (The $T$-reductions  built during the construction
have recursively bounded use 
 anyway, and  the 
proof of Lemma 6 [Lemma 9] goes through. In particular, if the
reduction procedures occurring in requirement $R_e$ are now
$wtt$-reductions
$[e_1]^{B_0}$ and $[e_2]^{B_1}$, then the step counting functions $g$
in the proofs of those lemmas can be computed from $B_0$ [$B_1$] with
recursively bounded use. So the weaker hypothesis $C_0 \not
\le_{wtt}B_0$ [$C_1 \not \le_{wtt}B_1$]
suffices.)

Here we use only the special case of the Theorem that $\b_0=\b_1=\p$.
If $nd[\p, \one]$ fails, then there are $\a_0,\a_1 < \one$ such that
$\a_0 \vee \a_1 = \one$ and $\a_0 \wedge \a_1= \p$. So for both $i=0$ 
and $i=1$, $\a_i$ is $\p$-cappable via $\c_i=\a_{1-i}$. \eop

We now prove the existence of appropriate EN-sets.

\begin{proposition} \label{Sacks} Suppose that $\u_0, \ldots, \u_m < 
\one$.
Then for each $n\ge 0$ there exist low $\v_0, \ldots, \v_n \in \Rwtt$ 
such that $\{\v_0, \ldots, \v_n\}$ is an EN-set and $\u_i \vee \v_j< 
\one$ for each $i\le m, j\le n$.
\end{proposition}

\pf Choose c.e.\ sets $U_i \in \u_i$. We construct c.e.\ sets $V_j$ such
that the statement of the theorem holds with $\v_j=\deg_{wtt}(V_j)$.

To achieve $\v_j \vee \v_{j'}= \one$ ($j' \neq j$) we ensure that
$K=V_j \cup V_{j'}$. For $nd[\v_j, \one]$, we  make each $V_j$ low and
apply Theorem \ref{ASP}.  We meet the standard 
lowness requirements

$$L_{ e,j }:  \ \ex^\infty s \ \{e\}^{V_{j}}(e)[s] \ \text{is defined}
\RRA \{e\}^{V_j}(e) \ \text{converges}.$$

Finally, for $\u_i \vee \v_j< \one$ 
($0 \le i\le m, 0 \le j \le n$) we meet the requirements

$$N_{e,i,j}: \ K \neq [e]^{U_i \oplus V_j},$$

by refraining from changing  $V_j$  till a permanent disagreement occurs. Let 
$(R_k)$ be some priority listing of the $L$-type and $N$-type 
requirements. If $R_k$ is $N_{e,i,j}$ let

$$\len(k,s) = \min\{x: \fa y < x \ K(y)= [e]^{U_i \oplus V_j}(y)[s],$$

and let $r(k,s)= \max\{[e]_s(y): y < \len(k,s)$.

If $R_k$ is a lowness requirement $L_{e,j}$, the restraint  associated  with 
$R_k$ is $$r(k, s)= u(V_{j,s}; e,e,s).$$

{\it Construction.}  At stage $s+1$, if $K_s= K_{s+1}$ do 
nothing.
Else, say $y$ is the unique element in $K_{s+1}-K_s$. Determine the 
minimal
$k$ such that $y < r(k,s)$. If $k$ fails to exist enumerate $y$ into all sets
$V_j$. Else let $j$ be the number such that
$R_k=L_{e,j}$ or $R_k=N_{e,i,j}$ for some $e,i$. Then $V_j$ is the set 
such that enumerating $y$ into
$V_j$ would violate $r(k,s)$. So enumerate $y$ into $V_{j'}$, for each 
$j'\neq j$. This completes the description of the construction.

\vspsmall

Clearly $K=V_j \cup V_{j'}$ for $j \neq j'$. By induction
on $k$ we prove:

\begin{lemma}  Let $k\ge 0$. 
\be
\item[(i)] The requirement $R_k$ is met.
\item[(ii)] $r(k)=\lim_s r(k,s)$ exists and is finite.
\ee
\end{lemma}

\pf Assume the Lemma  holds for all $h<k$.  Choose a stage $s_0$ such 
that
for all $h< k$, $r(h,s_0)$ has reached the limit, and $K$ does not change 
below $\max_{h< k} r(h)$ at any stage $s \ge s_0$. Then at no stage 
$s\ge s_0$
can any number $y < r(k,s)$ enter $V_j$, where $j$ is determined from 
$k$ as in the construction:   $j$ is  the number such that
$R_k=L_{e,j}$ or $R_k=N_{e,i,j}$ for some $e,i$.

If $R_k = L_{e,j}$, then $R_k$ is met, because
if ever $\{e\}^{V_{j}}[s]$ converges for $s\ge s_0$, then this
computation is preserved. Hence also $r(k,s)$ reaches its limit. Now 
suppose that $R_k=N_{e,i,j}$.

 For (i), assume for a contradiction
that $K=[e]^{U_i \oplus V_j}$. Then $$\limsup \len(k,s)=\infty.$$  We obtain
a $wtt$-reduction of $K$ to $U_j$ as follows: given an input $y$, 
compute $s \ge s_0$ such that $\len(k,s) > y$ and $U_i| [e](y) = U_{i,s}| 
[e](y)$.
Then $r(k,t)\ge [e](y)$ for all $t\ge s$, so (by he monotonicity of the 
function $[e]$) $[e]^{U_i\oplus V_j}|y+1$ is protected from changing at 
stages $\ge s$. So $K(y)=  [e]^{U_i\oplus V_j}(y)[s]$. Since $\u_i < \one$, 
we conclude that $N_{e,i,j}$ is met. 

For (ii), let  $x$
 be least such that $K(x) \neq [e]^{U_i\oplus V_j}(x)$. Let $s_1\ge s_0$ be 
least such that, $[e](x)$ is defined, then $K(x)$ and $U_i\oplus V_j|[e](x)$ 
have reached their final values at $s_1$. Then  $\len(k,s) \le x$ from 
$s_1$ on, hence $r(k,s)$ reaches it limit. \eop

Our next goal is to code relations between arbitrary finite EN-sets.

\begin{proposition} \label{SC} There is an object  scheme $S_C$ for
 coding objects of the form $(P_0,P_1,R)$ in $\Rwtt$, where $P_0,P_1$ 
are EN-sets, which has the following property: if $P_0,P_1$ are
finite, then for any $R\sub P_0 \times P_1$, $(P_0,P_1,R)$  can be 
coded.
\end{proposition}

\pf $S_C$ contains  parameters $\c_0,\d_0, \c_1,\d_1$  coding
$P_0,P_1$ and  further parameters for the relation $R$.  Suppose that
$P_0 = \{\p_0, \ldots, \p_n\}$ and $P_1 = \{\q_0, \ldots, \q_m\}$. 
First we  assume that, in addition, 

\begin{equation}
\label{restriction}
\p_i \vee \q_j < \one \ttext{for all} i,j.
\end{equation}

We will reduce the general case to this.

As in the proof of Lemma \ref{uniformly definableL}(ii) there are $\g,\h$ such that 

$$E(\g,\h) = E(\{\p_i \vee \q_j: R\p_i\q_j\}).$$ 

We claim that $$R\p_i\q_j \LLR \ex \z \in E(\g,\h)-\{\one\} [\p_i \vee \q_j\le 
\z].$$

For the direction from 
left to right, simply let $\z = \p_i \vee \q_j$. For the other direction, 
suppose that  the right hand side holds via $\z<\one$. By Lemma 
\ref{E3},
$\z = \inf\{\z \vee \p_r \vee \q_s: R \p_r\q_s\}$. But, if not 
$R\p_i\q_j$, then $\z \vee \p_r \vee \q_s= \one$ for each  pair $\p_r, 
\q_s$ in $R$, since $(i,j) \neq (r,s)$ and therefore $\p_i \vee \q_r =
\one$ or $\p_j \vee \q_s =
\one$. This contradicts $\z< \one$.

Now let 

\begin{eqnarray*}
  \label{phirelfm}
  \wt \phi_{rel}(x,y; c_0,d_0,c_1,d_1,g,h) &\LLR& \phi_P(x;c_0,d_0) \lland 
\phi_P(y;c_1,d_1) \lland \\ 
   &&\ex z <1[x,y\le z \lland z \in E(g,h)].
\end{eqnarray*}

Then in this special case each $R\sub P_0 \times P_1$ can be coded via 
$\wt \phi_{rel}$.

To remove the restriction (\ref{restriction}) we imposed, we interpolate with a third EN-
set.  By Proposition \ref{Sacks}, there is an EN-set $\v_0, \ldots, \v_n$ 
such that, for all $k \le n$, $\p_i \vee \v_k < \one$ and $\q_j\vee \v_k 
< \one $  ($i\le n, j\le m$).  Let $F: P_0 \mapsto \{\v_0, \ldots, \v_n\}$ 
be a bijection. Consider the relation $\wt R$ given by $\wt R \v_k\q_j 
\LLR RF^{-1}(\v_k) \q_j$. Both $F \sub P_0 \times \{\v_0, \ldots, 
\v_n\}$ and $\wt R \sub \{\v_0, \ldots, \v_n\} \times P_1$ can be 
coded by parameters via $\wt \phi_{rel}$. Then $R=F\wt R$ can be coded via the 
following formula (think of $z$ as $F(x)$):

\begin{eqnarray*} \phi_{rel}(x,y; \ol p) &\LLR \ex z& [\wt 
\phi_{rel}(x,z;c_0,d_0,c_2,d_2,g_0,h_0) \lland \\
&&\wt 
\phi_{rel}(z,y;c_2,d_2,c_1,d_1,g_1,h_1)], \end{eqnarray*}

where $c_2,d_2$ are parameters coding the auxiliary EN-set and $\ol p$
consists of all 10 parameters. \eop  

The following result 
 only has   the exact pair theorem \ref{EPTnew}, the technique of
 the Sacks
 splitting theorem and  Theorem
\ref{ASP} as recursion theoretic ingredients. 

\begin{theorem}[\cite{ANS:92}] $\Th(\Rwtt)$ is undecidable. \end{theorem}

\pf By Theorem \ref{Lavrov} the class $\CCC$ of finite directed graphs\id{graph!directed}
has a h.u.\ theory. Using Proposition \ref{SC}, $\CCC$ can be
uniformly coded in $\AAA= \{\Rwtt\}$. Hence by Fact \ref{indundec}, 
$\Th(\Rwtt)$ is undecidable.  \eop

Refining  the proof with the tools from Section \ref{fragments} yields
undecidability of  $\Pi_5-\Th(\Rwtt)$ as a partial order. In Lempp and Nies
\citelab{Lempp.Nies:96}{0}  a coding of finite bipartite graphs\id{graph!bipartite} based
on ID-sets is developed, which even yields
undecidability of  $\Pi_4-\Th(\Rwtt)$. The $\Pi_2$-theory of $\Rwtt$  as a partial
order is  decidable (Ambos e.a.\ \citelab{AFLL:96}{1}).

\section{Coding a copy of $\Nops$}

We use the same framework and similar notation as in the proof of
Theorem \ref{NinRm}.

\begin{theorem} A copy of $\Nops$ can be coded in $\Rwtt$ without
  parameters. \end{theorem}

We will use finite EN-sets to represent numbers. The scheme $S_C$ from 
Proposition
\ref{SC} enables us to express by a first-order condition on parameters
that EN-sets have the same cardinality, and also the arithmetical 
operations. In the end we face the harder problem to single out finite 
EN-sets. (Note that, even if our examples were all finite, there is no 
reason to 
believe that all sets defined via the scheme for EN-sets in Definition 
\ref{EN} are finite.)

We introduce the scheme without parameters to code $\Nops$. It consists of 
formulas $\phi_{num}(\ol x), \phi_=(\ol x, \wt y), \phi_+(\ol x, \ol  y, \ol 
z)$ and $\phi_\times(\ol x, \ol y, \ol z)$, where $\ol w$ stands for a pair
of variables $w_0,w_1$ which represent an exact pair needed to code an 
EN-set. The formula $\phi_{num}(\ol x)$ will be dealt with last, but of course
it implies the correctness condition for $S_P$, since $\ol x$ is thought
of as coding an EN-set. 

\vsp

{\it Equality and the arithmetical operations}

\vspsmall

Let $\phi_\equiv(\ol x, \ol y)$ be a formula expressing

$$ \ex C [C \ \text{is bijection} \ P_{\ol x} \mapsto 
P_{\ol y}],$$

using the scheme $S_C$ from  Proposition \ref{SC}. By that
proposition,  if $P_{\ol \a}$ and  $P_{\ol \e}$
are finite,
then $$|P_{\ol \a}|=|P_{\ol \e}|\LLR \Rwtt \models \phi_\equiv(\ol \a, \ol 
\e).$$

Next let $\phi_+(\ol x, \ol y, \ol z)$ be a formula expressing that
 can be partitioned into two sets of the same size as $P_{\ol x}$ and $P_{\ol y}$: 
 $$\ex \ol u \ex \ol v [\phi_\equiv(\ol x, \ol u) \lland  \phi_\equiv(\ol y, \ol v) 
\lland P_{\ol z}= P_{\ol u} \cup P_{\ol v} \lland P_{\ol u} \cap P_{\ol v} 
=\ES].$$

It can easily be checked that, for finite $P_{\ol \a}, P_{\ol \e}, P_{\ol \c}$

$$|P_{\ol \a}|+|P_{\ol \e}| = |P_{\ol \c}| \LLR \Rwtt \models \phi_+(\ol \a, \ol 
\e, \ol \c).$$

For the direction from left to right  one uses that subsets of $P_{\ol
  \c}$ are again EN-sets.

For  $\phi_\times(\ol x, \ol y, \ol z)$ we express in terms of definable
projection maps that $P_{\ol z}$ 
has the same size as the cartesian product $P_{\ol x} \times P_{\ol y}$. 
Thus $\phi_\times(\ol x, \ol y, \ol z)$ expresses

\begin{eqnarray*} \ex C_1 \ex C_2 
      && C_1: P_{\ol z} \mapsto P_{\ol x} \ttext{onto} \lland  C_2: P_{\ol z} \mapsto P_{\ol y} \ttext{onto} \lland \\
      && \fa a \in P_{\ol x} \ \fa b \in P_{\ol y} \ \ex!q \in P_{\ol z}
               [C_1(q) = a \lland C_2(q)=b].
\end{eqnarray*}

Then, for finite $P_{\ol \a}, P_{\ol \e}, P_{\ol \c}$

$$|P_{\ol \a}||P_{\ol \e}| = |P_{\ol \c}| \LLR \Rwtt \models \phi_\times(\ol \a, \ol 
\e, \ol \c).$$

\vsp

{\it Recognizing finiteness}

\vspsmall

To recognize in a first-order way that an EN-set 
coded by two parameters is finite, the idea   is to compare EN-sets to fragments of a uniformly definable 
subclass 
of the ID-sets. ID-sets are not as easy to construct as EN-sets, but a
more involved construction actually yields a u.c.e.\ {\it infinite} ID-set
$$Z^*= \{\a_i: i \in \NN\}.$$

To specify the uniformly definable  subclass of the class of ID-sets
we will impose  conditions on parameters
$\c,\d$ coding $Z= Z_{\c,\d}$ which are satisfied by $Z^*$ and 
imply that 

\be

\item when $\x$ ranges through degrees $\le \c,\d$, then $|Z \cap 
[\zer, \x]|$ assumes all finite cardinalities 
\item if $|P| = |Z \cap [\zer, \x]|$, $\x \le \c, \d$, then a bijection between the two sets
can be uniformly defined.

\ee

ID-sets $Z$ satisfying the conditions will be called {\it good}. For the 
special good ID-set $Z^*=Z^*_{\c,\d}$, $Z^* \cap [\zer, \x]$ is 
finite for $\x \le \c,\d$.  The formula $\phi_{num}$ implies  about $P$ that 
for each good $Z_{\c,\d}$, a bijection between $P$ and some $|Z \cap 
[\zer, \x]$, $\x \le \c,\d$, exist.

The set $Z^*$ is obtained by referring to  a rather hard theorem in 
Ambos-Spies and Soare \citelab{Ambos.Soare:89}{1}. To ensure property (2.) 
above, one has to make all the degrees $\a_i$ low. An  easier result 
in Lempp and Nies \citelab{Lempp.Nies:96}{1}  could also be used, but  has 
the disadvantage that the actual construction needs to be modified
in order to make the degrees $\a_i$ low.

\begin{specialmainlemma}[\cite{Ambos.Soare:89}] \label{ASo} There 
exists a u.c.e.\ sequence
$(A_i)_{i \in \NN}$  such that each $A_i$ is low, $A_i, A_j$ form  a $T$-
minimal pair for $i \neq j$
 and, where $\a_i = \deg_{wtt}(A_i)$, $nd[\zer, \a_i]$ for each $i$. Thus
$Z^*= \{\a_i\}$ is an ID-set.
\end{specialmainlemma}

\pf Recall that noncomputable c.e.\ set $C$ is \iid{non-bounding} if there 
is no minimal pair $A,B$ such that $A,B \le_T C$. This definition makes 
sense also for $wtt$-reducibility. Clearly,  $C$ is $wtt$-non-bounding
iff  $nd[\zer, \d]$ for each 
 $\d \le \c = \deg_{wtt}(C)$.

 In Ambos-Spies  e.a.\ 
\citelab{ANS:92}{3}, Lemma 6, it is proved that each non-bounding $C$
is also $wtt$-non-bounding. 
From Ambos-Spies and Soare \cite{Ambos.Soare:89} one obtains a u.c.e.\ 
sequence $(A_i)$ such that each $A_i$ is $T$-non-bounding and 
$A_i, A_j$ form  a $T$-minimal pair for $i \neq j$. Since there is a 
uniform construction to produce from a given c.e.\ set $A$ a low set 
$\wt A$ such that $\wt A$ is non-computable if $A$ is \cite{Soare:87}, we can assume 
that each set $A_i$ is low. \eopnospace

\begin{definition} \label{goodness} An ID-set $Z$ defined from parameters 
$\c,\d$ is {\rm good}\id{good}
if 

\be
\item[(i)] $\fa \x \le \c,\d (Z \not \sub [\zer, \x])$
\item[(ii)]  $\fa \x \le \c,\d \ \ex \wt P$

\begin{equation} \label{bijectionP} 
    \{\la \u,\v \ra: \u \le \v \lland \u \in Z\cap [\zer, \x] \lland \v \in 
\wt P \}
\end{equation}
$\ttext{is a bijection between} Z\cap [\zer, \x] \ttext{and} \wt P$.
\ee
\end{definition}

Clearly being good can be expressed by a first-order condition on $\c,\d$.
Moreover, (i) implies that $Z$ is infinite: else $Z \sub [\zer, \sup \
Z] \le \c,\d$. 

We will prove that any u.c.e ID-set $Z$ of low $wtt$-degrees is good, 
when defined from an exact pair for the $\SIII$-ideal generated by $Z$. 
In particular
the set $Z^*=\{\a_i\}$ from the Main Lemma \ref{ASo} is good. Assuming
this fact, we now give a first order condition on parameters expressing 
finiteness of an EN-set $P$.

\begin{lemma} $P$ is finite $\LLR$ $\fa \a,\b[ Z_{\a,\b} \ttext{good} 
\RRA \ex \x\le \a,\b$

$$ \ex \wt P [(\ref{bijectionP}) \ttext{is a bijection} \lland \ex C \ C 
\ttext{is bijection} P \leftrightarrow \wt P]].$$
\end{lemma}

\pf For the direction from left to right, assume that $P$ is finite. Because
good $ID$-sets are  infinite, we can choose $F \sub Z$ such that $|F|=|P|$. 
Let $\x= \sup F$ and choose $\wt P$ satisfying (\ref{bijectionP}). By 
Proposition
\ref{SC}, a bijection $P \leftrightarrow \wt P$ can be coded via $S_C$.
 
For the other direction, let $\a,\b$ be an exact pair coding the set 
$Z^*$ obtained from the Main Lemma \ref{ASo}. If $\x \le \a,\b$, 
then $\x \le \a_0, \ldots, \a_n$ for some $n$. By Lemma \ref{prell}, 
$\a_k\wedge \x=\zer$
for all $k > n$, so $Z \cap [\zer, \x]$ is finite. Thus $P$ is finite. \eop

Finally we prove that any infinite u.c.e.  ID-set $Z$ of low $wtt$-degrees is
good.
 Let $Z$ such a set, coded by an 
exact pair $\a,\b$. By a similar  argument as above, $Z\not \sub 
[\zer,\x]$ for any 
$\x \le \a,\b$. Since all degrees in $Z$ are low, it is now sufficient to 
prove the following.

\begin{lemma} Suppose that $\a_0, \ldots, \a_n$ are  low pairwise incomparable
degrees
in $\Rwtt$. Then there is an EN-set $\v_0, \ldots, \v_n$ such that
$$\a_i \le \v_j \LLR i=j.$$
\end{lemma} 

\pf  Choose c.e.\ sets $A_i \in \a_i$. We construct c.e.\ sets $V_j$ such
that the statement of the theorem holds with $\v_j=\deg_{wtt}(A_j
\oplus V_j)$. Clearly $\a_i \le \v_i$.  To ensure 
$\a_i \not \le \v_j$ for $i \neq j$,
we meet the requirements

$$N_{e,i,j}: \ A_i \neq [e]^{A_j \oplus V_j} \ (i \neq j),$$

by the same strategy as in the proof of Proposition \ref{Sacks}:
refrain from changing  $V_j$  till a permanent disagreement
occurs. We will define  some priority listing  $(R_k)_{k \in \NN}$
 of  all the 
requirements.  If $R_k$ is $N_{e,i,j}$ let

$$\len(k,s) = \min\{x: \fa y < x \ A_i(y)= [e]^{A_j \oplus V_j}(y)[s],$$

and let $r(k,s)= \max\{[e]_s(y): y < \len(k,s)\}$.

To achieve $\v_j \vee \v_{j'}= \one$ ($j' \neq j$) as in Proposition
\ref{Sacks} we ensure that
$K=V_j \cup V_{j'}$. For $nd[\v_j, \one]$, we  make each $A_j \oplus V_j$ low and
apply Theorem \ref{ASP}. Lowness is achieved by  the
side effects of the ``pseudo-
lowness requirements''

$$L_{ e,j }:  \ \ex^\infty s \ \{e\}^{A_j \oplus V_{j}}(e)[s] \ \text{is defined}
\RRA \{e\}^{A_j \oplus V_j}(e) \ \text{converges}.$$

While $L_{e,j}$ may fail to be met, it will produce enough restraint
to ensure $(A_j \oplus V_j)' \equiv_T \ES'$. We use a  standard technique
introduced by Robinson. By the recursion theorem, we can assume that
the sets $V_0, \ldots, V_n$ with
specific enumerations are {\it given} (see comment at the end). Since each set $A_i$ ($i \le
n)$ is low, the following property of $e,j$ and a stage number $\wt s$
can be checked with an oracle $\ES'$: 

\begin{equation} \label{lowcheck} 
\ex s \ge \wt s \ [\{e\}^{A_j \oplus V_{j}}(e)[s]\ttext{is defined via
  an} A_j \ \text{-correct computation}].
\end{equation}

By the Limit Lemma (\cite{Soare:87})  we can fix a computable function
$g(\wt s,e,j,t)$ such that $\lim_t \ g(\wt s,e,j,t)$
exists, has value $0$ or $1$, and the limit is $1$ iff (\ref{lowcheck})
holds. Let $(R_k)$  be some priority listing
 of  all the 
requirements.

\vspsmall

{\it Construction.} At Stage $0$ initialize all the lowness
requirements.

{\it Stage $s+1$.} First determine the restraint $r(k,s)$ for all
$k<s$ such that $R_k$ is a lowness requirement $L_{e,j}$. Let $\wt s <
s$ be greatest such that $R_k$ was initialized at $\wt s$. 
If  $\{e\}^{A_j \oplus V_{j}}(e)[s]$ is undefined, let
$r(k,s)=0$. Else let $u$ be the use of this computation and find the least $t \ge s$ such that either
\bi
\item[(1)] $A_{j,t+1}|u \neq A_{j,t}|u$, or
\item[(2)]  $g(\wt s,e,i,t)=1$.
\ei

Since $\lim_t \ g(\wt s,e,j,t) \ \LR$ (\ref{lowcheck}) holds and the
computation at $s$ seems to provide a witness for (\ref{lowcheck}),
one of the two cases has to apply.     
In Case (1) let $r(k,s)=0$, and in Case (2) $r(k,s)=u$.

\vsp

Now,  if $K_s= K_{s+1}$ terminate  stage $s+1$ here. Else, say $y$ is the unique element in $K_{s+1}-K_s$. Determine the 
minimal
$k$ such that $y < r(k,s)$. If $k$ fails to exist enumerate $y$ into all sets
$V_j$. Else let $j$ be the number such that
$R_k=L_{e,j}$ or $R_k=N_{e,i,j}$ for some $e,i$. Enumerate $y$ into $V_{j'}$, for each 
$j'\neq j$. Initialize all the lowness requirements $R_k'$,
$k'>k$. This completes the description of the construction.

\begin{lemma}  Let $k\ge 0$. 
\be
\item[(i)] If $R_k$ is $N_{e,i,j}$, then the requirement $R_k$ is met.
\item[(ii)] $r(k)=\lim_s r(k,s)$ exists and is finite.
\ee
\end{lemma}

\pf Assume the Lemma  holds for all $h<k$.  Choose a stage $s_0$ such 
that
for all $h< k$, $r(h,s_0)$ has reached the limit, and $K$ does not change 
below $\max_{h< k} r(h)$ at any stage $s \ge s_0$.

If $R_k$ is $N_{e,i,j}$, we can prove (i) and (ii) as in
Proposition \ref{Sacks}. In particular, if $A_i= [e]^{A_j \oplus
  V_j}$, then one can obtain a $wtt$ reduction procedure of $A_i$ to
$A_j$, contrary to the assumption that $\a_i,\a_j$ are incomparable. 

Now suppose that $R_k$ is $L_{e,j}$. We have to show that $\lim_s
r(k,s)$ is finite. Let $\wt s$ be the  greatest
stage    where $R_k$ is
initialized (necessarily $\wt s \le s_0$), and pick $s\ge  s_0$ where
$g(\wt s,e,j, s)$ has reached its limit. If the limit is $0$, then
$r(k,t)=0$ for all $t \ge s$. Else, by $(\ref{lowcheck})$ and the
definition of $g$ there is a least stage $t\ge \wt s$ such that  
$\{e\}^{A_j \oplus V_{j}}(e)[t]$ is defined via
  an $A_j$ -correct computation with use $u$ . Then at stage $t$ we define $r(k,t)=u$. 
 Since $R_k$ is not initialized at
  stages $> \wt s$, the computation $\{e\}^{A_j \oplus V_{j}}(e)[t]$
  is preserved. So $r(k,s)=u$ for all $s \ge t$. 
\eopnospace

\begin{lemma} $A_j \oplus V_j$ is low for each $j \le n$. 
\end{lemma}

\pf Given $e$, we have to determine with a $\ES'$-oracle whether
$\{e\}^{A_j \oplus V_j}(e)$ converges. Let $k$ be such that $R_k$ is $L_{e,j}$. 
Note that, in the preceding argument, we can determine $\wt s$ using a
$\ES'$-oracle. 
Then, by (\ref{lowcheck}),

 $$\lim_t g(\wt s,e,j,t)=0 \RRA\{e\}^{A_j \oplus V_j}(e) \ttext{diverges,}$$

 and by the
argument above,  

$$\lim_t g(\wt s,e,j,t)=1\RRA \{e\}^{A_j \oplus V_j}(e) \ttext{converges.}$$ 

The use of the recursion theorem deserves a comment: 
We are given some c.e.\ sets $V_0, \ldots, V_n$ via a
partial recursive enumeration function $\psi$ which maps $s$ to a
strong index for $V_0 \oplus \ldots \oplus V_n [s]$. From this
the construction produces a similar enumeration $\wt \psi$ for sets
$\wt V_0, \ldots,
\wt V_n$. By the recursion theorem, there must be  $\psi$ such that
$\wt \psi = \psi$, and in particular $V_j = \wt V_j$ for $j \le n$. The function $g$ actually contains an extra
argument, namely an index for $\psi$, and in the discussion above we
assume that the extra  argument is an   index such that $\wt \psi = \psi$.  \eop

\hfill And this, kids, is where the story ends. 

\hfill  Andr\'e Nies, 15 years later. 

\hfill Auckland,  2013.

{\small 
\printindex
\index{h.u.|see{hereditarily undecidable}}
\index{hh-simple|see{hyperhypersimple}}  
\index{split|see{splitting}}  
}
{\small

} 
\end{document}